\documentclass{article}
\usepackage{graphicx} \usepackage{amsmath,amsfonts,amsthm,color}
\newtheorem{defn}{Definition}
\newtheorem{ass}{Assumption}

\usepackage{chngcntr}
\counterwithin{figure}{section}

\newcommand{\xk}{\theta_k}

\newcommand{\xkk}{\theta_{k+1}}

\newcommand{\E}{\mathbb{E}}
\usepackage[T1]{fontenc}
\usepackage[utf8]{inputenc}
\usepackage{caption}
\usepackage{algorithm,algpseudocode}
\usepackage{graphicx,url}
\usepackage{amsmath,amssymb,amsopn}
\usepackage{tabularx}
\usepackage{color}
\usepackage{todonotes}
\usepackage{calc}
\usepackage{caption}
\usepackage{longtable}
\usepackage{multirow}
\usepackage{mathrsfs}
\usepackage{array}
\usepackage{hyperref}
\usepackage{bookmark}
\usepackage{tikz}
\usepackage{booktabs}
\usepackage{pgfplotstable}
\usepackage{pgfplots}
\usepackage[outline]{contour}
\usepackage{calrsfs}
\usepackage{graphicx}
\usepgfplotslibrary{groupplots}
\usetikzlibrary{matrix}
\usepackage{amsmath}
\usepackage{fancyhdr}
\usepackage{makecell}
\usetikzlibrary{automata, external, shapes.geometric, fit}
\usetikzlibrary{arrows.meta, shapes, positioning}
\usepackage{subcaption}
\usetikzlibrary{shapes.geometric}

\usepackage{mdframed}

\def \R{\mathbb{R}}						\def \E {\mathbb{E}}

\DeclareMathAlphabet{\pazocal}{OMS}{zplm}{m}{n}

\algnewcommand\algorithmiconput{\textbf{Constants:}}
\algnewcommand\algorithmicinput{\textbf{Input:}}
\algnewcommand\algorithmicoutput{\textbf{Output:}}
\algnewcommand{\algorithmicgoto}{\textbf{go to}}

\algnewcommand\Constants{\item[\algorithmiconput]}
\algnewcommand\Input{\item[\algorithmicinput]}\algnewcommand\Output{\item[\algorithmicoutput]}\algnewcommand{\Goto}[1]{\algorithmicgoto~\ref{#1}}

\usepackage{amsthm}

\newtheorem{theorem}{Theorem}

\newtheoremstyle{nodot}{}{}{}{}{\bfseries}{:}{ }{\thmname{#1}\thmnumber{ #2}\thmnote{ #3}}

\theoremstyle{nodot}
\newtheorem{exampleinner}{Numerical Example}

\makeatletter
\def\thmnote#1{{\bfseries\itshape #1}}
\makeatother

\newenvironment{num_example}[1][]
{\ifthenelse{\equal{#1}{}}{\begin{exampleinner}}{\begin{exampleinner}[#1]}}
{\end{exampleinner}
}

\definecolor{myblack}{RGB}{53, 53, 53}
\definecolor{myblue}{RGB}{40, 75, 99}
\definecolor{myred}{RGB}{192, 50, 33}
\definecolor{myyellow}{RGB}{255, 166, 48}
\definecolor{mywhite}{RGB}{240, 237, 238}
\definecolor{mygreen}{RGB}{0, 102, 0}

\definecolor{green1}{RGB}{9, 82, 86}
\definecolor{green2}{RGB}{8, 127, 140}
\definecolor{green3}{RGB}{6, 167, 125}
\definecolor{green4}{RGB}{79, 109, 122}
\definecolor{green5}{RGB}{192, 214, 223}
\definecolor{violet}{RGB}{26,69,131}

\definecolor{checkgreen}{rgb}{0,0.6,0}
\definecolor{phase1}{rgb}{0.008,0.655,1.000}
\definecolor{phase2}{rgb}{0.016,0.75,0.700}
\definecolor{phase3}{rgb}{0.929,0.35,0.700}
\definecolor{icsyellow}{cmyk}{0.00,0.11,0.53,0.00}

\definecolor{blackmy}{RGB}{38, 70, 83}
\definecolor{bluemy}{RGB}{39, 125, 161}
\definecolor{greenmy}{RGB}{42, 167, 143}
\definecolor{yellowmy}{RGB}{233, 196, 106}
\definecolor{brownmy}{RGB}{244, 162, 97}
\definecolor{redmy}{RGB}{249, 65, 68}

\definecolor{darkbluemy}{RGB}{65, 59, 147}
\definecolor{lightbluemy}{RGB}{71, 139, 194}
\definecolor{greenmy}{RGB}{98, 173, 153}
\definecolor{darkorangemy}{RGB}{230, 142, 52}
\definecolor{lightorangemy}{RGB}{217, 172, 59}

\definecolor{blue1}{RGB}{1, 42, 74}
\definecolor{blue2}{RGB}{1, 73, 124}
\definecolor{blue3}{RGB}{42, 111, 151}
\definecolor{blue4}{RGB}{44, 125, 160}
\definecolor{blue5}{RGB}{70, 143, 175}
\definecolor{blue6}{RGB}{137, 194, 217}

\definecolor{red1}{RGB}{204, 68, 75}
\definecolor{red2}{RGB}{218, 85, 82}
\definecolor{red3}{RGB}{227, 150, 149}
\definecolor{red4}{RGB}{228, 190, 171}

\definecolor{brown1}{RGB}{92,178,112}
\definecolor{brown2}{RGB}{130,194,110}		
\definecolor{brown3}{RGB}{163, 193, 173}

\usepackage[normalem]{ulem}
\usepackage{soul,xcolor}
\usepackage[many]{tcolorbox}
\usepackage{mathtools}
\usepackage{algorithm}
\usepackage{enumitem}
\usepackage{listings}
\usepackage{relsize}
\usepackage{nameref}

\definecolor{takeawaybg}{RGB}{235, 243, 250}
\definecolor{takeawayframe}{RGB}{30, 75, 130}
\definecolor{takeawaytitle}{RGB}{30, 75, 130}

\newtcolorbox{takeawaybox@inner}{
  breakable,
  colback=takeawaybg,
  colframe=takeawayframe,
  colbacktitle=takeawaytitle,
  coltitle=white,
  fonttitle=\bfseries\sffamily,
  title={Key Takeaways},
  boxrule=0.6pt,
  arc=3pt,
  left=12pt, right=12pt, top=8pt, bottom=10pt,
  toptitle=4pt, bottomtitle=4pt, lefttitle=10pt,
  before skip=14pt, after skip=14pt,
}
\newenvironment{takeawaybox}
  {\tikzexternaldisable\begin{takeawaybox@inner}}
  {\end{takeawaybox@inner}\tikzexternalenable}

\definecolor{codebg}{RGB}{248, 248, 250}
\definecolor{codeframe}{RGB}{205, 210, 220}
\definecolor{codekw}{RGB}{0, 80, 160}
\definecolor{codestr}{RGB}{170, 50, 100}
\definecolor{codecom}{RGB}{110, 110, 110}
\definecolor{codenum}{RGB}{160, 160, 160}

\lstdefinestyle{snippetstyle}{
  language=Python,
  basicstyle=\ttfamily\footnotesize,
  keywordstyle=\bfseries\color{codekw},
  commentstyle=\itshape\color{codecom},
  stringstyle=\color{codestr},
  numbers=left,
  numberstyle=\tiny\color{codenum},
  numbersep=8pt,
  showstringspaces=false,
  breaklines=true,
  breakatwhitespace=true,
  tabsize=2,
  columns=fullflexible,
  keepspaces=true,
  upquote=true,
  frame=single,
  rulecolor=\color{codeframe},
  backgroundcolor=\color{codebg},
  xleftmargin=12pt,
  xrightmargin=2pt,
  framexleftmargin=2pt,
  framexrightmargin=2pt,
  framextopmargin=4pt,
  framexbottommargin=4pt,
  abovecaptionskip=6pt,
  belowcaptionskip=2pt,
  aboveskip=12pt,
  belowskip=12pt,
  captionpos=b,
}
\usepackage{longtable}

\newcommand{\algorithmiccommentMine}[1]{\bgroup\hfill$\triangleright$~{ \textcolor{gray}{#1}}\egroup}
\newcommand\COMMENTmine[1]{\algorithmiccommentMine{#1}}

\newcommand\cancel[1]{}

\newcommand\oldtextt[1]{}
\newcommand\oldtext[1]{}
\newcommand\newtext[1]{#1}

\usetikzlibrary{spy}
\usepackage{pgfplots}
\usepgfplotslibrary{fillbetween}

\usepackage{array}
\newcommand{\PreserveBackslash}[1]{\let\temp=\\#1\let\\=\temp}
\newcolumntype{C}[1]{>{\PreserveBackslash\centering}p{#1}}
\newcolumntype{R}[1]{>{\PreserveBackslash\raggedleft}p{#1}}
\newcolumntype{L}[1]{>{\PreserveBackslash\raggedright}p{#1}}

\DeclareMathAlphabet\bpazocal{OMS}{cmsy}{b}{n}
\usetikzlibrary{patterns}

 \usetikzlibrary{arrows.meta}

\definecolor{c1}{HTML}{003F5C}
\definecolor{c2}{HTML}{444E86}
\definecolor{c3}{HTML}{955196}
\definecolor{c4}{HTML}{DD5182}
\definecolor{c5}{HTML}{FF6E54}
\definecolor{c6}{HTML}{FFA600}
\definecolor{cspatial}{HTML}{1F8D60}
\definecolor{call}{HTML}{A6231D}
\definecolor{cB4}{HTML}{0B6FA4}
\definecolor{cB8}{HTML}{4078A1}
\definecolor{cB16}{HTML}{1F8D60}
\definecolor{cB32}{HTML}{D27B11}
\definecolor{cB64}{HTML}{A6231D}
\definecolor{cB128}{HTML}{9D4E9F}
\definecolor{cB256}{HTML}{000000}

\title{Introduction to optimization methods for training SciML models}
\author{Alena Kopani\v{c}\'akov\'a\thanks{Toulouse-INP, IRIT-APO, ANITI, 2 Rue Charles Camichel, 31000 Toulouse, France; email: alena.kopanicakova@toulouse-inp.fr} and Elisa Riccietti\thanks{ENS de Lyon, CNRS, Inria, Universit\`e Claude Bernard Lyon 1, LIP, UMR 5668, 69342, Lyon cedex 07, France; email: elisa.riccietti@ens-lyon.fr}}

\begin{document}

\maketitle

\section{Introduction}
\newtext{Optimization is a key ingredient of modern machine learning (ML). 
Historically, optimization problems have been classified  as \emph{convex} or \emph{non-convex}, according to the \emph{properties} of the objective function and of the feasible set. 
Convex problems are generally more tractable, since global optimality guarantees can be established and any stationary point is also a global minimum.
Non-convex problems may possess multiple local minima and saddle points, making optimization significantly more challenging.}

\newtext{Optimization methods have been categorized by the degree to which they exploit derivative information.}
\emph{First-order} methods rely exclusively on gradient evaluations, whereas (approximate) \emph{second-order} methods incorporate curvature information via Hessian matrices or suitable approximations.
However, the massive scale of modern ML problems, combined with the principles of empirical risk minimization, has driven the field toward \emph{stochastic optimization}.
Stochastic optimization methods, such as Stochastic Gradient Descent (SGD)~\cite{robbins1951stochastic} and its variants (e.g.,~AdaGrad~\cite{duchi2011adaptive}, and Adam~\cite{kingma2017adammethodstochasticoptimization}), employ inexpensive noisy gradient evaluations that scale efficiently with the  data size.
Consequently, much of modern ML focuses on integrating first-order schemes, adaptive gradient methods, and curvature-aware techniques into stochastic optimization frameworks.

This familiar optimization landscape changes substantially in the context of \emph{scientific machine learning} (SciML).
Unlike classical ML, where abundant data allows stochastic approximation to perform well, SciML often operates in data-scarce regimes in which physical models must supplement, or even dominate, the available data.
\newtext{As a result, the optimization problems arising in SciML frequently take the form of
\emph{physics-informed} formulations, in which the objective function incorporates a partial differential equation (PDE) and boundary or initial conditions (BC or IC).
The PDE and BC/ICs are usually imposed in a \emph{soft way} through a penalization term.
An alternative is to use \emph{operator-constrained} formulations, in which the physical constraints are embedded directly into the model, typically through the neural network architecture,   and are therefore enforced in a \emph{hard way}.}

\newtext{
This changes the structure of the objective function compared to classical ML, where the loss decomposes into independent sample contributions, enabling efficient stochastic optimization. 
In SciML, the differential operator induces global spatio-temporal coupling, creating dependencies across the entire domain so that errors at different collocation points are not independent. 
This makes stochastic optimization significantly more challenging and highly sensitive to the choice of the samples.}

\newtext{Moreover, the resulting optimization landscape is typically shaped by the spectra of the underlying differential operators rather than by the statistics of a dataset. As a consequence, it is often highly \emph{anisotropic}, meaning that the curvature varies significantly across different directions, and \emph{stiff}, in the sense that the problem involves widely different scales and Hessian eigenvalues spanning several orders of magnitude. Such loss functions may therefore exhibit directions of very large curvature, for instance those induced by high-order derivatives, alongside nearly flat directions associated with weakly identifiable components of the SciML model.

This structure limits the effectiveness of first-order  methods, since it increases sensitivity to the choice of step size, which must remain small to ensure stability along directions of large curvature, thereby slowing progress in flatter directions. It also limits the effectiveness of  stochastic methods and often necessitates the use of deterministic or large-batch optimization algorithms, as well as curvature-exploiting approaches such as deterministic quasi-Newton methods.}

At the same time, SciML also includes a rapidly growing class of \emph{data-driven} formulations that more closely resemble classical supervised learning.
Prominent examples include operator-learning models~\cite{lu2021learning}, neural surrogate models trained on high-fidelity PDE simulation data~\cite{mishra2018machine}, autoencoder-based reduced-order models that learn low-dimensional latent dynamics~\cite{lazzara2022surrogate}, and neural differential equation frameworks that fit observational or simulated trajectories~\cite{chen2018neural}.
Because their training objectives decompose over samples and admit mini-batching, these models inherit much of the optimization behavior of standard ML and benefit directly from well-established stochastic optimization techniques.
As a result, effective optimization methods in SciML span two complementary regimes: physics-informed settings and data-driven settings that align more closely with large-scale stochastic ML training.

In this chapter, we first introduce the optimization problems that arise in both traditional ML and SciML, highlighting their structural differences and analytical properties.
We then present a foundational overview of deterministic and stochastic optimization methods, explaining how these algorithms can be adapted to address the specific challenges posed by SciML models.
Given the breadth of the optimization field, the material covered here does not aim to be exhaustive.
Instead, we present a concise introduction to common optimization techniques, demonstrate how they can be adapted to SciML training through practical examples\footnote{The accompanying notebooks are available at \\ \href{https://github.com/kopanicakova/intro_opt_SciML.git}{https://github.com/kopanicakova/intro\_opt\_SciML.git.}}, and point toward promising research directions.

\newtext{
\paragraph{Chapter's structure}
The rest of the chapter is organized as follows:
The first three sections form the foundations.
In particular, Section~\ref{sec_ERM} introduces the principles of empirical risk minimization and the finite-sum minimization problem, with examples spanning classical ML and SciML.
Section~\ref{sec:opt_problem} discusses the difficulties associated with solving arising minimization problems and introduces the neural tangent kernel (NTK), an essential tool for analyzing the training dynamics.
Section~\ref{sec:first_order} presents stochastic gradient descent (SGD), the most widely used method for ML optimization, together with a brief convergence analysis and numerical examples.}

\newtext{The following sections dive deeper into more advanced optimization methods.
More precisely, Section~\ref{sec:adaptive_methods} covers accelerated and adaptive stochastic gradient methods (momentum, AdaGrad, Adam).
Section~\ref{sec:second_order} introduces (approximate) second-order methods, including Hessian-free Newton-Krylov and quasi-Newton schemes, such as L-BFGS.}

\begin{figure}[h!]
\resizebox{\linewidth}{!}{
\includegraphics{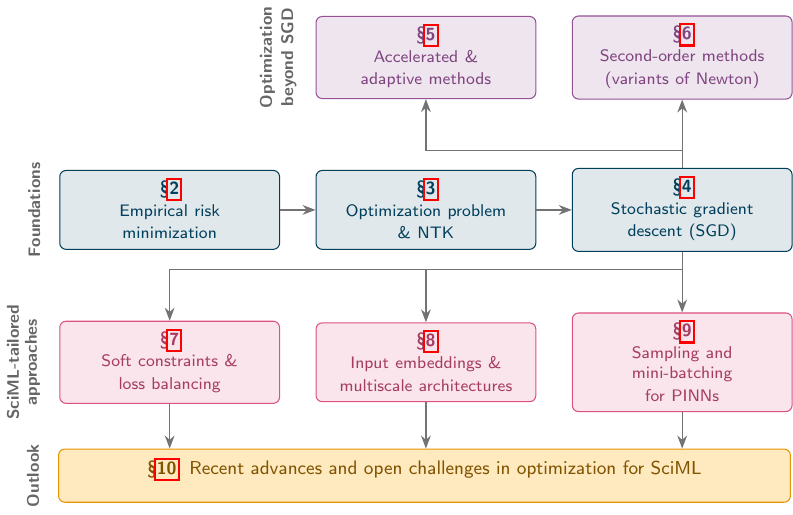}
}
\caption{\newtext{Roadmap of the chapter. Sections~\ref{sec_ERM}--\ref{sec:first_order} form foundations suitable for a first reading. 
Sections~\ref{sec:adaptive_methods}--\ref{sec:second_order} extend the foundations with general-purpose accelerated, adaptive, and second-order methods. 
Sections~\ref{sec:problem_def}--\ref{sec:collocation_minibatch} discuss aspects specific to the training of SciML/PINN models, namely balancing terms related to soft constraints, input embeddings, multiscale architectures, and choice of sampling and mini-batching strategies. 
Section~\ref{sec:outlook} surveys recent research directions and open challenges.}}
\label{fig:roadmap}
\end{figure}

\newtext{While Sections~\ref{sec_ERM}--\ref{sec:second_order} focus on standard optimization methods for ML, Sections~\ref{sec:problem_def}--\ref{sec:collocation_minibatch}, present more advanced concepts, which are specific to the SciML context.
In particular, Section~\ref{sec:problem_def} discusses the impact of regularization and soft-constraints on the performance of the optimization algorithms. 
Section~\ref{sec:data_precond} presents the strategies to mitigate \emph{spectral bias} through input embeddings and multiscale architectures.
Section~\ref{sec:collocation_minibatch} then focuses on sampling strategies that are critical for the training of PINNs.
Finally, Section~\ref{sec:outlook} gives an overview of current research directions and open challenges in optimization for SciML/PINNs.
Notably, the techniques presented in Sections~\ref{sec:problem_def}--\ref{sec:collocation_minibatch} can be applied in conjunction with any optimization method presented in Sections~\ref{sec:first_order}--\ref{sec:second_order} and can be therefore followed right after the reader becomes familiar with (S)GD algorithm. 
A visual roadmap of the chapter is given in Figure~\ref{fig:roadmap}.}

\newtext{\paragraph{Notations}
Throughout the chapter, $\|\cdot\|$ denotes the Euclidean norm, unless otherwise specified. }
  \section{Empirical risk minimization and the finite sum minimization problem}\label{sec_ERM}
In \newtext{ML}, a model is characterized by a parameter vector 
$\theta \in \mathbb{R}^n$ and a parameterized mapping  
$h_\theta: \pazocal{X} \to \pazocal{Y}$,
where $\pazocal{X}$ denotes the input (feature) space and $\pazocal{Y}$ is the output space.
The mapping $h_\theta$ transforms an input sample $x \in \pazocal{X}$ into an output prediction 
$h_{\theta}(x) \in \pazocal{Y}$, and its form is typically defined by a deep neural network (DNN) whose parameters are collected in the vector~$\theta$.

We aim to find a model $h_\theta$ that provides the best possible approximation of the desired mapping. 
The process of searching for such a model is called \emph{training} and \newtext{requires} solving an optimization problem defined by an expected loss over the data distribution.
Let $(x, y)$ denote a generic data pair drawn from an unknown probability distribution 
$\pazocal{P}$ on $\pazocal{X} \times \pazocal{Y}$. 
Let $\ell : \pazocal{Y} \times \pazocal{Y} \to \mathbb{R}_+$ denote a loss function measuring the model's performance. 
For example, $\ell$ may represent a prediction error (supervised learning), 
a reconstruction error (autoencoders), a divergence between distributions (generative models), 
or the violation of physical laws (SciML).

The general training objective is to find the parameters $\theta$ that minimize the \emph{expected risk} of the model, given by
\begin{align}
R(\theta)
  = \int_{\pazocal{X} \times \pazocal{Y}} 
        \ell(h_{\theta}(x), y) \, \mathrm{d}\pazocal{P}(x,y)
  = \mathbb{E}_{(x,y)\sim\pazocal{P}}
        \!\left[ \ell(h_{\theta}(x), y) \right].
        \label{eq:expected_risk}
\end{align}
However, since the distribution $\pazocal{P}$ is usually unknown in practice, learning relies on a finite dataset  $\pazocal{D} = \{(x_i, y_i)\}_{i=1}^m$,  consisting of $m$ samples drawn independently and identically from~$\pazocal{P}$,  i.e., $(x_i, y_i) \overset{\text{i.i.d.}}{\sim} \pazocal{P}$.

Using the dataset~$\pazocal{D}$, the expected risk~\eqref{eq:expected_risk} 
can be approximated, for example, using a Monte-Carlo estimate, yielding the 
\emph{empirical risk}
\[
\widehat{R}_m(\theta)
  =  \frac{1}{m}\sum_{i=1}^m \ell(h_\theta(x_i), y_i).
\]
Under standard assumptions (i.i.d.\ sampling, bounded variance),
$\widehat{R}_m(\theta)$ converges to $R(\theta)$ as $m \to \infty$; see,
e.g.,~\cite{Vapnik1998,ShalevShwartz2014}.

The \emph{Empirical Risk Minimization (ERM)} principle consists of minimizing $\widehat{R}_m$, 
which gives rise to the canonical \emph{finite-sum minimization problem}
\begin{align}
\min_{\theta \in \R^n} f(\theta) := \frac{1}{m} \sum_{i=1}^m f_i(\theta),
\qquad \text{where} \qquad
f_i(\theta) = \ell(h_\theta(x_i), y_i).
\label{eq:min_problem}
\end{align}
Here, each function $f_i(\theta)$ represents the contribution of one sample $(x_i, y_i)$ 
to the total objective $f(\theta)$. 
\newtext{In theoretical analyses, convergence of the iterates  generated by solving~\eqref{eq:min_problem} is typically 
established by showing that the gradient norm  $\|\nabla f(\theta)\|$ approaches zero, indicating that a  stationary point has been reached. 
In the ML context, however,  practitioners more commonly monitor the training loss  $f(\theta_k)$, the validation error on a held-out dataset,  or task-specific metrics such as accuracy. 
These practical metrics better reflect generalization  performance, which is the ultimate goal of ML training and  may not correlate directly with the gradient norm.}

\subsection{Examples of ERM across learning paradigms}
\label{sec:examples}
In this section, we illustrate how the abstract ERM formulation translates into practical learning settings by presenting three representative examples.

\begin{figure}[ht]
\centering
\includegraphics{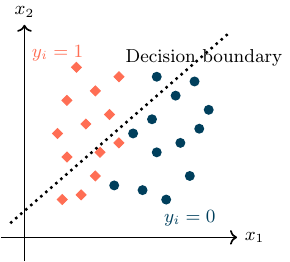}
\hfill
\includegraphics{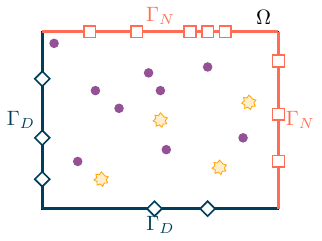}
\caption{
\emph{Left:} Example of binary classification in the $(x_1,x_2)$-plane: 
Samples from two classes $y_i \in \{0,1\}$ are shown as orange diamonds and blue circles, respectively. 
The dotted line indicates the decision boundary separating the two regions.
\emph{Right:} Example of a computational domain~$\Omega$ with sampled collocation points: 
Dirichlet boundaries~$\Gamma_D$ (\oldtext{blue}\newtext{diamond}), Neumann boundaries~$\Gamma_N$ (\oldtext{orange}\newtext{square}), 
interior points~$\pazocal{D}_\Omega$ (\oldtext{purple}\newtext{circle}), and optional empirical data samples~$\pazocal{D}_{\text{data}}$ (\oldtext{yellow}\newtext{star}).
}
\end{figure}

\subsubsection{{Example~1: Supervised learning -  classification}}
In supervised learning, we are given a labeled dataset $\pazocal{D} = \{(x_i, y_i)\}_{i=1}^m$,
where each input vector $x_i \in \mathbb{R}^d$ is associated with an output label
$y_i$, which indicates affiliation to a given class.
The goal is to learn the parameters $\theta \in \mathbb{R}^n$ of the model $h_{\theta} : \mathbb{R}^d \to \mathbb{R}$ to predict the label $y_i$ for a given input $x_i$ \newtext{and for unseen  examples}.
For simplicity, we restrict our attention to the binary classification setting \newtext{where $y_i\in \{0,1\}$.}

\paragraph{Logistic regression (linear model)}
In logistic regression, the model is
\[
h_{\theta}(x_i) = \rho(x_i^\top \theta),
\]
where $\rho(t) = \frac{1}{1 + e^{-t}}$ is the \emph{sigmoid} activation.
The discrepancy between prediction and label is measured by the logistic loss, given as
\begin{align}
\ell(h_{\theta}(x_i), y_i)
  = -\, y_i \log h_{\theta}(x_i)
    - (1 - y_i)\log\!\big(1 - h_{\theta}(x_i)\big).
\label{eq:logistic_loss}    
\end{align}
The ERM objective then becomes
\begin{align}
\widehat{R}_m(\theta)
 = \frac{1}{m} \sum_{i=1}^m
     \Big[
       -\, y_i \log \rho(x_i^\top\theta)
       - (1 - y_i)\log\!\big(1 - \rho(x_i^\top\theta)\big)
     \Big].
\label{eq:ERM_1}     
\end{align}
\newtext{This function is smooth and  \emph{convex} (strictly, cf.~Definition \ref{def:convex}, if the features are linearly independent and strongly,  cf.~Definition \ref{def:strongly_convex}, if  an $L^2$ regularization term is added).
Therefore, its minimization admits a unique minimum (and a unique minimizer in the strictly and strongly convex cases).}
The practical difficulty of minimizing~$\widehat{R}_m$ depends predominantly on the geometry of the data.
If the data are well-conditioned and the classes are separable, convergence is fast.
If, on the other hand, the samples are nearly collinear or the classes overlap, the optimization problem becomes more challenging.
\newtext{An example of implementing~\eqref{eq:ERM_1} in Python can be found in Code snippet \ref{lst:logistic_loss}.}

\begin{lstlisting}[language=Python, caption={\newtext{Example implementation of the logistic loss with $L^2$ regularization.}}, label=lst:logistic_loss]
def logloss(w, X, y, gamma=1e-4):
    # L2-regularized logistic loss averaged over samples
    # Logit z is the product of parameters (w) with input features (X)
    z = X @ w
    p = sigmoid(z)
    eps = 1e-12
    # Binary cross-entropy loss
    # Constant eps is added for numerical safety
    loss = -np.mean(y * np.log(p+eps) + (1-y)*np.log(1-p+eps)) 
    loss += 0.5 * gamma * np.sum(w**2)  # Add L2 penalty term
    return float(loss)
    \end{lstlisting}

\paragraph{Deep neural network (nonlinear model)}
If the data are not linearly separable, linear models may be insufficient to capture the underlying structure.
In such cases, DNNs provide flexible nonlinear alternatives capable of learning complex nonlinear relationships and feature interactions present in the data.

Let $h_\theta$ denote a feed-forward DNN composed of weight matrices $U_1,\, U_2,\, \dots,\, U_Q$ and bias vectors $b_1,\, b_2,\, \dots,\, b_Q$,
which can be collected into a parameter vector $\theta$ as $\theta = \mathrm{vec}(U_1,\dots,U_Q,b_1,\dots,b_Q) \in \mathbb{R}^n$.
A standard $Q$-layer network is defined recursively as
\[
z_0 = x,\qquad
z_j = \rho(U_j z_{j-1} + b_j),\quad j=1,\dots,Q-1,
\]
with the output
\begin{align}
h_\theta(x) = \rho(U_Q z_{Q-1} + b_Q).
\label{eq:dnn}
\end{align} 
Here, $\rho$ denotes an activation function of the user's choice, e.g., \emph{ReLU}, or \emph{tanh}. 
\newtext{A prototypical implementation of a DNN with two layers and \emph{ReLU} activation function in Pytorch can be found in  Code snippet~\ref{lst:DNN}}.

Given a dataset $\pazocal{D}$ and a loss function $\ell$, the empirical risk is now given as
\[
\widehat{R}_m(\theta)
  = \frac{1}{m} \sum_{i=1}^m 
      \ell\!\left(h_{\theta}(x_i),\, y_i\right).
\]
The choice of the loss function $\ell$ depends on the task at hand. 
For example, in a binary classification problem, one can use the logistic loss~\eqref{eq:logistic_loss}.
Here, we emphasize that the use of DNNs makes $\widehat{R}_m$ highly nonconvex.
The complexity of the resulting optimization problem  strongly depends on the data geometry and on the chosen DNN architecture.

\begin{lstlisting}[language=Python, caption={\newtext{Example of creating a small DNN in PyTorch.}}, label={lst:DNN}]
class DNN(nn.Module):
    # Small fully-connected ReLU network
    def __init__(self):
        super().__init__()
        # Linear layer with 2 input neurons and 20 neurons in hidden layer
        self.fc1 = nn.Linear(2, 20) 
        self.act = nn.ReLU() # Activation function
        self.fc2 = nn.Linear(20, 1) # Output layer with a single output

    def forward(self, x):
        # Forward pass through the network
        return self.fc2(self.act(self.fc1(x)))
\end{lstlisting}

\subsubsection{Example~2: Supervised learning - regression}
\label{sec:regression}
Regression is a supervised learning paradigm widely used in data-driven SciML.  
Given a dataset $\pazocal{D} = \{(x_i, y_i)\}_{i=1}^m$, the goal is to learn a mapping from inputs $x_i$ to real-valued outputs $y_i \in \mathbb{R}$ \newtext{that generalizes well to unseen data.}
For example, the input may consist of observational measurements, while the output might represent the solution of a PDE or a quantity derived from it.

In regression, the discrepancy between the model prediction and the true label is commonly measured using the squared $L^2$ loss.  
Defining the stacked data labels
$ y := \begin{bmatrix}
y_1, \dots, y_m
\end{bmatrix}^\top$ and model outputs
$
u(\theta) :=
\begin{bmatrix}
h_\theta(x_1), \dots, h_\theta(x_m)
\end{bmatrix}^\top$,
the ERM problem takes the following form:
\begin{equation}\label{eq:squared_loss}
\widehat{R}_m(\theta)
=
\frac{1}{2m}\,\|u(\theta) - y\|^2.
\end{equation}
This ERM formulation is commonly referred to as a \emph{least-squares} problem.

If the labels depend linearly on the inputs, this relationship can be modelled
using the linear model $h_\theta(x_i) = x_i^\top \theta$.
In this particular case, the ERM problem admits the compact representation $\widehat{R}_m(\theta) = \frac{1}{2m}\,\|X\theta - y \|^2$,
where $X \in \mathbb{R}^{m \times d}$ denotes the data matrix.
Minimizing this objective corresponds to solving a \emph{linear least-squares}
problem.
Thus, if $X$ has full column rank, the solution is given in closed form as $\theta = (X^\top X)^{-1} X^\top y$.

If, on the other hand, the model $h_\theta$ is nonlinear, e.g., a DNN  given in~\eqref{eq:dnn}, the ERM problem~\eqref{eq:squared_loss} becomes a \emph{nonlinear least-squares} problem.
Here, the term ``nonlinear'' refers to the nonlinear dependence of $h_\theta(x)$ on $\theta$, even though the loss remains quadratic in the residual (misfit).
The resulting optimization problem is generally \emph{nonconvex} and may admit multiple local minima and saddle points, making its minimization considerably challenging.

\subsubsection{Example~3: Physics-informed neural networks (PINNs)}
\label{sec:PINNs_definition}
In physics-informed learning, our goal is to approximate a solution of a differential equation using an ML model \newtext{\cite{raissi2019physics}}.
Let $\Omega \subset \mathbb{R}^d$ be a bounded domain with boundary 
$\partial \Omega = \Gamma_D \cup \Gamma_N$, where $\Gamma_D$ and $\Gamma_N$ 
represent \newtext{parts where Dirichlet and Neumann boundary conditions (BC) are imposed, respectively}.
We consider the following PDE:
\begin{equation}
\begin{aligned}
\pazocal{N}[u](x) &= q(x), && x \in \Omega, \\
u(x) &= g_D(x), && x \in \Gamma_D, \\
\partial_n u(x) &= g_N(x), && x \in \Gamma_N,
\end{aligned}
\label{eq:PDE}
\end{equation}
where $\pazocal{N}$ is a differential operator, $q$ is a source term, and $\partial_n u$ denotes the normal derivative on $\Gamma_N$.
The goal is to approximate the unknown solution $u:\Omega\to\mathbb{R}$  with a DNN~$h_\theta:\Omega\to\mathbb{R}$.
\newtext{The same formulation extends to time-dependent (non-stationary) PDEs by treating time as an additional input coordinate of $h_\theta$.
The computational domain $\Omega$ is then replaced by the space-time domain $\Omega \times [0, T]$, where $\Gamma_D$ and $\Gamma_N$ denote the parts of the spatial boundary where Dirichlet and Neumann conditions are enforced.
Moreover, an initial condition (IC) $u(x, 0) = u_0(x)$ is prescribed on $\Omega \times \{0\}$.}

In the continuous setting, model quality is quantified by an expected risk functional that penalizes the PDE residual
\begin{equation}\label{eq:pde_res}
    r_{\theta}(x):=\pazocal{N}[h_\theta](x) - q(x),
\end{equation}
boundary condition violations, and mismatch with available data, i.e.,
\begin{align}
R(\theta)
  &=  \gamma_{\Omega} \int_{\Omega} 
        \big|\pazocal{N}[h_\theta](x) - q(x)\big|^2 
        \, \mathrm{d}\pazocal{P}_\Omega(x)
   +  \gamma_{D} \int_{\Gamma_D} 
        \big|h_\theta(x) - g_D(x)\big|^2 
        \, \mathrm{d}\pazocal{P}_{\Gamma_D}(x) 
        \label{eq:pinn_expected_risk} \nonumber \\
  &\quad
   +  \gamma_{\Gamma_N} \int_{\Gamma_N} 
        \big|\partial_n h_\theta(x) - g_N(x)\big|^2 
        \, \mathrm{d}\pazocal{P}_{\Gamma_N}(x)\\
   & \quad + \gamma_{\Omega_{\text{data}}}  \int_{\Omega_{\text{data}}}
        \big|h_\theta(x) - y(x)\big|^2 
        \, \mathrm{d}\pazocal{P}_{\text{data}}(x),
        \nonumber
\end{align}
where $\gamma_{\Omega}, \gamma_{D}, \gamma_{\Gamma_N}, \gamma_{\Omega_{\text{data}}} \in \R^+$.

In practice, the integrals in~\eqref{eq:pinn_expected_risk} must be approximated. 
A common approach is to use Monte-Carlo (MC)  quadrature. Thus, we use the following datasets of randomly sampled collocation points
\begin{align*}
\pazocal{D}_\Omega      &= \{x_j^{(\Omega)}\}_{j=1}^{m_\Omega}, &
\pazocal{D}_{\Gamma_D}  &= \{x_r^{(D)}\}_{r=1}^{m_D}, \\
\pazocal{D}_{\Gamma_N}  &= \{x_k^{(N)}\}_{k=1}^{m_N}, &
\pazocal{D}_{\text{data}} &= \{(x_i^{(\text{data})}, y_i)\}_{i=1}^{m_{\text{data}}},
\end{align*}
to obtain the following empirical risk formulation:
\begin{align}
\widehat{R}_m(\theta)
  &= \frac{\gamma_\Omega}{m_\Omega} 
      \sum_{j=1}^{m_\Omega} 
        \big|\pazocal{N}[h_\theta](x_j^{(\Omega)}) - q(x_j^{(\Omega)})\big|^2
    + \frac{\gamma_D}{m_D} 
      \sum_{r=1}^{m_D} 
        \big|h_\theta(x_r^{(D)}) - g_D(x_r^{(D)})\big|^2 \nonumber\\
  &\quad
    + \frac{\gamma_N}{m_N} 
      \sum_{k=1}^{m_N} 
        \big|\partial_n h_\theta(x_k^{(N)}) - g_N(x_k^{(N)})\big|^2
    + \frac{\gamma_{\text{data}}}{m_{\text{data}}} 
      \sum_{i=1}^{m_{\text{data}}} 
        \big|h_\theta(x_i^{(\text{data})}) - y_i\big|^2\label{eq:risk_PINNs}.
\end{align}
\newtext{Code snippet~\ref{lst:pinn_loss} illustrates the  computation of the PDE residual loss for the 1D Poisson  equation using automatic differentiation in PyTorch.}

\begin{lstlisting}[language=Python, caption={\newtext{Pytorch implementation of the PDE residual loss  for the 1D Poisson equation $-u''(x) = f(x)$, using  automatic differentiation to compute the second order derivative  of the network output.}},label=lst:pinn_loss]
def pinn_second_derivative(model, x):
   # Compute second derivative by repeated automatic differentiation
   
    # Enable computation of gradient wrt. network inputs
    xg = x.clone().detach().requires_grad_(True)
    u = model(xg) # Forward pass through the network
    # Compute the first derivative
    du = torch.autograd.grad(u, xg, torch.ones_like(u), create_graph=True)[0]
    # Compute the second derivative
    d2u = torch.autograd.grad(du, xg, torch.ones_like(du), create_graph=True)[0]
    return d2u

def loss_fn(uxx, f_rhs):
    # Mean-squared of PDE-residual loss for -u''(x) = f(x)
    return ((uxx + f_rhs) ** 2).mean()
\end{lstlisting}

Here, we point out that more efficient quadrature techniques than MC can be used in order to reduce the approximation error.
For instance, \emph{quasi-Monte Carlo} (QMC) methods based on low-discrepancy sequences (e.g., Sobol’) can achieve approximation errors of order $\pazocal{O}(m^{-1}(\log m)^d)$, compared to the $\pazocal{O}(m^{-1/2})$ rate of standard MC under suitable regularity assumptions~\cite{dick2013high}.
Here, we emphasise that increasing the number of collocation points in PINNs is \emph{not} analogous to adding independent samples in classical ML, as the same deterministic PDE residual is evaluated at additional locations.  
Due to the global coupling induced by the differential operator, nearby points are often strongly correlated, so oversampling does not yield proportional accuracy gains~\cite{mishra2023estimates,lu2021deepxde}.

As before, in order to find the optimal parameters of model~$h_\theta$, we minimize~$\widehat{R}_m$. 
Although $h_\theta$ is usually a standard DNN, the choice of the architecture requires particular care because~$\widehat{R}_m$ involves \emph{spatial derivatives} of $h_\theta$.
In particular, the activation function~$\rho$ must be sufficiently smooth to allow stable differentiation through the PDE operator.
\newtext{Weak formulations of PINNs can relax these smoothness requirements~\cite{de2024wpinns}. 
However, for the strong formulation considered here, which is the most commonly used, smooth activations such as \emph{tanh}, \emph{sigmoid}, needs to be used. 
Alternatively, one might use \emph{SIREN} (Sinusoidal Representation Networks) architecture~\cite{sitzmann2020implicit} with sinusoidal activation functions. }

\begin{takeawaybox}
\newtext{
Training an ML model is typically formulated as minimizing an empirical risk function that aggregates individual sample losses into a finite-sum objective~\eqref {eq:min_problem}. 
The structure of this objective depends critically on the model and task.
For example, logistic regression yields a convex problem with a unique minimum, whereas DNNs yield highly nonconvex objectives with multiple local minima and saddle points. 
In PINNs, physical laws are incorporated into the loss as penalty terms, enforcing the PDE and boundary conditions in a soft manner. 
This introduces global spatio-temporal coupling across the domain, making the resulting optimization problem structurally different from classical ML objectives that naturally decompose over independent samples.
}
\end{takeawaybox}

  \section{Solving the optimization problem}
\label{sec:opt_problem}
Given a nonlinear, differentiable objective (loss) function $f: \mathbb{R}^n \rightarrow \mathbb{R}$,
the goal of numerical optimization \newtext{\cite{nocedal2006numerical}} is to solve
\[
\min_{\theta \in \mathbb{R}^n} f(\theta). 
\]
\newtext{In the ML context, this problem is often referred to as \emph{training}. }
\newtext{Such minimization problems are often solved approximately, until  a \emph{stationary point}, i.e., a point at which the gradient is zero, is reached. }
In most cases, this stationary point will be a \emph{local minimizer}, a point such that the function does not decrease in its neighborhood.
However, it is usually not possible to exclude it as a \emph{saddle point}.

To reach a stationary point, one typically employs \emph{iterative algorithms}, also referred to as \emph{optimizers}.  
Starting from an initial guess $\theta_0$, such methods generate a sequence $\{\theta_k\}_{k \in \mathbb{N}}$, where $k$ denotes the iteration index, that converges to a stationary point $\theta^\star$ through updates of the form
\[
\theta_{k+1} = \theta_k + \alpha_k\, p_k,
\]
where $p_k \in \mathbb{R}^n$ denotes a search direction and $\alpha_k > 0$ is the step size, often referred to as the \emph{learning rate}.

\newtext{
Many factors influence the convergence speed of the training.
On the problem side, these factors include the network architecture, the training data, and the non-convexity of the loss function. 
On the algorithmic side, they include the choice of the optimization algorithm, its hyperparameters, and the initialization strategy. 
The latter determines the initial network parameters $\theta_0$, which play the role of the initial guess for the optimization problem. 
This choice can strongly affect the optimization trajectory, especially in non-convex settings.
Several initialization strategies have been proposed in the literature for the weights of DNNs, e.g.,~\emph{Xavier} or \emph{Glorot initialization}~\cite{glorot2010understanding}, which is commonly used with $\tanh$ activations, and \emph{He initialization}~\cite{he2015delving}, which is typically used with ReLU activations.
}

\subsection{\newtext{A unifying perspective via preconditioning}}
\label{sec:precond_lens}
\newtext{
Throughout this chapter, we show that many optimization methods can be understood as forms of \emph{preconditioning}, directly linking algorithmic design choices to the conditioning of the underlying optimization problem.
Following~\cite{de2023operator}, we distinguish three types of preconditioning depending on the space in which they act.
Parameter-space preconditioners act directly on the gradient update,  e.g.,~by rescaling it (block) coordinate-wise (adaptive learning rates) or incorporating curvature information (second-order methods), see  Sections~\ref{sec:adaptive_methods},~\ref{sec:second_order}.
Data-space preconditioners transform the data, for example via input transformations, dataset resampling, or reweighting of the loss terms (Sections~\ref{sec:problem_def},~\ref{sec:data_precond},~\ref{sec:collocation_minibatch},~\ref{sec:outlook}).
Function-space preconditioners reshape the geometry of the function represented by the network, e.g.,\ via multiscale architectures  (Sections~\ref{sec:data_precond},~\ref{sec:outlook}).
}

\newtext{
Here, we note that this classification is not a strict partition.
Indeed, many strategies admit interpretations in more than one space (e.g.,\ some domain decomposition methods act in both parameter and function space). 
It should therefore be viewed as a helpful classification perspective rather than a strict taxonomy. 
Figure~\ref{fig:precond_lens} provides a schematic summary and serves as a roadmap for the strategies discussed in this chapter. Section~\ref{sec:outlook} returns to this perspective, surveying recent advances and open challenges in large-scale settings along each axis.
}

\begin{figure}[h!]
\centering
\tikzset{external/export=false}
\resizebox{\linewidth}{!}{\begin{tikzpicture}[
  font=\sffamily\small,
  node distance=3mm and 6mm,
  every node/.style={align=center},
  hdr/.style={
    rectangle, rounded corners=3pt,
    line width=0.6pt,
    inner xsep=8pt, inner ysep=5pt,
    minimum height=8mm, text width=46mm,
    font=\sffamily\bfseries,
  },
  cell/.style={
    rectangle, rounded corners=3pt,
    line width=0.4pt,
    inner xsep=6pt, inner ysep=4pt,
    minimum height=11mm, text width=46mm,
    font=\sffamily\footnotesize,
  },
]
\node[hdr, draw=c3, fill=c3!15, text=c3] (hP) {Parameter space};
\node[cell, below=of hP, draw=c3!70, fill=c3!8, text=c3!85!black] (P1)
  {Adaptive gradient methods\\ (AdaGrad, Adam) \S\ref{sec:adaptive_methods}};
\node[cell, below=of P1, draw=c3!70, fill=c3!8, text=c3!85!black] (P2)
  {(Approximate) second-order \\ methods (quasi-Newton, NK) \S\ref{sec:second_order}};
\node[cell, below=of P2, draw=c3!70, fill=c3!8, text=c3!85!black] (P3)
  {Domain decomposition \&\\ multilevel methods \S\ref{sec:outlook}};

\node[hdr, right=of hP, draw=c4, fill=c4!15, text=c4!85!black] (hD) {Data space};
\node[cell, below=of hD, draw=c4!70, fill=c4!8, text=c4!75!black] (D1)
  {Loss reweighting \&\\ soft-constraint balancing \S\ref{sec:problem_def}};
\node[cell, below=of D1, draw=c4!70, fill=c4!8, text=c4!75!black] (D2)
  {Adaptive / residual-based\\ sampling \S\ref{sec:collocation_minibatch}};
\node[cell, below=of D2, draw=c4!70, fill=c4!8, text=c4!75!black] (D3)
  {Input embeddings\\ (Fourier features) \S\ref{sec:data_precond}};

\node[hdr, right=of hD, draw=c6!90!black, fill=c6!25, text=c6!50!black] (hF) {Function space};
\node[cell, below=of hF, draw=c6!80!black, fill=c6!12, text=c6!50!black] (F1)
  {Multiscale architectures \S\ref{sec:data_precond}};
\node[cell, below=of F1, draw=c6!80!black, fill=c6!12, text=c6!50!black] (F2)
  {Sobolev training \S\ref{sec:outlook}};
\node[cell, below=of F2, draw=c6!80!black, fill=c6!12, text=c6!50!black] (F3)
  {Domain-decomposition\\ PINNs (XPINN, FBPINN) \S\ref{sec:outlook}};
\end{tikzpicture}}
\caption{\newtext{
Overview of the preconditioning classification. 
Optimization strategies discussed in this chapter are classified as preconditioners acting in one of the three spaces: parameter, data, or function space. 
Section references indicate where each strategy is presented in detail.
}}
\label{fig:precond_lens}
\end{figure}

{\subsubsection{Parameter-space preconditioning}}
{Assuming $f$ is sufficiently regular, most optimization methods build search directions from the derivatives of $f$, namely the \emph{gradient} $\nabla f(\theta)$ and the \emph{Hessian} $H(\theta)=\nabla^2 f(\theta)\in\mathbb{R}^{n\times n}$.
The update rule usually has the following form: 
\[
\theta_{k+1} = \theta_k - \alpha_k\, M_k^{-1} \nabla f(\theta_k),
\]
where $\alpha_k>0$ and  $M_k^{-1}$ is a preconditioner that rescales the gradient.}

Based on the information the preconditioner $M_k^{-1}$ uses, optimization methods are classically  divided into three large classes:
\begin{enumerate}
\item \emph{First-order methods}: 
These methods use only gradient information.  
As a consequence, their iteration cost is low, but they may require many iterations to converge. 
The canonical example is the gradient descent (GD) method, whose iterates are defined as follows:
\begin{align}
\theta_{k+1}=\theta_k-\alpha_k \nabla f(\theta_k),
\label{eq:first_order_it}
\end{align}
where $M_k=I$. 

\item \emph{Second-order methods}: 
These methods incorporate curvature information from the Hessian or its approximations, providing much faster local convergence at a higher per-iteration cost.
The most notable example is Newton's method, whose iterates are formally defined as follows:
\begin{equation}\label{eq:newton}
\theta_{k+1}=\theta_k-\alpha_k H(\theta_k)^{-1}\nabla f(\theta_k),
\end{equation}
where $M_k = H(\theta_k)$.
\newtext{The main cost associated with this method lies in evaluating the Hessian and solving the associated linear system
$
H(\theta_k)p_k = -\nabla f(\theta_k)
$
to compute the search direction $p_k$.} 
\newtext{The memory footprint is also a concern as storing $H(\theta_k)$ explicitly requires $\pazocal O(n^2)$ memory, which is prohibitive for large DNNs.}

\item \emph{Adaptive (preconditioned) gradient methods}:
Adaptive gradient methods address the challenge of balancing computational efficiency with curvature-awareness by performing updates of the form
\begin{align}
    \theta_{k+1} = \theta_k - \alpha_k M_k^{-1} \nabla f(\theta_k),
    \label{eq:adaptive_methods}
\end{align}
where $M_k$ is a diagonal preconditioner.
\newtext{As we will see in Section~\ref{sec:adaptive_methods},}
practical optimizers such as AdaGrad~\cite{duchi2011adaptive} and  Adam~\cite{kingma2017adammethodstochasticoptimization} construct diagonal $M_k$  from running estimates of the first and/or second moments of past  gradients, thereby avoiding the need to compute and store the exact Hessian,  while still exploiting curvature information.
\end{enumerate}

{\subsubsection{Data-space preconditioning}
A preconditioner $M^{-1}$ can also be applied to the model's input $x$, which alters both the search direction and its magnitude. 
For a model $h_\theta$ and a loss function $\ell$, the sample-wise objective becomes $f_i(\theta) = \ell\bigl(h_\theta(M^{-1}x),\, y\bigr)$, with gradient
\[
\nabla_\theta f_i(\theta)
=
J_{h_\theta}(M^{-1}x)\, M^{-1}\, \nabla_z \ell.
\]
Such transformations reshape the geometry of the input distribution and can significantly affect the conditioning and convergence behavior of the optimization process.
Loss reweighting (Section~\ref{sec:problem_def}), adaptive sampling (Section~\ref{sec:collocation_minibatch}), and input embeddings such as Fourier features (Section~\ref{sec:data_precond}) act in this space.}

{\subsubsection{Function-space preconditioning}
Function-space preconditioning refers to transformations that modify the geometry of the optimization problem directly in the space of functions represented by the DNN. 
This viewpoint is natural for SciML, where the learning objective is often defined through a differential operator~$\pazocal{N}$. Let $u$ denote the DNN output. 
For each sample, the objective reads $f(u)=\ell(\pazocal{N}(u),q)$, and the functional gradient is
\[
\nabla_u f = \pazocal{N}^*\!\left(\pazocal{N}(u) - q \right),
\]
where $\pazocal{N}^*$ denotes the formal adjoint of $\pazocal{N}$ with respect to the $L^2$ inner product,
defined by $\langle \pazocal{N} u, v\rangle = \langle u, \pazocal{N}^* v\rangle$ for all sufficiently smooth $u,v$. 
}

{A function-space preconditioner applies an operator $\pazocal{M}^{-1}$ to this gradient:
\[
u_{k+1} = u_k - \alpha_k\, \pazocal{M}^{-1}\nabla_u f,
\]
where ideally $\pazocal{M}\approx \pazocal{N}^*\pazocal{N}$. 
This mirrors classical PDE preconditioning and motivates several architectural and algorithmic strategies in SciML, in particular multiscale architectures (Section~\ref{sec:data_precond}) as well as Sobolev training and domain-decomposition PINNs (Section~\ref{sec:outlook}). Unlike in linear PDE solvers, however, the operator $\pazocal{N}$ is composed with a nonlinear neural representation, making the design and analysis of effective function-space preconditioners substantially more challenging~\cite{de2023operator}.}

\subsection{Neural tangent kernel}
\label{sec:NTK_illconditioning}
The Neural Tangent Kernel (NTK) is a theoretical tool introduced in \cite{Jacot2018_NTK} for studying the training dynamics of DNNs, especially in the “infinite‐width” limit (i.e., for networks with layers of infinite width). 
It connects neural network training with kernel methods,  providing a framework for understanding how GD shapes the function learned by a DNN.
In particular, the eigenvalue spectrum of the NTK governs the dynamics and conditioning of optimization methods and provides insights into the training behavior of DNNs.
For PINNs, NTK can be used to show that the optimization problem inherits stiffness (stiffness refers to the presence of widely separated scales - typically in time or space - that make numerical simulation difficult, since some components of the solution evolve very rapidly, while others evolve slowly) from the underlying differential operator~$\pazocal{N}$, as well as structural biases induced by the network architecture.
These effects make training notoriously challenging and often render first-order methods such as GD inefficient.

\subsubsection{Function-space view: GD as kernel flow}
In this section, we introduce the NTK from a function-space viewpoint.
Let $h_\theta : \Omega \to \mathbb{R}$ be a DNN with parameters $\theta \in \mathbb{R}^n$.
During GD training, the network updates its parameters, and the NTK describes how these parameter updates translate into changes in the network's output.

\paragraph{Continuous-time kernel dynamics}
The continuous-time version of GD, known as \emph{gradient flow}, updates the parameters according to
\[
\dot{\theta}(t) = - \nabla_\theta f(\theta(t)),
\]
where $\dot{\theta}(t) $ denotes the derivatives of ${\theta(t)}$ with respect to $t$. 
We are interested in the induced evolution of the function $h_\theta$.
Using the chain rule, and assuming that $f$ depends on 
$\theta$ only through $h_\theta$, we get
\begin{align}
\dot h_\theta(x) 
= \nabla_\theta h_\theta(x)^\top \, \dot\theta(t)
= - \nabla_\theta h_\theta(x)^\top \, \nabla_\theta f(\theta(t)).
\label{eq:ntk1}
\end{align}
If we consider the gradient of the loss, we have
\begin{align}
\nabla_\theta f = \int_\Omega \nabla_\theta h_\theta(x') \, \frac{\delta f}{\delta h(x')} \, dx',
\label{eq:ntk2}
\end{align}
where $\delta f / \delta h(x')$ is the functional derivative of $f$ with respect to $h(x')$.
Combining~\eqref{eq:ntk1} with~\eqref{eq:ntk2} gives
\[
\dot h_\theta(x) = - \int_\Omega K_\theta(x,x') \, \frac{\delta f}{\delta h(x')} \, dx',
\]
where 
\[
K_\theta(x,x') = \langle\nabla_\theta h_\theta(x), \nabla_\theta h_\theta(x')\rangle
\]
is the \emph{NTK}. 
Thus, when GD updates the parameters $\theta$, the corresponding evolution of the function $h_\theta$ is governed by the kernel $K_\theta$.  
In other words, \emph{gradient flow in parameter space induces a kernel gradient flow in function space with respect to the NTK metric}.

\paragraph{Discrete-time NTK dynamics}
The same argument can be developed in a discretized setting, leading to the empirical NTK. 
This can be simply shown by considering the empirical squared loss, given in \eqref{eq:squared_loss}. 
In this case, the GD step with fixed step size $\alpha>0$ reads as
\begin{equation}\label{eq:grad_erm}
\theta_{k+1} = \theta_k - \alpha\,\nabla f_m(\theta_k)
= \theta_k - \frac{\alpha}{m} J_k^\top \big(u(\theta_k)-y\big),
\end{equation}
where
\[
J_k =
\begin{bmatrix}
\nabla_{\theta} h_{\theta_k}(x_1)^\top \\
\vdots\\
\nabla_{\theta} h_{\theta_k}(x_m)^\top
\end{bmatrix}
\in\mathbb{R}^{m\times n}
\]
is the Jacobian of $h_{\theta_k}$ evaluated at all data points.

Following the linearization argument of~\cite{Jacot2018_NTK,lee2019wide}, we expand $h_\theta$ around $\theta_k$:
\begin{equation}\label{eq:NN_linear}
h_{\theta_{k+1}}(x)
  \approx h_{\theta_k}(x)
        + \nabla_\theta h_{\theta_k}(x)^\top (\theta_{k+1}-\theta_k).
\end{equation}
Since in our notation the predicted outputs evaluated on the training set satisfy $u(\theta_k) = h_{\theta_k}(x_{1:m})$, the linearization directly applies to $u(\theta_k)$ as well.
Substituting the GD update~\eqref{eq:grad_erm},
into this linearization then yields an approximate evolution equation for the predicted outputs:
\begin{equation}\label{eq:ntk_update}
u(\theta_{k+1}) \approx u(\theta_k) - \alpha\, \Theta_k\, (u(\theta_k)-y),
\qquad
\Theta_k := \frac{1}{m} J_k J_k^\top,
\end{equation}
where $\Theta_k \in \mathbb{R}^{m \times m}$ is the empirical NTK at iteration $k$.
 
This indicates that training evolves like kernel GD in the reproducing kernel Hilbert space (RKHS) induced by the NTK.
In sufficiently wide DNNs, or in the near-stationary regime (e.g., for small $\alpha$), the NTK remains nearly constant during training. 
In this case, the dynamics reduce to kernel GD with a fixed kernel, in turn providing a tractable mathematical description of deep learning training.
\newtext{We report below two examples of NTK evaluation in PyTorch, for a simple regression problem (Code snippet \ref{lst:NTK_naive}) and for a PINN trained to approximate the solution of the heat equation (Code snippet \ref{lst:NTK_PINN_naive}).
For PINNs, nested automatic differentiation is required, as the PINN loss itself involves spatial derivatives of the network output computed via autograd, and the NTK computation requires an additional differentiation with respect to the network parameters.
}

\begin{lstlisting}[language=Python, caption={\newtext{Example of computing the NTK for a simple regression problem. More elaborate implementations are discussed in Section~\ref{sec:ntk_practical_computations}.}},label=lst:NTK_naive]
def ntk_naive(net, X):
    # Compute Theta = (1/m) J J^T using naive loop
    params = [p for p in net.parameters() if p.requires_grad]
    rows = []
    # Loop over all data samples
    for i in range(X.shape[0]):
        out = net(X[i:i+1]).squeeze() # Forward pass through the network
        # Get the gradient of each output wrt. parameters 
        grads = torch.autograd.grad(out, params, retain_graph=False)
        rows.append(torch.cat([g.reshape(-1) for g in grads]))
    J = torch.stack(rows, dim=0)  # Stack together all gradients
    return (1.0 / X.shape[0]) * (J @ J.T), J
            
\end{lstlisting}

\begin{lstlisting}[language=Python, caption={\newtext{Example of evaluating the NTK for PINN (heat equation). Note that, since computing the loss requires the evaluation of the residual, computing the NTK requires nested autograd.}},label=lst:NTK_PINN_naive]
def f_source(xt):
    # Evaluate source term
    x, t = xt[..., 0], xt[..., 1]
    return torch.exp(-t) * (PI**2 - 1.0) * torch.sin(PI * x)

def residual_at(net, xt):
    # Evaluation of the PDE residual, which uses NESTED autograd, i.e.,   
    # the inner grad call must keep the graph alive (create_graph=True) so 
    # that  u_t and u_xx remain differentiable  in theta for computation of NTK
    xt = xt.detach().clone().requires_grad_(True)
    
    u = net(xt).squeeze(-1) # Evaluate the network output
    # Compute gradients with respect to input
    grads = torch.autograd.grad(u.sum(), xt, create_graph=True)[0] 
    u_x, u_t = grads[:, 0], grads[:, 1]
    u_xx = torch.autograd.grad(u_x.sum(), xt, create_graph=True)[0][:, 0]  
    return u_t - u_xx - f_source(xt) # Residual evaluation


def ntk_naive_pinn(net, xt):
    # Same as ntk_naive function, but per-sample function is the PDE residual   
    params = [p for p in net.parameters() if p.requires_grad]
    rows = []
    # Loop over all data samples
    for i in range(xt.shape[0]):
        # Evaluate residual for each data point
        r_i = residual_at(net, xt[i:i+1]).squeeze()
        # Compute gradients of residual wrt. parameters
        grads = torch.autograd.grad(r_i, params, retain_graph=False, allow_unused=True)
        flat = [(g if g is not None else torch.zeros_like(p)).reshape(-1)
                for g, p in zip(grads, params)]
        rows.append(torch.cat(flat))
    J = torch.stack(rows, dim=0)  # Stack together all gradients
    return (1.0 / xt.shape[0]) * (J @ J.T), J
    
\end{lstlisting}

\subsubsection{Parameter-space viewpoint}\label{sec_param_space}
The function-space analysis above can be complemented by a parameter-space perspective.  
For instance, for the empirical least-square loss \eqref{eq:squared_loss}, 
the gradient and Hessian are given as
\[
\nabla \hat{R}_m(\theta_k)=\tfrac{1}{m}J_k^\top(u(\theta_k)-y), 
\qquad
H(\theta_k)=\tfrac{1}{m}J_k^\top J_k + \text{(higher-order terms)}.
\]
If we are in the NTK regime, i.e., the network is wide and approximately linear in  $\theta$ \newtext{(cf.~\eqref{eq:NN_linear})}, the higher-order terms are small, and $H$ can be approximated by the Gauss-Newton (GN) approximation $\tfrac{1}{m}J_k^\top J_k$.
Since the empirical NTK is ${\Theta_k=\tfrac{1}{m}J_k J_k^\top}$,
it follows that \emph{$H(\theta_k)$ and $\Theta_k$ share the same nonzero eigenvalues}.  
This creates a direct correspondence between NTK conditioning in function space and Hessian conditioning in parameter space. 
In particular, it implies that small eigenvalues of $\Theta_k$ imply small eigenvalues of $H(\theta_k)$,  so that flat directions in the objective function landscape correspond exactly to slow-learning directions in function space.

As a consequence, if the NTK and the Hessian are ill-conditioned, this negatively impacts the convergence of the first-order methods. 
For example, let us consider the GD iteration~\eqref{eq:first_order_it} applied to a strongly convex quadratic model with constant Hessian $H$ and step size $0 < \alpha < 2/\lambda_{\max}(H)$.  
It is well known that the GD iterates satisfy the following classical contraction bound (see, e.g.,~\cite[Section~3.2]{nocedal2006numerical}):
\[
\|\theta_{k+1}-\theta^{\star}\|
  \le 
  \left(\frac{\kappa(H)-1}{\kappa(H)+1}\right)
  \|\theta_k-\theta^{\star}\|, 
\qquad 
\kappa(H)=\frac{\lambda_{\max}(H)}{\lambda_{\min}(H)}.
\]
Thus, if $\lambda_{\min}(H)$ is very small, the contraction factor approaches one and progress becomes arbitrarily slow.  
Similar observations extend to stochastic gradient methods, which we discuss in the upcoming sections.

\subsection{Explicit mode-wise error decay and spectral bias}
\label{sec:spectral_bias}
The eigenvalues of $\Theta_k$  determine the contraction rates of the corresponding error modes during GD training. 
To understand this behavior, let us again consider the empirical least-squares loss~\eqref{eq:squared_loss}. 
The function-space error is defined as $e_k := u(\theta_k)-y$. 
Using NTK-induced update~\eqref{eq:ntk_update}, we can obtain a linear recursion of the following form:
\begin{equation}\label{eq:lin_dyn}
e_{k+1} \approx (I - \alpha \Theta_k)\, e_k.
\end{equation}
In regimes where $\Theta_k$ varies slowly,
we may approximate~$\Theta_k$ locally with constant $\Theta$, i.e.,~$\Theta_k \approx \Theta$, yielding
\[
e_k \approx (I - \alpha \Theta)^k e_0,
\]
which is the standard linear kernel gradient-flow approximation
used in~\cite{Jacot2018_NTK,lee2019wide}.

Let the eigendecomposition of $\Theta$ be
\[
\Theta = Q \Lambda Q^\top, 
\qquad  \qquad
\Lambda=\mathrm{diag}(\lambda_1,\ldots,\lambda_m), 
\]
where $q_i$ denotes the $i$-th eigenvector (column of $Q$) and $\lambda_i$ is the associated eigenvalue.
Because $\{q_i\}_{i=1}^m$ form an orthonormal basis of $\mathbb{R}^m$, we may expand the initial error as
$e_0 = \sum_{i=1}^m c_i\, q_i$. 
Applying the recursion $e_{k+1} = (I - \alpha \Theta_k)e_k$ repeatedly, then yields
\[
e_k
  = \sum_{i=1}^m c_i (1-\alpha\lambda_i)^k q_i \qquad \text{ and }  \qquad
\|e_k\|^2
  = \sum_{i=1}^m c_i^2 (1-\alpha\lambda_i)^{2k}. 
\]
Thus, each coefficient $c_i$ associated with eigenvector $q_i$ is multiplied by the scalar factor  $(1 - \alpha\lambda_i)$ at every iteration.  
This implies that the $i$-th component decays geometrically at rate  $(1 - \alpha\lambda_i)^k$.
Hence, if $\lambda_i$ is large, the factor $1 - \alpha\lambda_i$ is significantly below one, which leads to a rapid decrease of the corresponding error mode.  
In contrast, if $\lambda_i$ is very small, then  $1 - \alpha\lambda_i \approx 1$, so the contraction per iteration is extremely small.

This imbalance in decay rates across different eigenmodes is referred to as \emph{spectral bias}~\cite{rahaman2019spectral,wang2021eigenvector}. 
DNNs learn components of the solution that align with large NTK  eigenvalues (typically low-frequency or smooth structures) much faster than the components associated with small eigenvalues (often high-frequency components).
As a consequence, the smallest \newtext{nonzero} NTK eigenvalue $\lambda_{\min}(\Theta)$ governs the slowest decaying error and directly controls the overall speed of learning.

Here, we emphasize that different DNN architectures induce different NTKs and, therefore, exhibit different implicit biases. 
Specialized architectures can be designed to mitigate spectral bias; \newtext{see, for example, Fourier features~\cite{tancik2020fourier}, discussed in Section~\ref{sec:data_precond} and Numerical Example \ref{sec:ff_multiscale_example}}.

\newtext{We report the Python code to evaluate the condition number of a kernel matrix $K$ in Code snippet \ref{lst:cond}. We remind that the matrix $K$ depends on the chosen optimization method, e.g., for  GD, $K=\Theta_k$.} 

\begin{lstlisting}[language=Python, 
caption={\newtext{Computing the condition number of the kernel matrix $K$.}},label=lst:cond]
 def cond_number(K):
     # Compute singular values (more numerically stable than eigenvalues)
     evals = torch.linalg.svdvals(K.cpu()).numpy()
     eps = 1e-9
     # Filter out numerical noise and zero eigenvalues
     pos = evals[evals > eps]
     if len(pos) == 0:
         return np.inf, pos
     # Condition number = sigma_max / sigma_min
     return pos.max() / pos.min(), pos
\end{lstlisting}

\begin{figure}[t]
  \centering
\includegraphics{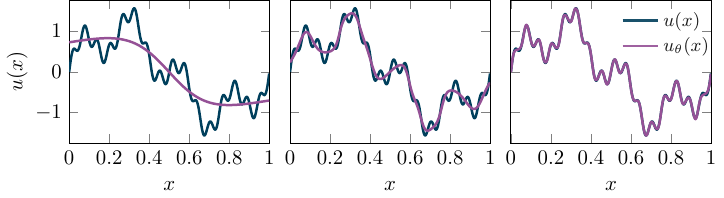}
\caption{Example of spectral bias to learn target function $u_{\theta}$ (purple) composed of low-, medium-,  and high-frequencies. 
    The network prediction $u_{\theta}$ (purple) is reported at training iterations 100 (left), 1,500 (middle), and 7,500 (right).}
    \label{fig:spectral_bias}
\end{figure}

\begin{num_example}{\textbf{\emph{Spectral bias.}}}\label{num_ex:spectral_bias}
To illustrate the spectral bias phenomenon in practice, we consider nonlinear regression as specified in Section~\ref{sec:regression} to learn the function $u(x) = \sin(2\pi x) \;+\; 0.5\,\sin(8\pi x) \;+\; 0.2\,\sin(32\pi x)$, which is a sum of sinusoids of increasing frequency. 
The training is performed using the GD method. 
As we can see from Figure~\ref{fig:spectral_bias}, in early training (iteration $100$), the network fits only the low frequencies. 
Medium frequencies appear later (iteration $1,500$), and only after many more iterations, the network begins to capture the high-frequency oscillations (iteration $7,500$). 
\end{num_example}

\subsection{Kernels of adaptive and second-order methods}
\label{sec:precond_kernel}
We have seen that the conditioning of the NTK and the Hessian plays a central role in determining the training dynamics of the GD method.  
A similar analysis extends to adaptive and second-order methods.  
Using the same linearization as in~\eqref{eq:ntk_update}, and writing  $ \nabla \hat{R}_m(\theta_k) = \tfrac{1}{m} J_k^\top \bigl(u(\theta_k)-y\bigr)$,
the update rule~\eqref{eq:adaptive_methods} gives
\[
\theta_{k+1} - \theta_k
  = -\frac{\alpha_k}{m}\, M_k^{-1} J_k^\top \bigl(u(\theta_k)-y\bigr),
\]
where $M_k^{-1}$ denotes the preconditioning matrix used by the optimizer.  
As before, we linearize the network output $h_\theta$ around $\theta_k$ at every training point and stack the resulting approximations into the vector $u(\theta)$.  
This yields the functional update
\[
u(\theta_{k+1})
  \approx 
  u(\theta_k)
  - \alpha_k\, \Theta_k^{(M_k)} \bigl(u(\theta_k)-y\bigr),
\]
with the \emph{adaptive (or preconditioned) kernel}
\begin{equation}
\Theta_k^{(M_k)}
  := \frac{1}{m}\, J_k\, M_k^{-1}\, J_k^\top.
\label{eq:precond_ntk}
\end{equation}

Thus, every adaptive and second-order method can be interpreted in function space as performing a kernel GD with kernel $\Theta_k^{(M_k)}$.  
The eigenvalues of $\Theta_k^{(M_k)}$ therefore govern the contraction rates of the error for a given optimizer.
In the special case of standard GD ($M_k = I$), we recover $\Theta_k^{(M_k)} = \Theta_k$, whose eigenvalues $\{\lambda_i\}$ produce per-iteration decay factors of the form $(1 - \alpha \lambda_i)^k$.  
Adaptive and second-order optimizers reshape these decay rates through the spectrum of $\Theta_k^{(M_k)}$.  
When $M_k$ is chosen well, the resulting spectrum is significantly more balanced, narrowing the gap between the fastest- and slowest-contracting modes, \newtext{so that the singular values are more similar in magnitude.
This phenomenon, often referred to as \emph{spectral flattening},} improves the effective conditioning of the optimization problem $\big({\kappa({\Theta_k^{(M_k)}}) < \kappa({\Theta_k}})\big)$ and typically leads to faster and more stable training.

\begin{figure}[t]
  \centering
\includegraphics{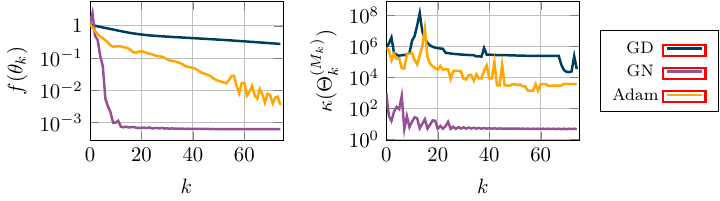}
\caption{Comparison of the convergence behaviors of GD, the adaptive gradient method (Adam), and the Gauss-Newton (GN) method. 
\emph{Left:} Evolution of the function value over the iterations. 
\emph{Right:} Condition numbers of the corresponding empirical kernels  $\kappa(\Theta_k^{(M_k)})$.}
\label{fig:kernel_condition_numbers}
\end{figure}

\begin{num_example}{\textbf{\emph{Preconditioning.}}}
To illustrate the effect of preconditioning on the optimization dynamics, we consider a nonlinear regression problem as specified in
Section~\ref{sec:regression}.
Our goal is to learn the function $\sin(\pi x_{1}) + \cos(\pi x_{2})$ from inputs $(x_1, x_2)$ using a simple DNN \newtext{with one hidden layer}.
The model is trained using three different optimization methods.
Specifically, we consider first-order GD, Adam as a representative of adaptive first-order methods \newtext{(Section~\ref{sec:adaptive_grad_algos} below)}, and a damped Gauss-Newton (GN)
method as an approximate second-order optimization \newtext{(Section~\ref{sec_param_space})}.
GD corresponds to the choice $M_k = I$, while Adam employs a diagonal preconditioner $M_k$, \newtext{see the update rule~\eqref{eq:adam_update})}.
The damped GN method (\newtext{or Levenberg-Marquardt method)~\cite[ch. 10.3]{nocedal2006numerical} } uses $M_k = J_k^\top J_k + \beta I$, where $J_k$ denotes the Jacobian of the model outputs with respect to the parameters and $\beta > 0$ is a damping parameter.

Figure~\ref{fig:kernel_condition_numbers} reports the evolution of the objective function and the condition number of the associated kernel~$\Theta_k^{(M_k)}.$ 
As we can see, GD produces the most ill-conditioned kernel, which correlates with slow convergence.
Adam's diagonal preconditioning partially alleviates this issue, leading to slightly faster convergence.
Notably, the GN method yields a substantially more balanced kernel spectrum, resulting in significantly smaller values of~$\kappa(\Theta_k^{(M_k)})$ and the fastest empirical decrease of the objective function among the methods considered.
\newtext{To demonstrate the computation of $\Theta_k^{(M_k)}$ in  practice, we provide in Code Snippet~\ref{lst:GN_K} an explicit  construction for the damped GN method.}

\begin{lstlisting}[language=Python, caption={\newtext{Example of constructing the kernel matrix $K=\Theta_k^{(M_k)}$ for the (damped) GN method.}},label=lst:GN_K]
def compute_GN_kernel(model, X, beta):
        # Compute GN kernel: K_GN = J * (J^T J + lambda I)^{-1} * J^T        
        # Damping is controlled by parameter beta
        B = J.t() @ J + beta * torch.eye(P, dtype=torch.float64)
        B_inv = torch.linalg.inv(B)
        K_gn = J @ B_inv @ J.t()
        
\end{lstlisting}

\end{num_example}
\subsection{Impact of the physical constraints on the loss landscape and the NTK of PINNs}
\label{sec:NTK_diff_op}

A key difficulty in training PINNs arises from the fact that the loss function does not act directly on the network output $h_\theta$, but rather on its derivatives through the PDE operator.  
This has a profound influence on the geometry of the objective function landscape and on the conditioning of the NTK, and ultimately on the behaviour of employed optimization methods.

To understand this effect, let us consider the decomposition of $h_\theta$ in the  Fourier components, i.e., 
$$
h_\theta(x)=\sum_{\omega}\hat{h}(\omega)e^{i\omega x}
$$
with $\hat{h}$ being the Fourier transform of $h_\theta$. 
Let us now consider a Fourier mode $e^{\mathrm{i}\omega x}$ present in $h_{\theta}(x)$. 

Applying $p$ derivatives multiplies this mode by $\omega^p$, so 
\[
\partial_x^p h_{\theta} \sim \omega^p, 
\qquad \qquad
J_k \sim \omega^p,
\qquad  \qquad
H(\theta)\newtext{\approx}\tfrac{1}{m}J_k^\top J_k \sim \omega^{2p}.
\]
Thus, each Fourier mode of frequency $\omega$ approximately satisfies
\[
H(\theta)\, q_\omega \approx \omega^{2p} q_\omega.
\]
This produces a Hessian spectrum whose eigenvalues span
\[
\lambda_{\min}(H(\theta_k)) \approx \omega_{\min}^{2p},
\qquad \qquad
\lambda_{\max}(H(\theta_k)) \approx \omega_{\max}^{2p},
\]
where $\omega_{\min}$ and $\omega_{\max}$ are the smallest and largest frequencies that the network effectively represents.
The  condition number thus scales as $\kappa(H)  \approx  \left(\frac{\omega_{\max}}{\omega_{\min}}\right)^{2p}$.
As a consequence, as the differential order $p$ increases, the exponent $2p$ causes the ill-conditioning to worsen dramatically.  
This explains why higher-order PDEs exhibit far more severe NTK/Hessian pathologies than first- or second-order PDEs.

The same scaling directly affects the curvature of the PINN objective function.
Low-frequency modes create shallow, nearly flat directions in the loss landscape, while high-frequency modes produce extremely steep ravines.  
This anisotropy makes the optimization problem stiff and highly challenging to solve.
First-order methods tend to move predominantly along low-frequency directions due to spectral bias, while second-order methods must cope with strong curvature anisotropy.

\begin{figure}
\centering
\includegraphics[scale=0.5]{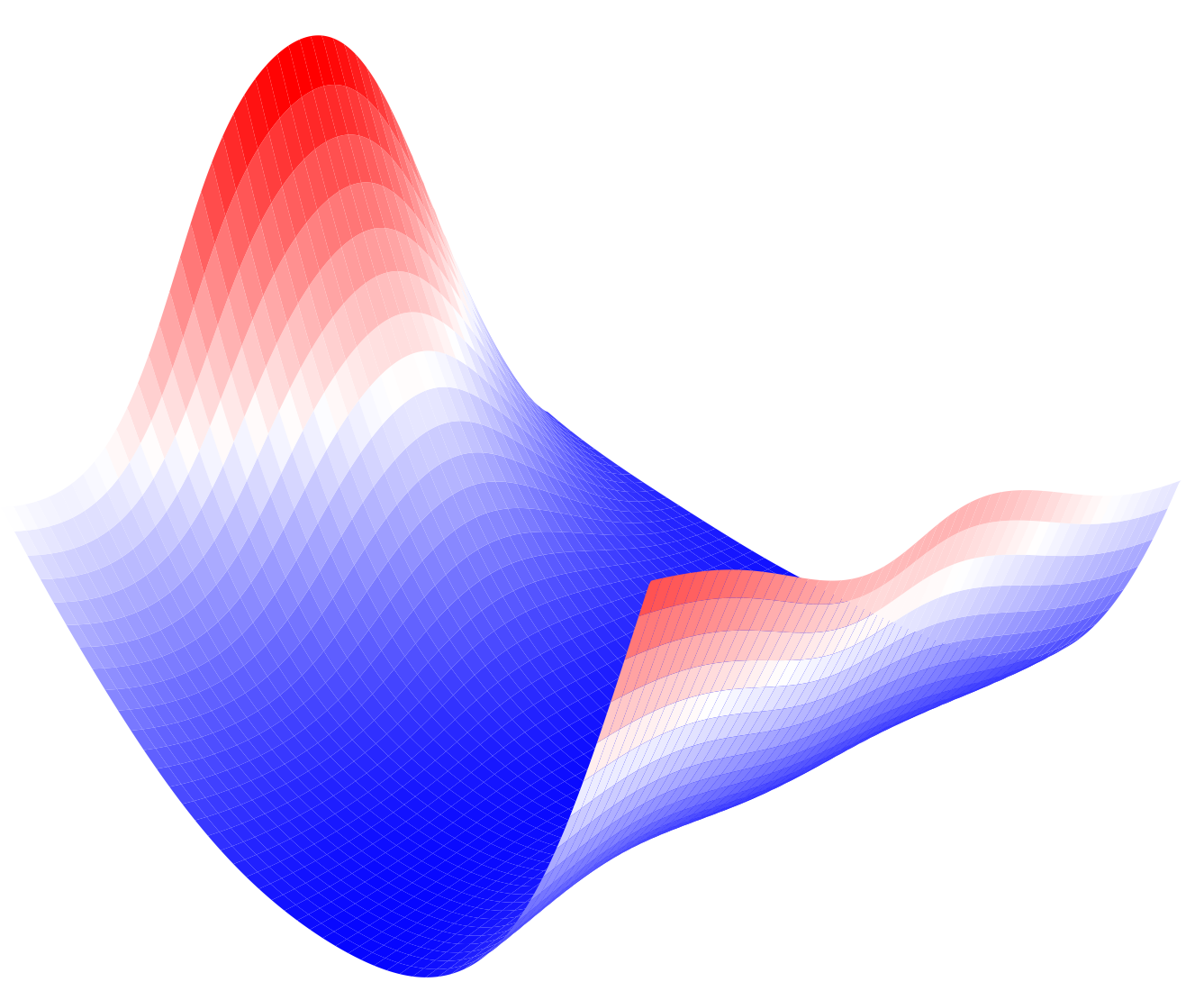}
\hfill
\includegraphics[scale=0.55]{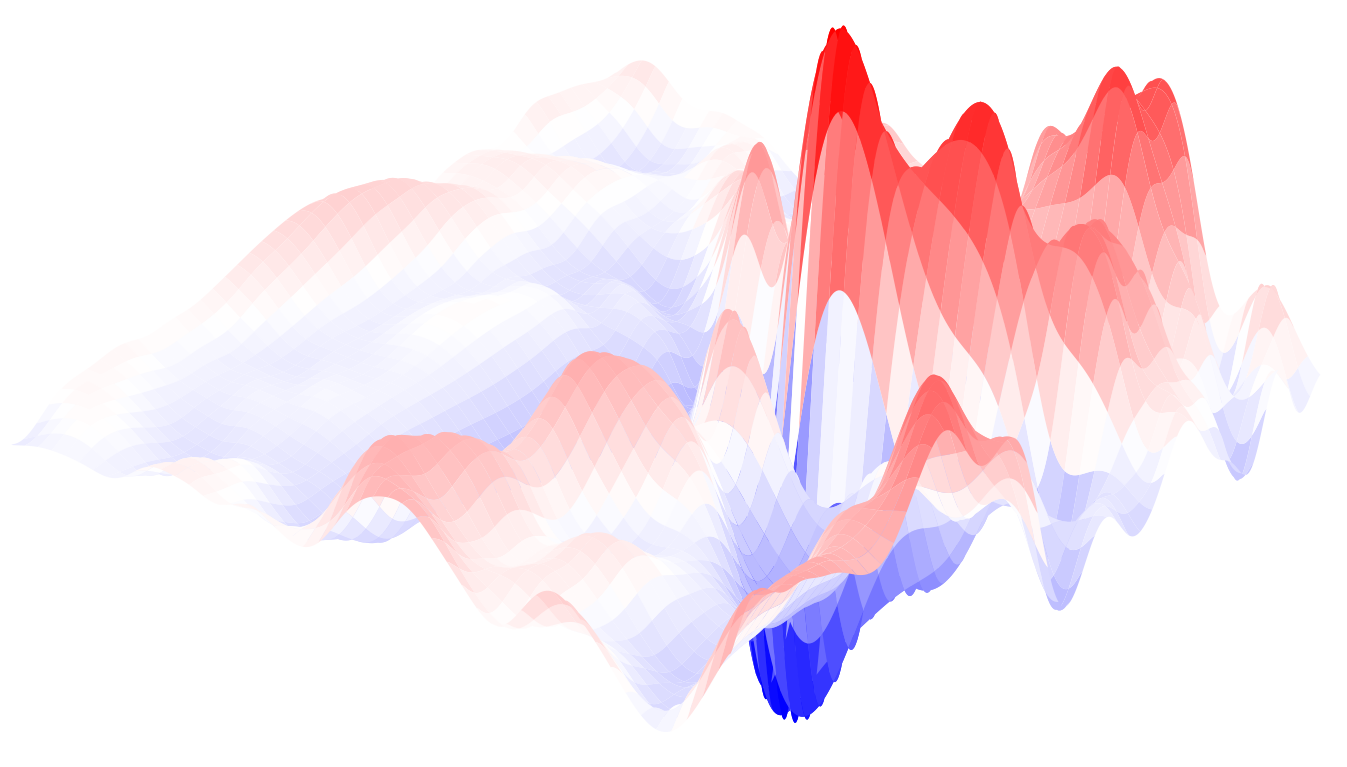}
\caption{Two-dimensional projection~\cite{li2018visualizing} of the objective function landscape for the same DNN under two training regimes for \newtext{1D Poisson problem, where we learn a map from right-hand side to a solution}.  
\emph{Left:} Data-driven nonlinear regression. 
\emph{Right:} PINN.}
    \label{fig:landscapes_pinn_data}
\end{figure}

To illustrate these concepts, we consider two DNNs with identical architectures and parameter values to learn the solution of a one-dimensional Poisson equation.
The first model is trained in a purely \emph{data-driven} setting, formulated as a nonlinear regression problem using samples of the exact solution (Section~\ref{sec:regression}).
The second model is trained as a \emph{PINN} (Section~\ref{sec:PINNs_definition}) by minimizing the PDE residual and the associated BC.
Figure~\ref{fig:landscapes_pinn_data} compares the associated loss landscapes.
Although both models approximate the same function $u$, the PINN loss landscape exhibits markedly more anisotropic curvature, while the data-driven objective remains comparatively smoother,  albeit with several flat or low-curvature regions.

\subsection{\newtext{Practical computations with the NTK}}
\label{sec:ntk_practical_computations}

\newtext{
All NTK arguments in this chapter rely on quantities derived from the empirical kernel
\begin{align}
\Theta_k \;=\; \tfrac{1}{m}\,J_k J_k^\top \in \mathbb R^{m\times m},
\qquad
(J_k)_{ij} \;=\; \tfrac{\partial h_{\theta_k}(x_i)}{\partial \theta_j},
\qquad J_k \in \mathbb R^{m\times n},
\label{eq:emp_ntk}
\end{align}
where $m$ is the number of inputs $\{x_i\}_{i=1}^m$ and $\theta\in\mathbb{R}^n$. 
While this representation is conceptually simple, computing the empirical NTK for modern DNNs is generally challenging due to the prohibitive computational and memory costs associated with forming and storing the Jacobian $J_k$, whose size scales as $mn$.
}

\newtext{In practice, one of the following three standard strategies  can be employed to compute $\Theta_k^{(M_k)}$ exactly, each suited  to a different computational regime:
\begin{itemize}
\item \textbf{Naive autograd loop}: We construct $J_k$ row by row by performing a backward pass for each sample. This approach is straightforward to implement and requires only standard automatic differentiation tools. However, it incurs $m$ backward passes, leading to significant runtime overhead due to repeated graph traversals and poor hardware utilization. 
Its memory cost is $\pazocal{O}(mn)$, if $J_k$ is stored explicitly, and it is typically impractical except for very small datasets.

\item \textbf{Jacobian contraction}: Instead of performing a naive loop over the samples, the full Jacobian $J_k$ can be assembled in a vectorized manner.
This can be achieved by using batched automatic differentiation, which applies the gradient computation to all the samples in the batch at once. 
The NTK is then obtained via the matrix product $\Theta_k = \frac{1}{m}J_k J_k^\top$. 
This approach is often significantly faster than the naive method. However, it still requires storing $J_k$, leading to $\pazocal{O}(mn)$ memory usage, and the matrix multiplication itself can become a bottleneck. 
This approach is preferred when both $m$ and $n$ are moderate and the matrix fits in the memory. An example of such calculation in PyTorch via the \texttt{torch.func.vmap} function is presented in Code snippet~\ref{lst:vmap}.

\item \textbf{Matrix-free NTK-vector products}~\cite{novak2022fastfinitewidth}: 
To reduce memory requirements, one can directly exploit the action of the kernel $\Theta_k$ on a vector $v \in \mathbb{R}^m$, i.e., the matrix-vector product $\Theta_k v$, without forming $\Theta_k$ explicitly.
This approach leverages forward- and reverse-mode automatic differentiation: a \texttt{vjp} computes $J_k^\top v$ in parameter space via a single backward pass, while a \texttt{jvp} maps this back to output space via a forward pass, yielding $J_k J_k^\top v = m \Theta_k v$.
Crucially, neither $J_k$ nor $\Theta_k$ need to be formed explicitly,  requiring only $\pazocal{O}(n)$ memory.
We refer the reader to Code Snippet~\ref{lst:vjp} for a PyTorch implementation.
\end{itemize}}

\newtext{The trade-off between Jacobian contraction and the matrix-free formulation is governed by the number of required matrix-vector products and memory constraints. 
If the full matrix $\Theta_k$ is needed (e.g., for direct eigendecomposition or visualization), Jacobian contraction is preferable. 
If only partial spectral information (e.g., extremal eigenvalues) or a small number of linear solves requiring $\Theta_k$ is desired, the matrix-free approach becomes advantageous. 
In practice, a few dozen matvecs often suffice when Krylov methods are used for the linear solves, making this approach cheaper than the contraction method once the cost of forming $J_kJ_k^\top$ dominates.
}

\begin{lstlisting}[language=Python, caption={\newtext{Example of evaluating the empirical NTK via Jacobian contraction using \texttt{torch.func.vmap} function.}},label=lst:vmap]
def fnet_single(params, x):
    # Scalar network output for a single input
    return functional_call(net_reg, params, (x.unsqueeze(0),)).squeeze()

def ntk_contraction(net, X, fn_single):
    # Vectorized approach: assemble J in one pass, then form Theta = (1/m) J J^T
    params = {k: v.detach() for k, v in net.named_parameters()}
    # Apply vmap over samples
    # Function jacrev returns the gradient of the scalar h_theta wrt theta
    jac = vmap(jacrev(fn_single, argnums=0), (None, 0))(params, X)
    cols = [jac[k].reshape(X.shape[0], -1) for k in jac] 
    J = torch.cat(cols, dim=1) 
    return (1.0 / X.shape[0]) * (J @ J.T), J

# Get NTK using contraction approach
Theta_contr, J_contr = ntk_contraction(net_reg, x_reg, fnet_single)
\end{lstlisting}

\begin{lstlisting}[language=Python, caption={\newtext{Example of matrix-free computation of NTK-vector product.}},label=lst:vjp]
def make_ntk_matvec(net, X, fn_batch):
    # Return a closure Theta v
    # Argument fn_batch must take a params dict and return the batched output
    
    params = {k: v.detach() for k, v in net.named_parameters()}
    m_local = X.shape[0]
    def matvec(v):
        # Step 1: J^T v via reverse-mode
        _, vjp_fn = vjp(lambda p: fn_batch(p, X), params)
        (jtv,) = vjp_fn(v)
        # Step 2: J (J^T v) via forward-mode
        _, theta_v = jvp(lambda p: fn_batch(p, X), (params,), (jtv,))
        return theta_v / m_local
    return matvec

def fnet_batch(params, X):
    # Vector-valued network prediction 
    # Returns the m predictions stacked into (m,)
    return functional_call(net_reg, params, (X,)).squeeze(-1)
    
# Get NTK*v product
matvec = make_ntk_matvec(net_reg, x_reg, fnet_batch)
# Build full Theta from the matrix-free form by applying it to the canonical basis (we evaluate m matvecs)
Theta_implicit = torch.stack([matvec(torch.eye(m)[i]) for i in range(m)], dim=1)
\end{lstlisting}

\newtext{
Finally, we note that in large-scale settings, even matrix-free  NTK computations may become prohibitively expensive.
In such cases, several approximations have been proposed, for instance via random parameter projections~\cite{halko2011finding}, last-layer  NTKs~\cite{lee2019wide}, or closed-form infinite-width  NTKs~\cite{arora2019exact,novak2019neural}.
}

\begin{num_example}{\textbf{\emph{Evaluating the empirical NTK.}}}
\label{sec:ntk_practical_example}
\newtext{
Figure~\ref{fig:ntk_scaling} reports the wall time (left) and the dominant tensor memory (right) of the three exact approaches for computing the NTK.
Experiments are performed using a four-layer tanh-MLP with 32 width ($n=3{,}265$ parameters).
We consider the number of inputs $m$ from $32$ to $1,024$. 
The reported memory is associated with the dominant tensor allocation, i.e., $J_k\in\mathbb R^{m\times n}$ for the naive loop and contraction approach. 
For the matrix-free approach, only a constant working set of size $\sim\!4n$ is required.
}

\newtext{
As we can see from the observed results, the naive loop is never competitive.
Jacobian contraction is the practical default whenever $mn$ fits in memory, scaling almost linearly in $m$ in both time and memory. 
Moreover, we also observe that the matrix-free approach is the qualitative outlier in terms of the memory requirements. 
Its size is set by the network size, not by $m$, and stays constant throughout the sweep.
 Combined with the flat per-matvec cost, this makes it the only viable algorithm for medium- to large-scale DNNs.
 }

\begin{figure}[h!]
\centering
\includegraphics{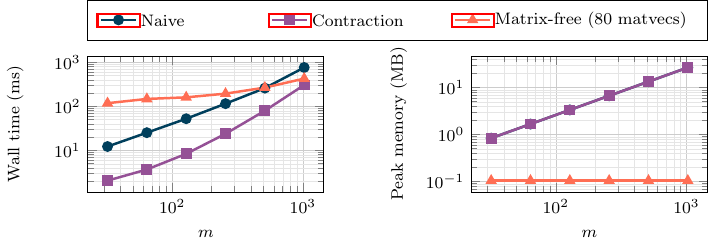}
\caption{\newtext{Cost of computing the empirical NTK for a four-layer tanh-MLP of width 32 ($n=3{,}265$ parameters) as a function of the number of inputs $m$. 
\emph{Left:} Wall time. \emph{Right:} Dominant tensor memory.}}
\label{fig:ntk_scaling}
\end{figure}
\end{num_example}

\begin{takeawaybox}
\newtext{
The NTK is a fundamental tool for understanding the difficulties related to the solution of problems involving parametric models. 
It connects parameter-space optimization to function-space learning dynamics: the eigenvalues of $\Theta_k$ govern how quickly each error mode contracts during training.
An ill-conditioned NTK, caused by widely separated eigenvalues, leads to slow convergence, as low-frequency modes learn much faster than high-frequency ones (\emph{spectral bias}).
Adaptive and second-order methods can improve convergence by reshaping the spectrum of the preconditioned NTK kernel $\Theta_k^{(M_k)}$, reducing its condition number and accelerating training.
}
\end{takeawaybox}

  \section{Stochastic gradient descent (SGD)}
\label{sec:first_order}
Stochastic Gradient Descent (SGD) is the most widely used optimization method in modern ML, and plays an equally central role in data-driven SciML.
Its popularity stems from its simplicity, scalability, and ability to train high-dimensional models on massive datasets.
This widespread use is largely due to the limitations of classical GD, which, when applied to the finite-sum problem~\eqref{eq:min_problem}, requires evaluating the full gradient:
\[
\nabla f(\theta)
  = \frac{1}{m} \sum_{i=1}^m \nabla f_i(\theta).
\]
When the dataset size ($m$) is large, evaluating the full gradient quickly becomes computationally prohibitive.
However, real-world datasets are highly redundant, and computing every per-sample gradient at each iteration usually provides little additional benefit.

These considerations motivate the shift from classical deterministic optimization to \emph{stochastic optimization}.  
The finite-sum structure of ML problems allows one to replace the full gradient with a cheaper approximation.
The stochastic gradient (SGD)
 method randomly chooses a sample  at each iteration and utilizes just the corresponding gradient.  
Its iterations thus read as follows: 
\begin{equation}\label{eq:sgd}
\xkk=\xk-\alpha_k\nabla f_{i_k}(\xk),
 \end{equation}
 where $\alpha_k>0$ is a given step size. The symbol $i_k$ denotes an index drawn uniformly at random from $\{1,\dots,m\}$ at  iteration $k$.

 Each iteration of the SGD method is thus very cheap, involving only the computation of the gradient $\nabla f_{i_k}(\theta_k)$, corresponding to one sample, \newtext{rather than all the $m$ terms}. 
Since the index $i_k$ is chosen randomly, the generated sequence  depends on the random sequence $\{i_k\}$ and is not uniquely  determined from $\theta_0$.
Note that each direction $-\nabla f_{i_k}(\theta_k)$ might not be one of descent for $f$ (in the sense of
yielding a negative directional derivative), but we can show that if it is a descent direction
in expectation (the expected value of the scalar product with $\nabla f(\theta_k)$ is negative), then the sequence $\{\xk\}$ can be guided toward a minimizer of $f$. 
This comes from the fact that $\nabla f_{i_k}(\theta_k)$ is an unbiased estimator of $\nabla f(\theta_k)$, i.e., 
$$
\mathbb{E}[\nabla f_{i_k}(\theta_k)\mid \theta_k] = \sum_{i=1}^m
\frac{1}{m}\nabla f_{i}(\theta_k)=\nabla f(\theta_k).$$

The SGD method is summarized in Algorithm~\ref{alg:sgd}.
Note that progress is in practice measured in epochs rather than  iterations, where one epoch represents a complete pass over the  training dataset, corresponding to $m$ iterations.
\newtext{In contrast, the standard GD method, which evaluates the  gradient using all $m$ samples, performs exactly one iteration per  epoch, which is why it is often referred to as batch gradient descent.}

\begin{algorithm}[H]
\caption{Stochastic Gradient (SGD) Method}
\label{alg:sgd}
\begin{algorithmic}[1]
\State \textbf{Given:} A dataset $\{x_i,y_i\}_{i=1}^m$, the objective functions $\{f_i\}_{i=1}^m$, an initial iterate $\theta_0 \in \R^n$
\For{$k = 0, 1, 2, \dots$}
    \State Randomly choose $i_k \in \{1, \dots, m\}$ \COMMENTmine{Sample-index selection}
    \State $p_k = -\nabla f_{i_k}(\theta_k)$ \COMMENTmine{Search direction computation}
    \State Choose $\alpha_k > 0$ \COMMENTmine{Step size selection}
    \State $\theta_{k+1} = \theta_k + \alpha_k p_k$ \COMMENTmine{Iterate update}
\EndFor
\end{algorithmic}
\end{algorithm}

The SGD method often exhibits rapid progress during the early stages of training but may stagnate as training proceeds. \newtext{This behavior is commonly attributed to the inherent noise in stochastic gradient estimates, typically measured in terms of the variance
\begin{equation}\label{eq:variance}
\sigma^2(\theta)
:=
\mathbb{E}\!\left[
\|\nabla f_{i_k}(\theta)-\nabla f(\theta)\|^2
\right].
\end{equation}

This noise can help the iterates escape local minima and saddle points, but it also limits the attainable optimization accuracy, as discussed below. However, since ERM only approximates the expected risk, excessively optimizing the empirical objective may be detrimental and lead to overfitting, i.e., degraded performance on unseen data. Consequently, this limitation is often not critical in practice.}

\emph{Mini-batch methods} provide an intermediate approach between full-gradient methods and SGD.
Whereas SGD relies on a single sample at each iteration, mini-batch methods compute gradient estimates using a small subset $ \pazocal{I}_k\subset \{1,\dots,m\}$ of the samples at each iteration. 
\newtext{The mini-batch $\pazocal{I}_k$ is typically drawn uniformly  without replacement at each epoch,} giving rise to the following  update rule:
\begin{align*}
\xkk=\xk-\frac{\alpha_k}{\lvert \pazocal{I}_k \rvert}\sum_{i\in \pazocal{I}_k}\nabla f_{i}(\xk).
\end{align*}
This allows one to exploit some parallelism when computing mini-batch gradients, as each sample contribution can be evaluated independently.
Moreover, a straightforward generalization of the variance  formula~\eqref{eq:variance} to the mini-batch setting shows that  the variance of the gradient estimate decreases as  $\lvert\pazocal{I}_k\rvert$ increases.
As a consequence, mini-batch methods are easier to tune in terms of  step size selection.
\newtext{Indeed, as we will see in the following section, the product of the step size and of the variance impacts how closely the iterates  can approach an optimum.
Since the variance decreases with larger mini-batch sizes, larger step sizes can be used, yielding faster convergence with lower optimization error.

}

In practice, selecting an appropriate step size is challenging, as its value directly affects the method's convergence.
The stepsize~$\alpha_k$ is often chosen heuristically. 
\newtext{For example, a small constant value may be fixed, which, however, can ensure just convergence to a neighborhood of an optimum.
Most often, a step size schedule is used to decrease $\alpha_k$ progressively, see~Numerical Example \ref{sec:example_code}.}

\subsection{Theoretical results}
\label{sec:theory_sgd}
In this section, we present key convergence results for SGD applied to problem~\eqref{eq:min_problem}, and outline one representative proof to illustrate the main techniques.
Additional results and proofs can be found in~\cite{garrigos2023handbook,bottou2018optimization}.
Throughout, we denote by $\theta^\star$ a solution of~\eqref{eq:min_problem}.
Note, the discussed results remain valid if the stochastic gradient $\nabla f_{i_k}$ is replaced by any unbiased gradient estimator $g_k$ satisfying $\mathbb{E}[g_k\mid \theta_k]=\nabla f(\theta_k)$.

The stated results mainly depend on the assumptions we make about the objective function. We always assume the function to be bounded from below. 
One quite common assumption is smoothness.

\begin{defn}
 A function $f:\mathbb{R}^n\rightarrow \mathbb{R}$ is said to be $L$-smooth if there exists $L>0$ such that, for all $\theta_1,\theta_2\in \mathbb{R}^n$, 
 $$
 \|\nabla f(\theta_1)-\nabla f(\theta_2)\|\leq L \|\theta_1-\theta_2\|. 
 $$
\end{defn}

\begin{ass}\label{hp_smooth}
Each function $f_i:\mathbb{R}^n\rightarrow \mathbb{R}$ in problem \eqref{eq:min_problem} is $L_i$-smooth. 
We denote $L_{\max}:=\max_{i=1,\dots,m} L_i$.
\end{ass}

Another important assumption, which, however, restricts the class of functions under study, is the convexity assumption. \newtext{This assumption generally does not hold in the most general neural network training problem, but there are contexts in which it plays an important role, see for instance \cite{bengio2005convex,bach2017breaking}.}
\begin{defn}\label{def:convex}
A function $f:\mathbb{R}^n\rightarrow \mathbb{R}$ is said to be convex if for all $\theta_1,\theta_2\in \mathbb{R}^n$, for all $t\in [0,1]$
$$
f(t\theta_1+(1-t)\theta_2)\leq tf(\theta_1)+(1-t)f(\theta_2).
$$
$f$ is said to be \emph{strictly} convex if for all $\theta_1,\theta_2\in \mathbb{R}^n$, for all $t\in (0,1)$
$$
f(t\theta_1+(1-t)\theta_2)< tf(\theta_1)+(1-t)f(\theta_2).
$$
\end{defn}
\begin{defn}\label{def:strongly_convex}
A function $f:\mathbb{R}^n\rightarrow \mathbb{R}$ is said to be $\mu$-strongly convex if there exists $\mu>0$ such that for all $\theta_1,\theta_2\in \mathbb{R}^n$
$$
f(\theta_2)\geq f(\theta_1)+\nabla f(\theta_1)^T (\theta_2-\theta_1)+\frac{\mu}{2}\|\theta_1-\theta_2\|^2.
$$
\end{defn}
When needed, we will assume the following. 
\begin{ass}\label{hp_convex}
Each function  $f_i:\mathbb{R}^n\rightarrow \mathbb{R}$ in  problem \eqref{eq:min_problem} is convex. 
\end{ass}

Convexity of the function usually makes the solution of the problem easier, since the landscape will be more benign (for instance, all the local minima are global minima for convex functions). With this assumption, we can obtain stronger convergence results. For instance, we will see that when $f$ is convex, we can prove a convergence result on the sequence of functions (even on the sequence of iterates when $f$ is strongly convex), while in the nonconvex case, we can just show a result on the norm of the gradient. 

Another important factor influencing convergence is the variance of the gradient estimate.
\begin{ass}\label{hp_noise}
    Assume that the variance of the gradient estimate in \eqref{eq:variance} is bounded:  $\sigma^2(\theta)\leq \sigma^2$ for all $\theta$. 
\end{ass}
We will see that a non-zero variance will negatively impact the convergence of optimization algorithms, which is why variance reduction techniques~\cite{bottou2018optimization} may be used to obtain better stochastic directions and relax assumptions on the step size, which also affects the results.  
In this chapter, we consider both constant ($\alpha_k=\alpha$ for all $k$) and vanishing ($\alpha_k=\pazocal{O}(1/k)$) step sizes.

We will see that in all the results, the key constants that determine the convergence rates
(smoothness $L_{\max}$, strong convexity $\mu$, and the SGD variance)
are directly related to the eigenvalues of the Hessian $H(\theta)$.
In particular, we recall that, if $f$ is $\mu$-strongly convex, then
\[
L_{\max} =\sup_{\theta}\lambda_{\max}(H(\theta)), 
\quad 
\mu=\inf_{\theta} \lambda_{\min}(H(\theta)),
\quad
\kappa(H(\theta))=\frac{\lambda_{\max}(H(\theta))}{\lambda_{\min}(H(\theta))}.
\]
Thus, the magnitudes of $\lambda_{\min}(H(\theta))$ and $\lambda_{\max}(H(\theta))$, as well as their ratio $\kappa(H(\theta))$, 
play a crucial role in interpreting the obtained theoretical convergence guarantees.
The convergence theorems show in particular how SGD deteriorates with increasing $\kappa(H(\theta))$.

\begin{figure}
  \centering
  \includegraphics{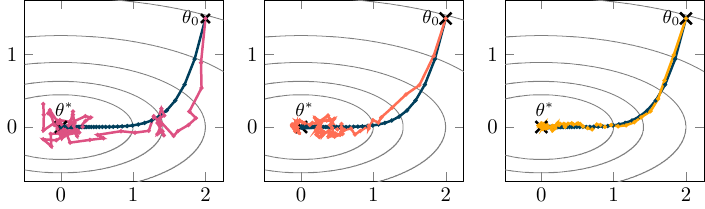}
\caption{Trajectories of GD (blue) and SGD on a smooth quadratic function for different noise levels 
$\sigma \in \{\sigma_0,\;\sigma_0/2,\;\sigma_0/4\}$ (purple, orange, yellow). 
Higher noise magnitudes enlarge the region in which SGD 
oscillates around the minimizer $\theta^\ast$, while lower noise confines the iterates 
closer to the GD path.}
  \label{fig:SGD_threepanel}
\end{figure}

\subsubsection{Convergence for strongly convex and smooth functions}
The first result applies to strongly convex functions, i.e.,~the ones with the most favorable landscape. 

\begin{theorem}\label{thm:strongly}
Let assumptions  \ref{hp_smooth}, \ref{hp_convex}  and \ref{hp_noise} hold, and assume that  $f$ is $\mu$-strongly convex. Assume that $\{\theta_k\}_{k\in\mathbb{N}}$ is the sequence generated by SGD \eqref{eq:sgd}  
 with a constant step size $0<\alpha_k=\alpha< \frac{1}{2 L_{\max}}$
for all $k$. Then, for $k\geq 1$,
one has 
$$
\E{[\|\theta_k-\theta^{\star}\|^2]}\leq (1-\alpha \mu)^k \|\theta_0-\theta^{\star}\|^2+\frac{2\alpha}{\mu } \sigma^2.
$$
\end{theorem}

\begin{proof}
From the definition of the step \eqref{eq:sgd} it holds that
\begin{align*}
\|\theta_{k+1}-\theta^\star\|^2
&= \|\theta_k-\alpha \nabla f_{i_k}(\theta_k) - \theta^\star\|^2\\
&= \|\theta_k-\theta^\star\|^2
 -2\alpha \langle \nabla f_{i_k}(\theta_k), \theta_k-\theta^\star\rangle
 + \alpha^2 \|\nabla f_{i_k}(\theta_k)\|^2 .
\end{align*}
Taking the conditional expectation given $\theta_k$ yields
\begin{align*}
\mathbb{E}[\|\theta_{k+1}-\theta^\star\|^2\mid \theta_k]
&= \|\theta_k-\theta^\star\|^2
   - 2\alpha\langle \nabla f(\theta_k), \theta_k-\theta^\star\rangle
   + \alpha^2 \mathbb{E}[\|\nabla f_{i_k}(\theta_k)\|^2\mid \theta_k].
\end{align*}
From the definition of strong convexity it follows 
\[
\mathbb{E}[\|\theta_{k+1}-\theta^\star\|^2\mid \theta_k]
\le (1 - \alpha \mu) -2\alpha (f(\theta_k)-f(\theta^\star))
  + \alpha_k^2 \sigma^2.
\]
Taking full expectation and  using a variance transfer (see \cite[lemma 4.20]{garrigos2023handbook})  gives the main recursion:
\begin{align*}
\mathbb{E}[\|\theta_{k+1}-\theta^\star\|^2]
&\le (1 - \alpha_k \mu) \,\mathbb{E}[\|\theta_k-\theta^\star\|^2]
 + 2\alpha_k^2 \sigma^2 +2\alpha(2\alpha L_{\max}-1)\mathbb{E}[f(\theta_k)-f(\theta^\star)]\\
&\le (1 - \alpha_k \mu) \,\mathbb{E}[\|\theta_k-\theta^\star\|^2]
 + 2\alpha_k^2 \sigma^2,
\end{align*}
where we have used the assumption on the step size. 
Recursively applying the
above and summing up the resulting geometric series gives
\begin{align*}
\mathbb{E}[\|\theta_{k+1}-\theta^\star\|^2]&\leq (1-\alpha\mu)^k\|\theta_0-\theta^\star\|^2+2\sum_{j=0}^{k-1}(1-\alpha\mu)^j\alpha^2\sigma^2\\
&\leq(1-\alpha\mu)^k\|\theta_0-\theta^\star\|^2+\frac{2\alpha\sigma^2}{\mu}.
\end{align*}

\end{proof}

\newtext{
This result tells us that SGD is not guaranteed to converge with a constant step size, since the upper bound on the error does not go to zero when $k$ increases. 
Instead, the error stabilizes on the constant value $\frac{2\alpha\sigma^2}{\mu}$}.  
If the variance of the gradient estimate were zero, the generated sequence would converge in expectation to the solution at a fast linear rate (i.e., the error from one iterate to the next would decrease exponentially). 
A gradient estimate with nonzero variance leads to the addition of a constant positive term on the right-hand side, which breaks the monotone decrease of the error and leads the sequence to converge to a ball around 
the optimum, rather than exactly to the optimum itself. 

This means that the method exhibits linear convergence, similar to classical GD, at a rate depending on the learning rate $\alpha$.
This phase continues until the iterates enter a neighborhood of the optimum of radius proportional to $\alpha \sigma^2$; see  Figure~\ref{fig:SGD_threepanel} for an illustration.
In addition, a large $\alpha$ leads to a fast initial decrease but to a larger steady-state error, while a small value leads to a slighter decrease, but we can eventually come closer to the minimum. 

Note that the eigenvalues of the Hessian also impact the convergence. 
A small $\mu\approx \lambda_{\min}(H(\theta_k))$ makes the linear rate $(1-\alpha\mu)$ close to~$1$, while it simultaneously amplifies the variance term $(2\alpha/\mu)\sigma^2$.
Thus, SGD cannot make progress in flat directions, where the curvature is very small 
($\lambda_{\min}(H(\theta_k))\approx 0$), and instead accumulates noise there.

\subsubsection{Convergence for convex and smooth functions}
In the following theorem we state an ergodic  result for SG, applied to convex functions, with fixed step size \cite{garrigos2023handbook}. 
\begin{theorem}
Let assumptions  \ref{hp_smooth}, \ref{hp_convex} and \ref{hp_noise} be satisfied. Assume that $\{\theta_k\}_{k\in\mathbb{N}}$ is the sequence generated by SGD \eqref{eq:sgd} with constant step size $0<\alpha_k=\alpha\leq \frac{1}{4 L_{\max}}$, for all $k$. Define, for every $K\geq 1$, 
$$\bar \theta_K:= \frac{1}{K}\sum_{k=1}^K\theta_k.$$ 
Then, 
\begin{equation}\label{eq:convex}
\E{[f(\bar{\theta}_K)-f^{\star}]}\leq \frac{\|\theta_0-\theta^{\star}\|^2}{\alpha K}+2\alpha \sigma^2,
\end{equation}
\newtext{where $f^\star=\inf f$.} In particular, if for a fixed $K$ we set $\alpha=\frac{\alpha_0}{\sqrt{K}}$ with $\alpha_0\leq \frac{1}{4 L_{\max}}$, the two terms are balanced, giving 
$$
\E{[f(\bar{\theta}_K)-f^{\star}]}=\frac{\|\theta_0-\theta^{\star}\|^2}{\alpha_0 \sqrt{K}}+\frac{2\alpha_0 \sigma^2}{\sqrt{K}}=\pazocal{O}\left(\frac{1}{\sqrt{K}}\right).
$$
\end{theorem}

We can see that the decrease is slower than in the strongly convex case (it is sublinear) and that the variance has the same effect as for strongly convex functions. In this case, the largest eigenvalue of the Hessian plays a role. Since the step size $\alpha$ is limited by \newtext{$\alpha \le \frac{1}{4L_{\max}}\le \frac{1}{4\lambda_{\max}(H(\theta_k))}$}, 
a large $\lambda_{\max}(H(\theta_k)) $ (and hence a large condition number) forces very small step sizes.  
Moreover, reducing $\alpha$ also increases the $1/(\alpha K)$ term in \eqref{eq:convex}, 
delaying convergence.

Once again, convergence cannot be achieved with constant step sizes, but the effect of the variance can be counteracted by a vanishing step size, as shown in the following theorem \cite{garrigos2023handbook}. 

\begin{theorem}
Let assumptions  \ref{hp_smooth}, \ref{hp_convex} and \ref{hp_noise} be satisfied. Assume that $\{\theta_k\}_{k\in\mathbb{N}}$ is the sequence generated by SGD \eqref{eq:sgd} with a vanishing  step size $\alpha_k=\alpha_0/\sqrt{k+1}$ with $\alpha_0\leq \frac{1}{4 L_{\max}}$, and define 
$$
\bar{\theta}_K:=\frac{1}{\sum_{k=0}^{K-1}\alpha_k}\sum_{k=0}^{K-1}\alpha_k \theta_k.
$$
Then, \newtext{for $f^\star=\inf f$,}

$$
\E{[f(\bar{\theta}_K)-f^{\star}]}\leq \frac{5\|\theta_0-\theta^{\star}\|^2}{4\alpha_0 \sqrt{K}}+ \sigma^2\frac{5\alpha_0\log(K+1)}{\sqrt{K}}=\pazocal{O}\left(\frac{\log(K+1)}{\sqrt{K}}\right). 
$$
\label{th:convex}
\end{theorem}

\subsubsection{Convergence for nonconvex and smooth functions}

In this section, we just assume that the stochastic gradients are Lipschitz continuous. Note that in this nonconvex setting,  we cannot prove global optimality results.  Nevertheless, we can still obtain bounds on the stationarity of the algorithm, stated in the following theorem \cite{bottou2018optimization}. 
\begin{theorem}
Let Assumptions \ref{hp_smooth} and \ref{hp_noise} hold. Assume that $\{\theta_k\}_{k\in\mathbb{N}}$ is the sequence generated by SGD \eqref{eq:sgd} with a constant  step size 
$0<\alpha_k=\alpha\le \frac{1}{L_f}$, 
where $L_f$ is the Lipschitz constant of the gradient of $f$.  
Then, for any  $K\ge1$,
\[
\frac1K
\sum_{k=0}^{K-1}
\mathbb E[\|\nabla f(\theta_k)\|^2]
\le
\frac{2(f(\theta_0)-f^\star)}{\alpha K}
+
L_f\alpha\sigma^2,
\]
where $f^\star=\inf f$. In particular, choosing
$
\alpha
=\min\{\frac{1}{L_f},
\frac{
1
}{\sqrt{
K
}}\},
$
yields
\[
\frac1K
\sum_{k=0}^{K-1}
\mathbb E[\|\nabla f(\theta_k)\|^2]
\leq\frac{2(f(\theta_0)-f^\star)}{\sqrt{K}}
+
\frac{L_f\sigma^2}{\sqrt{K}}=
\pazocal O\!\left(\frac{1}{\sqrt K}\right).
\]
Consequently, for a given $\epsilon>0$, if $K=\pazocal{O}(\epsilon^{-2})$ then 
$$
\frac1K
\sum_{k=0}^{K-1}
\mathbb E[\|\nabla f(\theta_k)\|^2]=\pazocal{O}(\epsilon).
$$
\label{th:nonconvex}
\end{theorem}

This means that we need $\pazocal{O}(\epsilon^{-2})$ iterations to drive the mean gradient norm below the tolerance $\epsilon$. 
Note that the constant hidden in the $\pazocal{O}(1/\sqrt{K})$ term depends on $L_f=\sup_\theta  |\lambda_{\max}(H(\theta))|$.
Thus, a large curvature in some directions (large $L_f$) makes
$\pazocal{O}(1/\sqrt{K})$ progress much slower in practice.

\newtext{
As we have seen, the Hessian's eigenvalues critically affect the theoretical results across all considered settings.
Notably, small eigenvalues $\lambda_{\min}(H(\theta))$ lead to slow contraction in the convex setting, while the variance scales as  $1/\lambda_{\min}(H(\theta))$, amplifying noise as $\lambda_{\min}(H(\theta))$  decreases. 
At the same time, large eigenvalues enforce very small step sizes.}
These observations explain why first-order methods, such as (S)GD, struggle when training ill-conditioned problems, such as PINNs.

\subsection{\newtext{Comparison of computational cost: GD vs.\ SGD}}
\newtext{Since the iterations of SGD are much cheaper than those of classic GD, one may expect that it is always preferable to use SGD. 
However, the inexactness introduced in the gradient evaluation comes at the cost of a slower convergence rate. 
For instance, if $f$ is $\mu$-strongly convex, GD is known to converge linearly~\cite[Theorem 10.29]{beck2017first}, meaning that}
\newtext{
	\begin{equation*}
	f(\xk)-f^{\star}\leq \left(1-\frac{\mu}{L}\right)^k(f(\theta_0)-f^{\star}).
	\end{equation*}  
	The total number of iterations to drive the training error below a given threshold $\epsilon>0$ is thus proportional to $\log(1/\epsilon)$. 
	Since the per-iteration cost is proportional to $m$ (due to the need to compute the full gradient), the total work required to obtain the accuracy $\epsilon$ is proportional to $m\log(1/\epsilon)$.

The convergence rate of the SGD method is slower than that of the GD method.   
For instance, in the strongly convex case with $i_k$ drawn uniformly  from $\{1,\dots,m\}$ and step size $\alpha_k = \pazocal{O}(1/k)$,  SGD converges sublinearly in expectation~\cite[Theorem 4.7]{bottou2018optimization}:
	\begin{equation*}
	\mathbb{E}(f(\xk)-f^{\star})=O(1/k).
	\end{equation*}
However, the per-iteration cost of SGD does not depend on the dataset size $m$, making it attractive in large-scale settings.
Specifically, the total work to reach accuracy $\epsilon$ is proportional to $1/\epsilon$ for SGD, compared to $m\log(1/\epsilon)$  for GD.
Thus, while SGD may be less efficient for small to moderate $m$, it becomes increasingly advantageous as $m$ grows large.
}
	
\newtext{
If the gradient estimate is computed using a mini-batch of data,  a different trade-off between cost and convergence rate arises.
The per-iteration cost scales with the mini-batch size, while the variance of the gradient estimate decreases relative to plain SGD,  leading to improved convergence.
Thus, the main difference between these methods lies in step quality:  GD performs fewer but more accurate and more expensive steps, whereas 
SGD takes many cheap but noisy steps, and mini-batch methods provide a middle ground.
}

\newtext{
The same comparison applies to memory requirements. 
In particular, both GD and SGD require $\pazocal O(n)$ memory for the parameters and the gradient estimate. 
Storing the full dataset requires $\pazocal{O}(md)$ memory, where  $d$ is the per-sample feature dimension.
In the mini-batch setting, only the current $|\pazocal I_k|$ samples need to reside in memory, which is the practical reason why mini-batching is the default strategy once $m$ samples no longer fit in GPU memory.
}

\begin{num_example}{\textbf{\emph{SGD implementation.}}}
\label{sec:example_code}
This numerical example demonstrates how the theoretical results discussed in Section~\ref{sec:theory_sgd} manifest in practice.  
To illustrate the behaviour of SGD, we consider the logistic regression problem introduced in Section~\ref{sec:examples} with a synthetic dataset.  
Figure~\ref{fig:variance_SGD} shows the effect of the mini-batch size on the behaviour of the SGD method.  
With a fixed learning rate and very small batches, SGD initially makes rapid progress but soon stagnates as the stochastic variance in the gradient estimates becomes dominant, preventing a consistent decrease in the loss.  
In contrast, larger batches reduce the sampling noise and produce search directions  that more closely approximate the full gradient. 
\newtext{A prototypical implementation of the SGD in NumPy can be found in Code snippet~\ref{lst:SGD}.}

\begin{lstlisting}[language=Python, caption={\newtext{Example of SGD implementation in NumPy.}},label=lst:SGD]
def run_sgd(X, y, gamma=1e-4, epochs=300, batch_size=8,
    alpha0=1e-3, scheduler=None, seed=0, f_opt=None):
   # Run SGD with a fixed batch size and (optional) learning-rate scheduler
   
    n, d = X.shape # Number of sample (n) and feature dimension (d)
    rng = np.random.default_rng(seed) # Set seed for reproducibility
    w = np.zeros(d) # Initial guess for parameters

    for ep in range(epochs): # Run over epochs
    	  #  Use learning-rate scheduler. If None, fixed at alpha0
        lr = alpha0 if scheduler is None else scheduler(ep, alpha0)

        idx = rng.permutation(n)
        Xs, ys = X[idx], y[idx] # Suffle samples
        
        for start in range(0, n, batch_size): # Run over all mini-batches
            end = min(start + batch_size, n) # Get mini-batch
            # SGD step direction
            d_step = -lr * grad_logloss(w, Xs[start:end], ys[start:end], gamma)
            w = w + d_step  # Update the parameters
        
\end{lstlisting}

\begin{figure}
\centering
  \includegraphics{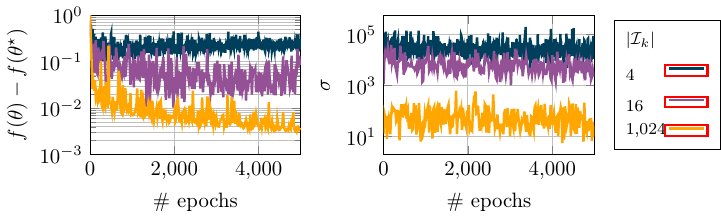}
\caption{Effect of batch size on convergence of the mini-batch GD method with fixed $\alpha_k=10^{-3}$.
\emph{Left:} Suboptimality $f(\theta_k)-f(\theta^\star)$ over epochs for three batch sizes ($| \pazocal{I} | = \{4; 16; 1,024 \}$).
\emph{Right:} Corresponding relative variance of subsampled gradient.}
\label{fig:variance_SGD}
\end{figure}

In practice, two standard mechanisms can be used to \newtext{mitigate the effect of the variance of the} stochastic search directions, and to \newtext{drive the residual variance-related term in the convergence results to zero.}
Their effect is illustrated in Figure~\ref{fig:schedulers}.  
The first mechanism is \emph{learning rate scheduling}, where the step size $\alpha_k$ is incrementally reduced
during training, \newtext{see the first function in Code snippet \ref{lst:LRS_BSS} for a possible implementation.}
In our experiments, we employ a \emph{polynomial decay} scheduler of the following form:
\[
\alpha_k \;=\; \frac{\alpha_0}{(1 + \tau k)^{pw}},
\]
with initial learning rate $\alpha_0 = 10^{-3}$, decay coefficient $\tau = 2\!\cdot\!10^{-2}$, and exponent $pw = 1$.  
Reducing the learning rate decreases the magnitude of stochastic perturbations and,
as predicted theoretically, it enforces a gradual reduction of variance~\cite{robbins1951stochastic,bottou2018optimization}.  
Note that other schedules, such as step decay, exponential decay, or cosine annealing, are commonly used in modern practice~\cite{loshchilov2017sgdr}.

\begin{lstlisting}[language=Python, caption={\newtext{Example of implementation of a polynomial learning-rate scheduler and a progressive batch-size scheduler.}},label=lst:LRS_BSS]
def poly_scheduler(epoch, alpha0, decay=5e-3, power=1):
    # Polynomial decay learning-rate schedule
    return alpha0 / ((1 + decay * epoch) ** power)
    
def batch_size_scheduler(epoch, bs_start, bs_end, frac):
	 # Linear batch scheduler from bs_start to bs_end at frac rate
        frac = epoch / (epochs - 1) if epochs > 1 else 1.0
        bs = int(bs_start + frac * (bs_end - bs_start))
        return max(1, min(bs, n))
\end{lstlisting}

\newtext{The second mechanism is intended to reduce the variance of the gradient estimates by a \emph{batch-size scheduler}, which increases the accuracy of each gradient estimate by progressively increasing the size of the batches, \newtext{see the second function provided in Code snippet \ref{lst:LRS_BSS} for a possible implementation.} 
In this way, early steps are cheap, computed on a few samples, and become increasingly accurate but more expensive as the iterations approach a minimum, allowing for a more precise final solution.
} 
In our experiment, the mini-batch size evolves linearly from 
$|\pazocal{I}_0| = 16$ to $|\pazocal{I}_T| = 1,024$ over $T$ epochs as
\[
|\pazocal{I}_k|
\;=\;
|\pazocal{I}_0|
\;+\;
\frac{k}{T}\bigl(|\pazocal{I}_T| - |\pazocal{I}_0|\bigr),
\qquad k = 0,\dots,T.
\]
The learning rate is kept constant at $\alpha_k = \alpha_0$, enabling progress in regimes where small-batch SGD with a fixed step size would stagnate~\cite{smith2017don,shallue2019measuring}.

Figure~\ref{fig:schedulers} demonstrates the performance of these two approaches, compared to SGD with $\alpha_k=10^{-3}$ and $| \pazocal{I}_k | = 16$.
As we can see, both learning rate decay and batch-size growth strategies effectively reduce the variance of the search direction, thereby restoring convergence.  
However, they differ subtly in practice: learning rate decay improves stability but slows down the method permanently, whereas batch-size growth preserves the step magnitude and reduces noise without sacrificing asymptotic speed, but progressively increases the iteration cost.

\begin{figure}
\centering
\includegraphics{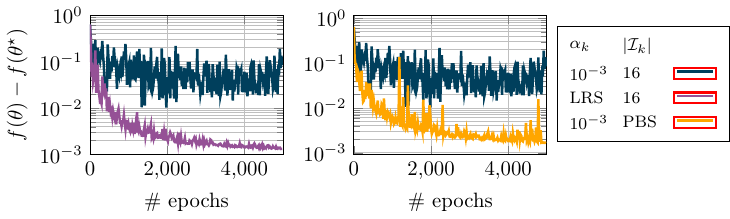}
\caption{Example of reducing SGD noise by using learning rate scheduling (LRS) and progressive batch-size (PBS) increase.
The SGD (blue) utilizes ${|\pazocal{I}_k|=16}$ and $\alpha_k=10^{-3}$. 
The SGD with LRS (purple) utilizes polynomial decay of~$\alpha_k$, while SGD with PBS progressively increases the batch size $|\pazocal{I}_k|$.}
\label{fig:schedulers}
\end{figure}
\end{num_example}

\begin{takeawaybox}
\newtext{
Stochastic gradient descent (SGD) and its mini-batch variants minimize finite-sum objectives by replacing the full gradient with a cheap estimate computed on small  random subsets of data.
The \emph{variance} of this estimate governs convergence. 
In particular, large variance  limits the attainable accuracy with fixed step sizes, while small variance  permits larger steps and faster progress.
Two standard strategies to counteract the effect of variance are  \emph{step size decay} and \emph{batch size growth}, both of which restore convergence in expectation, at the cost of slower iterations  or higher per-iteration cost, respectively.
}
\end{takeawaybox}

 \section{Accelerated and  adaptive stochastic gradient methods}
\label{sec:adaptive_methods}
The main limitation of standard first-order methods is their slow convergence, often attributed to the difficulty of selecting an appropriate step size.
This issue becomes even more pronounced in large-scale SciML, where physical constraints, large models, noisy gradients, and heterogeneous data make choosing a suitable step size both challenging and essential for stable and efficient training.
To address these challenges, a wide range of accelerated and adaptive variants of GD and SGD have been developed.
In this section, we present the most widely used ones.
\newtext{From the unifying preconditioning viewpoint introduced in Section~\ref{sec:precond_lens}, the adaptive gradient methods discussed below act as diagonal preconditioners in parameter space.}

\begin{figure}
\centering
  \includegraphics{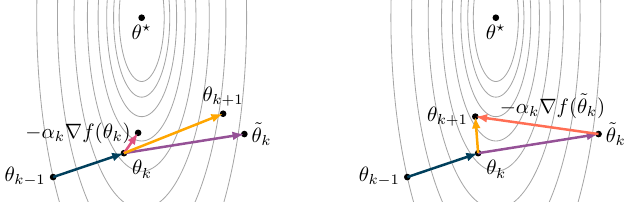}
\caption{Sketch of heavy-ball momentum (left) and NAG (right) iteration.}
\label{fig:momentum_sketch}
\end{figure}

\subsection{Momentum}
The purpose of momentum is to mitigate a key limitation of GD, namely that the parameter update at iteration $k$ only depends on local information on the gradient, and does not exploit the history of previous updates.
Incorporating past information can provide valuable insight into the objective function landscape along the optimization trajectory and can significantly improve convergence.

Momentum techniques keep memory of past gradients and accumulate them. 
To account for time, they use an exponentially weighted average to give more weight to recent gradients. 
The update at iteration $k$ then becomes
\[
p_k = \beta\, p_{k-1} + (1-\beta)\, g_k,
\]
where $p_k$ is the new search direction and $g_k$ is any gradient approximation at iteration $k$. 
The momentum parameter  $\beta \in (0,1)$ controls the weighting between past and current gradients. 
\newtext{The gradient approximation $g_k$ is typically obtained by sub-sampling the training set, and may correspond to the full gradient ($|\pazocal{I}_k|=m$), a mini-batch gradient  ($1<|\pazocal{I}_k|<m$), or a stochastic gradient  ($|\pazocal{I}_k|=1$).}

The term $(1-\beta)$ is often replaced with the learning rate, and a dynamically adjusted parameter $\beta_k$ may be used instead of a fixed $\beta$ to scale $p_{k-1}$, so the general update  is given as
\begin{equation}
\label{momentum}
\theta_{k+1}=\theta_{k}-\alpha_k g_k + \beta_k(\theta_{k} - \theta_{k-1}),
\end{equation}
where the scalar sequence $\{\alpha_k\}$ is either predetermined or set dynamically. 

\paragraph{Heavy ball}  If $\alpha_k = \alpha$ and $\beta_k = \beta$ for some constants $\alpha > 0$ and $\beta \in (0,1)$ and for all $k \in \mathbb{N}$, the update~\eqref{momentum} is referred to as the \emph{heavy ball method}.
An alternative interpretation of this method is obtained by expanding the update:
\[
\theta_{k+1}
=
\theta_k
-
\alpha \sum_{j=1}^{k} \beta^{k-j} g_j.
\]
Thus, each step can be viewed as an exponentially weighted average of past gradients.
Notably, larger values of $\beta$ place greater emphasis on gradients from earlier iterations.
In practice, $\beta$ is typically chosen close to $0.9$.

\paragraph{Nesterov accelerated gradient (NAG)}
Another widely used variant of momentum is the Nesterov accelerated gradient (NAG). 
The core idea behind NAG is that if the current parameter vector is $\theta_{k}$, then the momentum term alone (i.e., ignoring the term with the gradient) is about to nudge the parameter vector by $\beta_k (\theta_{k}-\theta_{k-1})$. 
Therefore,  it makes sense to compute the gradient at $\theta_{k} + \beta_k (\theta_{k}-\theta_{k-1})$, instead of at the old position $\theta_{k}$.

Written as a two-step procedure, the NAG update is given as
\begin{equation*}
\tilde{\theta}_k=\theta_{k} + \beta_k ( \theta_{k} - \theta_{k-1} ) \qquad \text{and} \qquad \theta_{k+1}=\tilde{\theta}_k-\alpha_k\nabla f(\tilde{\theta}_k),
\end{equation*}
which can be reformulated as
\begin{equation}\label{nesterov}
\theta_{k+1}= \theta_{k}-\alpha_k\nabla f(\theta_{k} + \beta_k(\theta_{k} - \theta_{k-1})) + \beta_k(\theta_{k} -\theta_{k-1}).
\end{equation}
The main difference with classical momentum is thus the order of computation (momentum step versus gradient step), see also Figure~\ref{fig:momentum_sketch}.

\paragraph{Cost and convergence}
\newtext{The cost of (S)GD methods with momentum is comparable to that of classical (S)GD, but their convergence is generally faster. For example, in the deterministic strongly convex case, the convergence rate of GD (cf.~Theorem \ref{thm:strongly} with $\sigma=0$) improves from $\pazocal{O}((1-\frac{\mu}{L})^k)$ to $\pazocal{O}((1- \sqrt{\frac{\mu}{L}})^k)$, see\cite{polyak1964some} for details.
In the stochastic case, the convergence neighborhood size increases from $\pazocal{O}(\alpha \sigma^2)$ (cf.~Theorem \ref{thm:strongly}) to $\pazocal{O}(\frac{\alpha \sigma^2}{1-\beta})$, where $\alpha$ and $\beta$ denote the step size and the momentum parameter, respectively. 
Thus, momentum improves the convergence rate but enlarges  the noise-dominated convergence neighborhood, which can degrade the final accuracy~\cite{wang2026generalized}.}

\subsection{Adaptive gradient algorithms: AdaGrad, Adam}
\label{sec:adaptive_grad_algos}

To overcome the sensitivity of SGD to the choice of the learning rate, several adaptive optimization algorithms have been proposed. These methods adjust the learning rate of each parameter individually based on past gradient information. The most notable examples are AdaGrad and Adam.

\paragraph{Adaptive Gradient(AdaGrad)~\cite{duchi2011adaptive}} AdaGrad updates parameters $\theta$ at iteration $k$ in a component-wise manner, i.e.,
\begin{equation}\label{eq:adagrad}
\theta_{k+1}(i) = \theta_{k}(i) - \frac{\alpha}{\sqrt{G_k(i) + \epsilon}} \, g_k(i),
\end{equation}
where the symbol $(i)$ denotes the i-th component of a vector.
The symbol  $\alpha$ denotes the global learning rate, and $\epsilon$ is a small constant to avoid division by zero. 
The vector $g_k$ stands for the gradient approximation at iteration $k$ and $G_k$ contains the sum of squares of past gradients for each parameter, i.e.,~${G_k(i)=\sum_{j=1}^k g_j(i)^2}$.
\newtext{All together, this gives rise to a preconditioned gradient method with diagonal $M_k$, where $(M_k)(i,i)=\sqrt{G_k(i)}+\epsilon$}.

AdaGrad gives larger updates to infrequent parameters and smaller updates to frequent ones, which can be particularly useful for sparse data. 
However, the accumulated squared gradients in $G_k$ increase indefinitely and can cause the effective learning rate to decrease excessively over time.

\paragraph{Adaptive Moment Estimation (Adam)~\cite{kingma2017adammethodstochasticoptimization}}  
Adam addresses AdaGrad's vanishing learning rate by maintaining exponentially decaying averages of past gradients (first moment) and squared gradients (second moment):
\begin{align*}
m_k(i) &= \beta_1 m_{k-1}(i) + (1-\beta_1) g_k(i), \\
v_k(i) &= \beta_2 v_{k-1}(i) + (1-\beta_2) (g_k(i))^2,
\end{align*}
which are then bias-corrected as
\begin{align*}
\hat{m}_k(i) = \frac{m_k(i)}{1-\beta_1^k}, &&
\hat{v}_k(i) = \frac{v_k(i)}{1-\beta_2^k},
\end{align*}
where $\beta_1$ and $\beta_2$ are decay rates for the moving averages.
\newtext{The bias correction is necessary because $m_k$ and $v_k$ are initialized at zero, causing their early estimates to be biased toward zero. The factors $(1-\beta_1^k)^{-1}$ and $(1-\beta_2^k)^{-1}$ compensate for this effect by rescaling the estimates to the appropriate magnitude; their influence vanishes as $k \to \infty$.}

The moments $\hat{m}_k$ and $\hat{v}_k$ are then used to update the current iterate in the following, component-wise, manner: 
\begin{equation}\label{eq:adam_update}
\theta_{k+1}(i)= \theta_{k}(i) - \alpha \frac{\hat{m}_k(i)}{\sqrt{\hat{v}_k(i)} + \epsilon},
\end{equation}
\newtext{giving a diagonal preconditioned update with $(M_k)(i,i) = (\sqrt{\hat v_k(i)}+\epsilon)/\hat m_k(i)$.}
Adam maintains an exponentially decaying moving average of squared gradients, which remains bounded and thus avoids the vanishing learning-rate issue inherent to AdaGrad’s cumulative accumulation of squared gradients.
As a result,  Adam combines the benefits of AdaGrad's per-parameter learning rates with momentum’s acceleration, making it robust to noisy or sparse gradients.  
It is nowadays implemented in all the ML libraries, and it is one of the most widely used optimizers.

\paragraph{Cost and convergence}
\newtext{Both AdaGrad and Adam have the same per-iteration time complexity as (S)GD, with a small memory overhead. 
AdaGrad stores the accumulated squared-gradient sum $G_k$, while Adam stores the first and second moments $m_k$ and $v_k$.}

\newtext{
Regarding convergence, AdaGrad is fairly well understood. 
For convex Lipschitz functions the rate is $\pazocal{O}(K^{-1/2})$~\cite{duchi2011adaptive,zaheer2018adaptive}, improving to $\pazocal{O}(\log(K)/K)$ for strongly convex functions. 
In the nonconvex setting, AdaGrad achieves $\min_{k\le K}\|\nabla f(\theta_k)\|^2=\pazocal{O}(K^{-1/2})$~\cite{defossez2020simple}, matching the optimal SGD rate under comparable assumptions.

The convergence theory for Adam is more involved. 
Even for convex problems, the original Adam scheme can fail to converge~\cite{reddi2019convergence},  as aggressive decay of the second-moment estimate can discard useful historical gradient information.
Modern variants exist, such as AMSGrad~\cite{reddi2019convergence}, which enjoy convergence guarantees but tend to underperform Adam empirically.
Despite its incomplete theoretical foundation, Adam remains one of the most effective and widely used optimizers in practice.}

\begin{num_example}{\textbf{\emph{SGD, NAG, AdaGrad and Adam.}}}
We now empirically compare the performance of SGD, NAG, AdaGrad, and Adam, on the same logistic-regression problem as considered in Numerical Example~\ref{sec:example_code}.
All methods are run in the full-batch regime ($|\pazocal{I}_k| = 6,000$), ensuring that differences in convergence arise solely from the update rule rather than stochastic sampling effects.
Step sizes are chosen to guarantee stable behavior: SGD and NAG use a small step size $\alpha_k = 10^{-4}$, whereas AdaGrad and Adam employ significantly larger values $\alpha_k = 10^{-2}$.
\newtext{Code snippet~\ref{lst:acc_grad} provides a prototypical implementation of the four optimizers used in this comparison.}

\begin{lstlisting}[language=Python, caption={\newtext{Implementation of SGD, momentum, NAG, AdaGrad, and Adam in NumPy.}}, label=lst:acc_grad]
def run_optimizer( X, y, lam=1e-4, epochs=300, batch_size=16, optimizer="SGD", eta0=1e-3,beta=0.9, beta1=0.9, beta2=0.999, eps=1e-8, seed=0, f_opt=None): 
    # Unified driver that runs SGD / Momentum / NAG / AdaGrad / Adam
    # At each iteration a mini-batch gradient g is computed and an update rule is applied according to `optimizer`
    
    n, dim = X.shape # Dimension of input features and parameters
    rng = np.random.default_rng(seed) # Seed for reproducibility
    w = np.zeros(dim) # Initialization of parameters

    # State variables
    v = np.zeros_like(w)        # Momentum / NAG velocity
    G = np.zeros_like(w)        # AdaGrad sum of squared grads
    m = np.zeros_like(w)        # Adam first moment
    v2 = np.zeros_like(w)       # Adam second moment
    k = 0                       # Global step counter 


    for ep in range(epochs): # Run over epochs
        idx = rng.permutation(n) # Suffle dataset
        Xs, ys = X[idx], y[idx] 
 
        for start in range(0, n, batch_size): # Run over all batches	
            end = min(start + batch_size, n)
            Xb, yb = Xs[start:end], ys[start:end] # Select mini-batch
            k += 1
            
            if(optimizer != "NAG"):	 
            	# Evaluate gradient using current iterate w and mini-batch
                g = grad_logloss(w, Xb, yb, lam) 
 
            if optimizer == "SGD":
                step = -eta0 * g # SGD search direction evaluation

            elif optimizer == "momentum":
                v = beta * v + g # Momentum evaluation
                step = -eta0 * v  # SGD with momentum step

            elif optimizer == "NAG":
                # Computation of look-ahead point
                w_look = w - eta0 * beta * v 
                # Evaluation of gradient at look-ahead point
                g = grad_logloss(w_look, Xb, yb, lam) 
                v = beta * v + g # Momentum computation
                step = -eta0 * v # NAG step

            elif optimizer == "AdaGrad":
                G += g * g # Updating gradient history
                step = -eta0 * g / (np.sqrt(G) + eps) # Adagrad step

            elif optimizer == "Adam":
                m = beta1 * m + (1 - beta1) * g # First-moments
                v2 = beta2 * v2 + (1 - beta2) * (g * g) # Second-moments
                m_hat = m / (1 - beta1 ** k) # Bias corrected first-moments
                v_hat = v2 / (1 - beta2 ** k) # Bias corrected second-moments
                step = -eta0 * m_hat / (np.sqrt(v_hat) + eps) # Adam step

            else:
                raise ValueError("Unknown optimizer")

            w += step # Parameter update with search direction (step)
\end{lstlisting}

Figure~\ref{fig:adaptive_methodp_comparison} shows that standard SGD exhibits the slowest convergence.
NAG and AdaGrad achieve faster progress during the initial phase, but both ultimately remain limited by the ill-conditioning.
Adam, by contrast, exhibits the fastest overall convergence, demonstrating that adaptive preconditioning can stabilize training while enabling the use of larger constant learning rates.

We emphasize that optimal step size choices and the relative performance of these optimizers are problem-dependent and must be tuned in practice.
Moreover, in deep-learning applications, methods that progress slowly or conservatively (such as SGD) may explore the landscape more thoroughly and often yield solutions with better generalization performance.
As a result, identifying the most suitable optimization method and tuning it appropriately remains a key practical challenge in ML workflows.

\begin{figure}
\centering
  \includegraphics{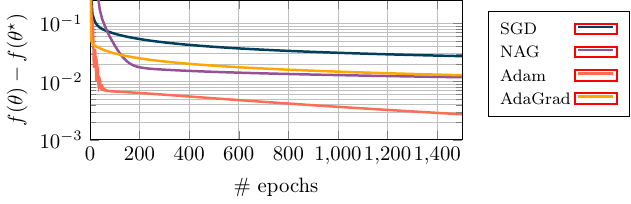}
\caption{Convergence of SGD, NAG, AdaGrad, and Adam on the logistic regression problem of Numerical Example~\ref{sec:example_code}.}
\label{fig:adaptive_methodp_comparison}
\end{figure}
\end{num_example}

\begin{takeawaybox}
\newtext{
Momentum techniques are designed to improve the convergence rate of (S)GD by incorporating past information in the update. 
Adaptive gradient methods can be thought of as preconditioned gradient methods, i.e., their updates are obtained by the gradient one, premultiplied by an invertible matrix intended to approximate the inverse of the Hessian, to include curvature information and speed up convergence. 
In adaptive gradient methods, the preconditioner (the invertible matrix) is diagonal and built from only first-order information (first- and second-order moments of the gradients).
}
\end{takeawaybox}

 \section{Beyond first-order optimization}
\label{sec:second_order}
While first-order methods such as SGD and Adam dominate modern ML, they often struggle in SciML settings, motivating the use of curvature-aware approaches.
Second-order methods, such as Newton's method, can converge much faster near a local minimizer, which is particularly beneficial in SciML, where the loss landscape is often stiff and highly anisotropic due to underlying differential operators. 
\newtext{Coupled with \emph{globalization strategies}, such as line search~\cite{nocedal2006numerical} or trust-region methods~\cite{conn2000trust}, these methods are ensured to convergence from arbitrary starting points}. 
However, their high computational cost makes them impractical at scale, leading to the development of scalable alternatives such as Hessian-free inexact Newton,  subsampled Newton, or quasi-Newton methods.

\subsection{Hessian-free inexact Newton--Krylov methods} 
Recall that the Newton step~\eqref{eq:newton} requires solving the linear system
\begin{align}
    H(\theta_k)\,p_k = -\nabla f(\theta_k)
    \label{eq:search_dir_newton}
\end{align}
for the search direction $p_k$, which is computationally challenging due to both the memory and computational cost of forming and inverting $H(\theta_k)$.

Hessian-free inexact Newton--Krylov methods~\cite{knoll2004jacobian} provide a scalable means of incorporating curvature information without explicitly forming or storing the Hessian.
Rather than solving~\eqref{eq:search_dir_newton} exactly, these methods compute an approximate Newton step 
\newtext{satisfying
\[
    \|H(\theta_k)p_k + \nabla f(\theta_k)\|
    \leq \eta\,\|\nabla f(\theta_k)\|,
\]
with $\eta>0$, using an iterative Krylov solver (e.g., conjugate gradient), which requires only matrix-vector products $H(\theta_k)v$. 
If the tolerance $\eta$ is progressively tightened as the iterations proceed, inexact Newton methods achieve superlinear convergence~\cite{dembo1982inexact}.
A concrete example is the \emph{Eisenstat--Walker} choice of $\eta$, given as~\cite{eisenstat1996choosing}
\[
    \eta_k = \frac{|\,\|\nabla f(\theta_k)\| -
    \|\nabla f(\theta_{k-1}) + H(\theta_{k-1})p_{k-1}\|\,|}
    {\|\nabla f(\theta_{k-1})\|},
\]
which adapts $\eta_k$ automatically based on how well the previous linear system was solved, avoiding both over-solving and under-solving.}

A key advantage of Newton--Krylov methods is that Krylov solvers require only Hessian-vector products.
This is particularly attractive in ML, where Hessian-vector products can often be computed efficiently using automatic differentiation or adjoint techniques.
Such ideas have recently been employed for regression problems in \newtext{SciML}~\cite{zampini2024petscml}, demonstrating higher model accuracy or lower computational cost than adaptive first-order methods.

\paragraph{Cost and convergence}
\newtext{
Each Newton--Krylov iteration is more expensive than a GD iteration, as it requires several inner Krylov iterations, each involving two Hessian-vector products of cost comparable to a gradient evaluation.
In practice, only a few inner iterations suffice to achieve the required accuracy on the Newton direction, and the improved convergence behavior typically compensates for the additional per-iteration cost.
The memory footprint remains comparable to that of first-order methods, since the Hessian is never explicitly formed, making Newton--Krylov one of the few practically scalable second-order approaches for large DNNs.

Since fast local convergence is guaranteed only in a neighborhood of a solution, these methods are often used in combination with 
globalization strategies to ensure convergence from arbitrary starting points}.

\subsection{Subsampled Hessian-free inexact Newton methods} 
Analogous to mini-batch gradient methods, the computational cost of Hessian-free inexact Newton--Krylov methods can be further reduced by employing subsampled Hessian~\cite{roosta2019sub}. 
In this case, \newtext{ mini-batches are used both for the gradient and for the Hessian approximations, and} the step is given as 
\[
    H_{\pazocal{I}^H_k}(\theta_k)p_k = -\nabla f_{\pazocal{I}_k}(\theta_k),
\]
where we used index sets $\pazocal{I}_k\subset \{1,\dots,m\}$ and $\pazocal{I}^H_k\subset \{1,\dots,m\}$ to obtain the stochastic gradient and stochastic Hessian estimates as
$$
\nabla f_{\pazocal{I}_k}(\xk)=\frac{1}{\lvert \pazocal{I}_k \rvert}\sum_{i\in \pazocal{I}_k}\nabla f_{i}(\xk),
\qquad \qquad 
H_{\pazocal{I}_k^H}(\xk)=\frac{1}{\lvert \pazocal{I}^H_k \rvert}\sum_{i\in \pazocal{I}^H_k} \nabla^2 {f_{i}}(\xk). 
$$
Here, we assumed that \newtext{the mini-batch} $\pazocal{I}^H_k$ is uncorrelated with \newtext{the mini-batch} $\pazocal{I}_k$ and  $\lvert \pazocal{I}_k^H\rvert < \lvert \pazocal{I}_k\rvert$, since the iteration scheme is more tolerant to noise in the Hessian estimate than it is to noise in the gradient estimate.

If one chooses the subsample size $\lvert \pazocal{I}_k^H\rvert$ small enough, then the cost of each product involving the Hessian approximation can be reduced significantly, thus reducing the cost of each Krylov iteration. 
Algorithm~\ref{alg:subsampled_newton} summarizes the subsampled Newton method. 
The choice of the step size in this context is often done via line-search strategies \cite[ch.3]{nocedal2006numerical}, adaptive techniques that ensure convergence of the method for any starting point. 

Similar to stochastic first-order methods,  subsampled Newton methods are not particularly well-suited for physics-informed problems due to the global coupling induced by the PDE constraints.
However,  they may be well suited for data-driven SciML problems, where the loss decomposes over independent training samples, e.g., operator learning, DNN surrogate models trained on simulation data, or latent-space representations learned from large collections of solution snapshots.
In such cases, the Hessian admits a natural sample-wise decomposition, making subsampled Newton methods a potentially effective yet largely unexplored alternative to (adaptive) first-order optimization.

\begin{algorithm}[H]
\caption{Subsampled Hessian-free Inexact Newton Method}
\label{alg:subsampled_newton}
\begin{algorithmic}[1]
\State \textbf{Given:} A dataset $\{x_i, y_i \}_{i=1}^m$, an objective function $f:\mathbb{R}^n\rightarrow\mathbb{R}$, an initial iterate $\theta_0 \in \mathbb{R}^n$ 
\For{$k = 0, 1, 2, \dots$}
\State Choose $\pazocal{I}_k^H$ and construct $H_{\pazocal{I}_k^H}(\xk)$  \COMMENTmine{Hessian approximation}
    \State Choose $\pazocal{I}_k$ and $\nabla f_{\pazocal{I}_k}(\xk)$ \COMMENTmine{Gradient approximation}
    \State Obtain $p_k$ by approximately solving \COMMENTmine{Inexact  step computation}
    \[
    H_{\pazocal{I}_k^H}(\xk) p_k = -\nabla f_{\pazocal{I}_k}(\xk)
    \]
     
    \State Choose $\alpha_k$ \COMMENTmine{Step size selection}
    \State Set $\theta_{k+1} = \theta_{k} + \alpha_k p_k$ \COMMENTmine{Iterate update}
\EndFor
\end{algorithmic}
\end{algorithm}

\subsection{Quasi-Newton methods}
Quasi-Newton methods use the same update rule as Newton's method, but instead of the exact Hessian $H(\theta_{k})$, they employ an approximation $B_k$.  
The update is then given as
$
\theta_{k+1} = \theta_{k} +\alpha_k p_k,
$
where $ p_k$ is such that $B_k p_k=-\nabla f(\theta_{k})$. 
The most well-known quasi-Newton method is perhaps BFGS ~\cite{nocedal2006numerical,Liu1989_LBFGS}, which allows for directly approximating the inverse of the Hessian $\tilde{B}_k$, so that the iteration update further simplifies by avoiding the need to solve a linear system. 
Thus, the search direction is obtained as 
$$
p_k=-\tilde{B}_k \nabla f(\theta_{k}). 
$$

At iteration $k$, the inverse Hessian approximation $\tilde B_k$ is constructed recursively, starting from an initial approximation $\tilde B_0$, according to
\begin{equation}\label{eq:bfgs_inverse}
\tilde B_{k+1}
=
(I - \varrho_k s_k y_k^T)\,\tilde B_k\,(I - \varrho_k y_k s_k^T)
+
\varrho_k s_k s_k^T,
\qquad
\varrho_k = \frac{1}{y_k^T s_k},
\end{equation}
where
\[
s_k := \theta_{k+1} - \theta_k,
\qquad
y_k := \nabla f(\theta_{k+1}) - \nabla f(\theta_k).
\]
If $\tilde B_k$ is symmetric positive definite and the curvature condition $y_k^T s_k > 0$ holds, then $\tilde B_{k+1}$ is also symmetric positive definite.

\paragraph{L-BFGS}
\begin{algorithm}[H]
\caption{L-BFGS two-loop recursion}
\label{alg:bfgs_twoloop}
\begin{algorithmic}[1]
\State \textbf{Given:} Current gradient $\nabla f(\theta_k)$, pairs $\{s_i,y_i\}_{i=1}^S$, $\varrho_i=1/(y_i^Ts_i)$, initial approximation of the inverse of the Hessian $\tilde{B}_k^0$ 
\State Set $q=\nabla f(\theta_k)$
\For{$i = k-1, k-2,\dots,k-S $}
    \State Compute $\alpha_i=\varrho_is_i^Tq$ 
    \State Update $q=q-\alpha y_i$ 
    \EndFor
    \State Set $r=\tilde{B}_k^0q$
    \For{$i=k-S,k-S+1,\dots,k-1$}
    \State Compute $\beta=\varrho_iy_i^Tr$
    \State Update $r=r+s_i(\alpha_i-\beta)$
\EndFor\\
\Return $r$
\end{algorithmic}
\end{algorithm}

For large-scale problems, storing the full matrix $\tilde B_k$ is infeasible.
This motivates the use of a limited-memory variant, termed L-BFGS, which stores only $S$ most recent secant pairs $\{(s_i, y_i)\}_{i=k-S}^{k-1}$ instead of the full $\tilde B_k$.
These pairs suffice to approximate the application of $\tilde B_k$ to $\nabla f(\theta_k)$ in order to compute the search direction $p_k$, in turn reducing storage requirements from $O(n^2)$ to $O(nS)$, \newtext{where $n$ is the number of variables}.

To derive an efficient procedure for applying $\tilde B_k$ to $\nabla f(\theta_k)$, we rewrite the BFGS inverse update~\eqref{eq:bfgs_inverse} in the compact form as
\begin{equation}\label{eq:bfgs_compact}
\tilde B_{k+1} = V_k^T \tilde B_k V_k + \varrho_k s_k s_k^T,
\qquad
V_k := I - \varrho_k y_k s_k^T .
\end{equation}
As a consequence, the product $\tilde B_k \nabla f(\theta_k)$ can be computed using only inner products and vector additions involving the stored secant pairs $\{(s_i,y_i)\}_{i=k-S}^{k-1}$.
Thus, by repeatedly applying~\eqref{eq:bfgs_compact}, we obtain 
\begin{align*}
\tilde B_k
&= (V_{k-1}^T \cdots V_{k-S}^T)\, \tilde B_k^{0}\, (V_{k-S} \cdots V_{k-1}) \\
&\quad + \varrho_{k-S} (V_{k-1}^T \cdots V_{k-S+1}^T)\, s_{k-S} s_{k-S}^T\, (V_{k-S+1} \cdots V_{k-1}) \\
&\quad + \varrho_{k-S+1} (V_{k-1}^T \cdots V_{k-S+2}^T)\, s_{k-S+1} s_{k-S+1}^T\, (V_{k-S+2} \cdots V_{k-1}) \\
&\quad + \dots + \varrho_{k-1} s_{k-1} s_{k-1}^T .
\end{align*}
This expression leads to a recursive scheme for computing $\tilde B_k \nabla f(\theta_k)$ efficiently, namely the two-loop recursion summarized in Algorithm~\ref{alg:bfgs_twoloop}.
Assuming $\tilde B_k^{0}$ is diagonal, the two-loop recursion requires approximately $4Sn + n$ multiplications.
A common choice for $\tilde B_k^{0}$ is a scaled identity $\tilde B_k^{0} = \tau_k I$, with $\tau_k > 0$.
\newtext{Parameter $\tau_k$ is often chosen as  $\tau_k = \frac{s_{k-1}^{\top} y_{k-1}}{y_{k-1}^{\top} y_{k-1}}$, which matches the local curvature information from the most recent quasi-Newton pair, see~\cite{nocedal1980updating} for details.}

Algorithm~\ref{alg:bfgs} outlines the L-BFGS procedure.
Note that when using L-BFGS, full gradients are typically required, as the method is sensitive to noise in gradient estimates.
While noise-robust variants of BFGS and L-BFGS have been proposed \cite{shi2022noise,sp_bfgs}, this is not a significant limitation for physics-based models, which are usually trained in deterministic settings due to global coupling induced by PDE constraints.
Consequently, L-BFGS has emerged as one of the most widely used optimizers for PINNs, see a comprehensive benchmark of quasi-Newton methods for PINNs provided in~\cite{kiyani2025optimizing}.

\newtext{
\paragraph{Cost and convergence}
Each L-BFGS iteration costs $\pazocal{O}(nS)$, compared to $\pazocal{O}(n)$ for GD, as it applies the two-loop recursion over $S$ stored curvature pairs $(s_i, y_i)\in\mathbb{R}^n\times\mathbb{R}^n$.
The memory footprint is $\pazocal{O}(2nS)$, which for a typical choice of $S=2$--$5$ is comparable to Adam's within a small constant factor, while avoiding the $\pazocal{O}(n^2)$ cost of full BFGS.

Despite modest overhead, L-BFGS typically converges much faster than GD, often superlinearly in practice.
As the method relies on the positive definiteness of the Hessian approximation, a line search satisfying the Wolfe conditions~\cite{nocedal2006numerical} is standard, ensuring both sufficient decrease and the curvature condition required for the two-loop recursion to remain well-defined.
}

\begin{algorithm}[H]
\caption{L-BFGS}
\label{alg:bfgs}
\begin{algorithmic}[1]
\State \textbf{Given:} A dataset $\{x_i, y_i \}_{i=1}^m$, an objective function $f:\mathbb{R}^n\rightarrow\mathbb{R}$, an initial iterate $\theta_0 \in \mathbb{R}^n$, a memory size $S>0$
\For{$k = 0, 1, 2, \dots$}
    \State Choose $\tilde{B}_k^0$  \COMMENTmine{ Initial inverse Hessian  approximation}
    \State Compute $p_k=-\tilde{B}_k\nabla f(\theta_k)$ by using Algorithm  \ref{alg:bfgs_twoloop}\COMMENTmine{Step computation}

    \State Choose $\alpha_k$
    \COMMENTmine{Step size selection}
    \State $\theta_{k+1} = \theta_{k} + \alpha_k p_k$ 
    \COMMENTmine{Iterate update}
    \If{$k>S$} \COMMENTmine{Secant pairs update}
    \State Discard $\{(s_{k-S},y_{k-S})\}$ from the storage 
    \State Compute and save $s_k=\theta_{k+1}-\theta_k$, $y_k=\nabla f(\theta_{k+1})-\nabla f(\theta_k)$
    \EndIf
\EndFor
\end{algorithmic}
\end{algorithm}

 \subsection{Switching from first- to second-order optimization}
Adaptive first-order optimizers such as Adam are highly effective during the early stages of training due to their ability to handle poor gradient scaling, strong anisotropy, and stochastic noise.
However, their adaptivity can bias the optimization trajectory and limit the asymptotic accuracy of the learned solution~\cite{Wilson2017_AdamGeneralization,Keskar2017_LargeBatch,Heusel2017_GANs}.  
In contrast, quasi-Newton methods are often more effective in later stages of training, when the terms in~\eqref{eq:pinn_expected_risk} vary more smoothly along the optimization trajectory and curvature information can be exploited reliably.

This complementarity motivates hybrid training schedules that combine the robustness of first-order adaptive optimizers with the precision and fast local convergence of quasi-Newton methods.  
The most common strategy in the PINN literature is the two-stage \emph{Adam~$\rightarrow$~L-BFGS} scheme~\cite{Raissi2019_PINN,Wang2023_ExpertGuide}.  
Adam provides a robust initialization phase, bringing the parameters~$\theta$ into a favorable region of the loss landscape.
Once the optimization enters a smoother regime, L-BFGS then refines the solution using curvature-aware updates.

Because the optimal switching point is problem-dependent, switching too early leads to unstable curvature estimates, whereas switching too late wastes iterations.  
To address this, several heuristics have been proposed to determine the suitable transition point in practice. 
For example, one might observe that the gradient norm $\|\nabla_\theta f(\theta_k)\|$ is typically noisy during early Adam iterations but stabilizes once the loss landscape becomes smoother.  
Hence, a persistent plateau in $\|\nabla_\theta f(\theta_k)\|$ over several epochs indicates readiness to switch to L-BFGS~\cite{Wang2023_ExpertGuide,McClenny2020_SA}.
Alternatively, one might explore approaches based on PDE residual saturation~\cite{Wang2023_ExpertGuide,McClenny2020_SA} or stabilization of adaptive learning rates~\cite{Wang2023_ExpertGuide,Ruiz2022_PINN_training_dynamics}.

\begin{num_example}{\textbf{\emph{Training of PINNs.}}}
Figure~\ref{fig:opt_switch} illustrates the contrasting optimization behaviors of Adam, L-BFGS, and a hybrid Adam$\to$L-BFGS strategy for a representative PINN example (1D Poisson problem).
\newtext{An implementation of this hybrid strategy is provided in Code Snippet~\ref{lst:adam_BFGS}.}
As we can see, Adam often achieves rapid loss reduction during the early, noisy, and highly anisotropic phase of training. However, its progress may slow down once the optimization enters a more stable regime where curvature effects dominate.
Standard L-BFGS, when started from the same initialization, fails to make progress at a certain point due to unstable curvature estimates.
In contrast, the hybrid Adam$\rightarrow$L-BFGS schedule automatically transitions when the gradient norm stabilizes, enabling L-BFGS to exploit accurate curvature information and achieve substantially faster convergence and improved accuracy.  
\newtext{Here, the switch from Adam to L-BFGS is triggered when the average gradient norm over $200$ consecutive iterations decreases by less than $0.5\%$.}

\begin{lstlisting}[language=Python, caption={\newtext{Example of using Pytorch's Adam and L-BFGS optimizers for training of PINN (1D Poisson). Automatic switch between Adam and L-BFGS is also implemented.}},label=lst:adam_BFGS]
    def train_adam_then_lbfgs_automatic_switch(model, x, f_rhs, u_exact, adam_max_iters=8000, adam_min_iters=500, lbfgs_iters=1000, lr_adam=1e-4, window=200, plateau_tol=0.995, history_size=50, tol_grad=1e-12, tol_change=1e-12):
    # Train PINN with Adam first, then automatically switch to L-BFGS on convergence plateau

    # Use Adam optimizer of Pytorch to train parameters of model
    opt_adam = optim.Adam(model.parameters(), lr=lr_adam)

    def grad_norm():
        # Compute L2 norm of gradients
        return sum(p.grad.norm().item()**2 for p in model.parameters()
                   if p.grad is not None)**0.5

    # ----- Adam stage: early exploration -----
    for k in range(adam_max_iters):
        opt_adam.zero_grad() # Zero gradient storage in Adam
        # Eval derivative required to evaluate PDE
        uxx = pinn_second_derivative(model, x) 
        loss = loss_fn(uxx, f_rhs) # Compute PINN loss
        loss.backward() # Perform back-propagation to get gradient
        
        # Adam optimizer step
        opt_adam.step()

        L, E = compute_loss_and_error(model, x, f_rhs, u_exact)
        grad_norms.append(grad_norm()) # Store gradient norms

        # Check for convergence plateau using sliding window
        if k > adam_min_iters + window:
            r = np.mean(grad_norms[-window:]) / np.mean(grad_norms[-2*window:-window])
            if r > plateau_tol:  # Stagnation detected
                print(f">>> Switching to L-BFGS at iter {k}")
                break

    # ----- L-BFGS stage: final refinement -----
    # Use L-BFGS optimizer of Pytorch to train parameters of model
    opt_lbfgs = optim.LBFGS( model.parameters(), max_iter=1, history_size=history_size, line_search_fn="strong_wolfe", tolerance_grad=tol_grad, tolerance_change=tol_change)

    for k in range(lbfgs_iters):
        # L-BFGS requires a closure, which return loss, in order to perform line-search to determine optimal step-size
        def closure():        
            opt_lbfgs.zero_grad() # Zero gradient storage in L-BFGS
            # Eval derivative required to evaluate PDE
            uxx = pinn_second_derivative(model, x) 
            loss = loss_fn(uxx, f_rhs) # Compute PINN loss 
            loss.backward() # Perform back-propagation to get gradient
            return loss # Return loss

        opt_lbfgs.step(closure)  # Perform L-BFGs optimizer step
        
\end{lstlisting}

\begin{figure}[t]
  \centering
  \includegraphics{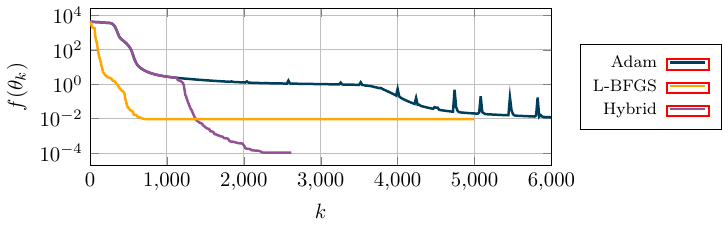}
\caption{PINN training for Poisson problem in 1D using Adam, L-BFGS, and a hybrid Adam$\rightarrow$L-BFGS strategy. The hybrid strategy initially follows Adam, then automatically switches to L-BFGS once the gradient norm stabilizes, yielding significantly faster convergence in the later stages.}
\label{fig:opt_switch}
\end{figure}
\end{num_example}

\begin{takeawaybox}
\newtext{
Second-order methods exploit curvature information contained in the Hessian, enabling faster convergence than first-order methods near a solution.
In their classical form, however, they are too expensive for large-scale problems due to the $\pazocal{O}(n^2)$ cost of forming and storing the Hessian.
Scalable alternatives, such as Hessian-free Newton--Krylov, subsampled Newton, and quasi-Newton (L-BFGS) methods, recover much of this benefit at a cost that remains roughly linear in the number of parameters.
For PINNs, a practical and widely used strategy is to begin training with Adam for robust early-stage progress, then switch to L-BFGS to exploit curvature information and achieve faster final convergence.
}
\end{takeawaybox}

 \section{Soft constraints and loss balancing}
\label{sec:problem_def}
The objective function $f$ is frequently composed of $N$ terms as ${f(\theta) = \sum_{i=1}^N \gamma_i\, \widehat{R}_i(\theta)}$, where $\gamma_i \in \mathbb{R}^+$.
This is due to the fact that typically several terms are incorporated into the loss function in order to improve generalization, enhance the well-posedness of the optimization problem, and enable more stable and efficient training.
Such terms may include regularization or physical constraints, enforced softly through penalization (recall Section~\ref {sec_ERM}).
The following subsections discuss the key ideas behind these formulations as well as the strategies for balancing heterogeneous loss terms.
\newtext{From the preconditioning viewpoint introduced in Section~\ref{sec:precond_lens}, the discussed loss-reweighting strategies act as data-space preconditioners that reshape the empirical loss and the effective Hessian/NTK spectrum.}

\subsection{Loss functions with multiple terms}
\subsubsection{Regularization}
Regularization terms are often used to ensure that the learned model generalizes well beyond the training samples/collocation points.
They augment the objective with an additional term $\widehat{R}_{\text{reg}}$, leading to the modified objective
\[
f(\theta) := \widehat{R}_m(\theta) + \gamma \widehat{R}_{\text{reg}}(\theta),
\]
where $\gamma > 0$ is a regularization parameter. 
A common strategy is to incorporate $L^2$ or $L^1$  penalties on the parameters or on the network's output gradients, e.g., 
\[
\widehat{R}_{\text{reg}}(\theta)
  = 
\|\theta\|^2,
  \qquad \qquad
\widehat{R}_{\text{reg}}(\theta)
  = 
\|\nabla_{\theta} h_{\theta}\|^2.
\]
The first term encourages smoother parameter trajectories, while the second directly penalizes sharp variations in the network output.

In PINNs, one can additionally employ physics-guided regularization to control high-frequency components.  
A frequently used formulation is \newtext{\cite{czarnecki2017sobolev,son2021sobolev}}
\[
\widehat{R}_{\text{reg}}(\theta)
  = 
\int_{\Omega} 
       \big|\nabla^{p} h_{\theta}(x)\big|^2 \, \mathrm{d}x,
\]
where 
$p$ is the differential order of the operator $\pazocal{N}$.  
Such derivative- or spectrum-based penalties enforce physically consistent smoothness and often improve both training stability and out-of-distribution generalization.

\subsubsection{Soft constraints}
\newtext{A canonical example of introducing soft constraints arises in PINNs (cf.~Section~\ref{sec_ERM}), where the PDE residual and the BC conditions are added to the loss via penalty terms, i.e., 
\[
\widehat{R}_m(\theta) =
   \gamma_\Omega\, \widehat{R}_\Omega(\theta)
   + \gamma_D\, \widehat{R}_{\Gamma_D}(\theta)
   + \gamma_N\, \widehat{R}_{\Gamma_N}(\theta)
   + \gamma_{\text{data}}\, \widehat{R}_{\text{data}}(\theta),
\]
where $\widehat{R}_\Omega, \widehat{R}_{\Gamma_D}, \widehat{R}_{\Gamma_N}$, and $\widehat{R}_{\text{data}}$ denote, respectively, the PDE residual, the Dirichlet and Neumann boundary terms, and the data-misfit term as defined in \eqref{eq:risk_PINNs}. 

The weights $\gamma_\Omega, \gamma_D, \gamma_N, \gamma_{\text{data}}\in\mathbb{R}^+$ control how strongly each constraint is enforced. In the limit $\gamma\to\infty$, the penalty formulation formally recovers the corresponding hard-constraint problem.
However, the choice of these weights strongly influences the optimization dynamics. 
If they are too small, the constraint may fail to be enforced effectively. 
If they are too large, the loss landscape becomes stiff, with one term dominating the overall gradient and hindering progress in the remaining components.
The next subsection discusses how to mitigate this difficulty by systematically balancing the contributions of the heterogeneous loss terms.
}

\subsection{Balancing the loss terms}
Selecting appropriate weights $\gamma = (\gamma_1, \ldots, \gamma_N)$ is often challenging. 
Loss-balancing strategies aim to rescale the contributions of each component so that gradients and curvature are more homogeneous.
\newtext{
From the NTK viewpoint of Section~\ref{sec:NTK_illconditioning}, each $\gamma_i$ scales the eigenvalues of the $i$-th NTK block, directly controlling that component's convergence rate.
However, imbalanced weights widen the overall spectrum and degrade the conditioning of the composite problem.
Loss-balancing strategies can therefore be understood as hand-tuned data-space preconditioners (cf.~Section~\ref{sec:precond_lens}).}

\subsubsection{\newtext{Example: Impact of penalty terms on NTK's conditioning}}
\label{sec:example_weights}
\newtext{Let us consider the simplified PINN setting with two loss components, the PDE residual $\widehat{R}_\Omega$ and a boundary term $\widehat{R}_\Gamma$, with corresponding weights $\gamma_\Omega$ and $\gamma_\Gamma$.
In the NTK regime, the empirical NTK acquires a block structure aligned with these two components.
More precisely, neglecting the off-diagonal coupling blocks (which are typically small when the two loss components involve disjoint collocation points), the weighted NTK reads as
\begin{align}
    \Theta^{(\gamma)}
    = \begin{pmatrix}
        \gamma_\Omega\, \Theta_{\Omega\Omega} & 0 \\
        0 & \gamma_\Gamma\, \Theta_{\Gamma\Gamma}
      \end{pmatrix},
\label{eq:block_NTK}           
\end{align}
where $\Theta_{\Omega\Omega}, \Theta_{\Gamma\Gamma}$ denote the per-component empirical NTK blocks.}

\newtext{The condition number of the NTK is then given as
\begin{align}
    \kappa\bigl(\Theta^{(\gamma)}\bigr)
    = \frac{\max\bigl(\gamma_\Omega\,\lambda_{\max}^{(\Omega)},\;
                       \gamma_\Gamma\,\lambda_{\max}^{(\Gamma)}\bigr)}
           {\min\bigl(\gamma_\Omega\,\lambda_{\min}^{(\Omega)},\;
                       \gamma_\Gamma\,\lambda_{\min}^{(\Gamma)}\bigr)},
\label{eq:cond_NTK_weights}                       
\end{align}
where $\lambda_{\max}^{(i)}, \lambda_{\min}^{(i)}$ are the largest and smallest eigenvalues of the NTK blocks. 
As we can see, the NTK spectrum is jointly controlled by the spectra of the individual blocks (problem-dependent) and by the relative scaling of the weights $\gamma_\Omega,\gamma_\Gamma$. 
In particular, even if each block $\Theta_{ii}$ is individually well-conditioned, a large mismatch between $\gamma_\Omega$ and $\gamma_\Gamma$ can severely inflate $\kappa(\Theta^{(\gamma)})$ and slow convergence.
Our overall goal is therefore to explore loss-balancing strategies, which select $\{\gamma_i\}_{i=1}^N$ to reduce $\kappa(\Theta^{(\gamma)})$.}

\subsubsection{Bilevel approaches}
A principled strategy for choosing the weights~$\gamma$ is to treat them as \emph{hyperparameters} that are selected by optimizing some outer performance measure~$\pazocal{C}$, such as the validation loss or the mean PDE residual.  
This leads to the following bilevel optimization problem \newtext{\cite{hao2022bi}}:
\[
\min_{\gamma>0}\; \pazocal{C}(\theta(\gamma)),
\qquad
\text{s.t.}\quad 
\theta(\gamma)
  = \arg\min_{\theta}\, \sum_{i=1}^N \gamma_i\, \widehat{R}_i(\theta),
\]
where the inner problem determines the model parameters $\theta$ for a fixed weights $\gamma$, while the outer problem adjusts $\gamma$ based on the performance measure $\pazocal{C}$.

At the inner optimum~$\theta^\star(\gamma)$, the stationarity condition
\[
\nabla_\theta \widehat{R}_m(\theta^\star,\gamma)=0
\]
implicitly defines the dependence of $\theta^\star(\gamma)$ on $\gamma$.  
Differentiating this condition with respect to~$\gamma$ yields the \emph{hypergradient} $\nabla_\gamma \pazocal{C}$ whose components are given as
\[
\frac{\partial\pazocal{C} }{\partial \gamma_j} 
  = -\,(\nabla_\theta \pazocal{C}(\theta^\star))^\top 
      H^{-1}(\theta^\star,\gamma)\,\nabla_{\theta}^2 \widehat{R}_j(\theta^\star),
\qquad
H = \nabla_\theta \widehat{R}_m(\theta^\star,\gamma),
\]
which shows that the response of~$\theta^\star$ to changes in~$\gamma$ is governed by the
inverse Hessian~$H^{-1}$ of the inner objective.

\newtext{In the NTK regime, the inner Hessian $H$ shares its nonzero eigenvalues with the empirical NTK $\Theta_k$. 
The hypergradient therefore implicitly uses $\Theta_k^{-1}$ to update~$\gamma$, so that bilevel balancing is in fact an \emph{NTK-aware} update of the weights. 
It accounts simultaneously for inter-block scale mismatch and intra-block stiffness in $\Theta^{(\gamma)}$. 
In particular, it drives $\gamma$ toward values that reduce $\kappa(\Theta^{(\gamma)})$. 
However, computing or inverting $H$ is prohibitive in practice, especially for large ML models.
The heuristic strategies discussed below can be viewed as inexpensive, first-order approximations of this exact update, each capturing only part of its effect on the NTK spectrum.}

\subsubsection{Adaptive weighting strategies}
Adaptive weighting strategies are heuristic techniques that can be used to adjust the weights~$\gamma$ during the optimization process so that the different loss components contribute more evenly.

\paragraph{Gradient-norm balancing}
The common approach is to employ \emph{global heuristics}, that assign a single coefficient to each loss component.
One of the most widely used approaches in this category is \emph{Gradient-norm balancing (GN)}~\cite{Chen2018_GradNorm}, which updates the weights~$\gamma$ in order to equalize the magnitudes of the per-loss gradients with respect to the shared parameters.
At each iteration, \emph{GN} computes for each component~$i$ the partial gradients  $ \|\gamma_i \nabla_\theta  \widehat{R}_i(\theta)\|$, measuring the contribution of the $i$-th loss component to the overall gradient.
A target value for these quantities is derived from the relative training speed of the losses:
\[
\rho_{i} =
  \frac{\widehat{R}_i(\theta) / \widehat{R}_i(0)}
       {\tfrac{1}{N}\sum_{j=1}^N \widehat{R}_j(\theta) / \widehat{R}_j(0)}.
\]
\emph{GN} approach is designed such that it minimizes the surrogate objective
\[
L_{\mathrm{grad}}(\gamma)
  = \sum_{i=1}^N \big|\, \|\gamma_i \nabla_\theta  \widehat{R}_i(\theta)\| -  \frac{1}{N}\sum_{j=1}^N \|\gamma_i \nabla_\theta  \widehat{R}_i(\theta)\|\, \rho_i^\zeta \,\big|,
\]
with respect to the weights~$\gamma$, where $\zeta>0$ controls how strongly the method enforces equalized relative descent rates among the losses.
In practice, this step is inexpensive as only a single gradient-descent update is applied to each~$\gamma_i$ using $\partial L_{\mathrm{grad}}/\partial\gamma_i$, followed by a normalization step to maintain weights on a similar scale.
This procedure equalizes the influence of the different loss terms without computing second-order information, thus providing a practical and efficient approximation to bilevel loss balancing.
Other heuristics in this category include for example uncertainty-based weighting~\cite{Kendall2018_Uncertainty, Xiang2021_lbPINN} and stage-dependent curricula learning~\cite{Wang2023_ExpertGuide}.

\newtext{In the two block-NTK example discussed above (Section~\ref{sec:example_weights}), GN rescales the \emph{between-block} contributions. 
By enforcing that $\|\gamma_i\nabla_\theta\widehat{R}_i\|$ is equal across components, GN approximately drives the weights toward
\[
\gamma_\Omega\,\lambda_{\max}^{(\Omega)} \;\approx\; \gamma_\Gamma\,\lambda_{\max}^{(\Gamma)} \;=:\; \mu_{\max},
\]
so that the dominant eigenvalues of the two blocks become comparable. 
Substituting into~\eqref{eq:cond_NTK_weights}, the condition number of the GN-rescaled NTK simplifies to
\[
\kappa\bigl(\Theta^{(\gamma_{\text{GN}})}\bigr)
\;=\; \frac{\mu_{\max}}{\min\!\bigl(\gamma_\Omega\,\lambda_{\min}^{(\Omega)},\;\gamma_\Gamma\,\lambda_{\min}^{(\Gamma)}\bigr)}
\;=\; \max\!\bigl(\kappa^{(\Omega)},\;\kappa^{(\Gamma)}\bigr),
\]
where $\kappa^{(i)} := \lambda_{\max}^{(i)}/\lambda_{\min}^{(i)}$ denotes the within-block condition number. 
GN therefore removes the inter-block scale mismatch in $\Theta^{(\gamma)}$ and brings the overall condition number down to the worst within-block stiffness. However, since each $\Theta_{ii}$ itself is unchanged, GN does not address the intra-component stiffness.}

 \paragraph{Spatially adaptive weighting}
Using a single global weight for each loss component presumes that all samples contribute comparably.
However, in PINNs, the PDE residual is often highly non-uniform, frequently differing by several orders of magnitude across the domain~\cite{McClenny2020_SA}.
This mismatch can cause the optimizer to overfit well-resolved regions while failing to sufficiently reduce the error in stiff or under-resolved parts of the domain.
\emph{Spatially adaptive weighting} strategies were designed to mitigate this difficulty by assigning location-dependent weight~$\xi_j$ to all collocation points, thus replacing the standard residual term by
\[
\widehat{R}_{\Omega}(\theta,\xi)
  = \frac{1}{m_\Omega} 
    \sum_{j=1}^{m_\Omega} 
      \xi_j \, |\pazocal{N}[h_\theta](x_j) - q(x_j)|^2,
\qquad
\sum_{j=1}^{m_\Omega} \xi_j = m_\Omega,\quad \xi_j \ge 0.
\]

A widely used approach for constructing spatially adaptive weights $\xi$ is residual-based weighting~\cite{Anagnostopoulos2023_RBA, Chen2024_BRDR},
which constructs the weights directly from the pointwise residual ${r_{\theta}(x_j) = \pazocal{N}[h_\theta](x_j) - q(x_j)}$.
A typical example includes power-law normalization, where the $j$-th weight takes on the following form:
\begin{align}
\xi_j = 
\frac{|r_{\theta}(x_j)|^\beta}{\tfrac{1}{m_\Omega}\sum_{k=1}^{m_\Omega}|r_{\theta}(x_k)|^\beta},
\qquad \beta>0.
\label{eq:ww_scaling}
\end{align}
This promotes more uniform residual reduction across the domain~$\Omega$ by up-weigting under-resolved regions, i.e., regions with high residual.

\newtext{Recall the two-block-NTK example discussed in Section~\ref{sec:example_weights}. 
The per-point weights $\xi_j$ rescale rows and columns of the residual block. 
In particular, let $D_\xi = \mathrm{diag}(\xi_1,\ldots,\xi_{m_\Omega})$.  
The NTK of residual loss becomes
\[
\widetilde{\Theta}_{\Omega\Omega} \;:=\; D_\xi^{1/2}\, \Theta_{\Omega\Omega}\, D_\xi^{1/2}.
\]
Using residual-based weights~\eqref{eq:ww_scaling}, collocation points with larger PDE residuals receive larger weights. 
Since these points often correspond to directions that are learned slowly, they are associated with small eigenvalues of $\Theta_{\Omega\Omega}$. 
Thus, re-weighting them increases their impact and can make the spectrum of $\widetilde{\Theta}_{\Omega\Omega}$ more balanced.
Thus, the within-block condition number is reduced, i.e., 
\[
\kappa^{(\Omega)}_{\xi} \;:=\; \kappa\bigl(\widetilde{\Theta}_{\Omega\Omega}\bigr)
\;\leq\; \kappa\bigl(\Theta_{\Omega\Omega}\bigr) \;=\; \kappa^{(\Omega)}.
\]
}

\newtext{
As we can see, GN-type and spatially adaptive strategies are complementary. 
The former addresses inter-block imbalance, while the latter tackles intra-block stiffness.
They can therefore be combined, yielding the following condition number:
\[
\kappa\bigl(\Theta^{(\gamma)}\bigr) \;\approx\; \max\!\bigl(\kappa^{(\Omega)}_{\xi},\;\kappa^{(\Gamma)}\bigr) \;<\; \max\!\bigl(\kappa^{(\Omega)},\;\kappa^{(\Gamma)}\bigr).
\]}

\begin{num_example}{\newtext{\textbf{\emph{Convergence of GD with adaptive weights}.}}}\label{sec:weighting_example}
\newtext{We illustrate the analysis of Section~\ref{sec:example_weights} on a small numerical example with two loss components, the PDE residual $\widehat{R}_\Omega$ and a boundary term $\widehat{R}_\Gamma$, with corresponding weights $\gamma_\Omega$ and $\gamma_\Gamma$. 
Following the block-NTK structure of Equation~\eqref{eq:block_NTK}, we assume that the per-component NTK blocks have spectra
\[
\Lambda(\Theta_{\Omega\Omega}) = \text{diag}(0.1,\, 0.04,\, 0.02,\, 0.01),
\qquad
\Lambda(\Theta_{\Gamma\Gamma}) = \text{diag}(10,\, 5),
\]
so that the within-block condition numbers are $\kappa^{(\Omega)} = 10$ and $\kappa^{(\Gamma)} = 2$. 
The boundary block $\Theta_{\Gamma\Gamma}$ thus has eigenvalues two to three orders of magnitude larger than the residual block $\Theta_{\Omega\Omega}$.
This situation commonly occurs while training PINNs, since the boundary functional is often much better conditioned than the differential operator.
}

\newtext{To show the dynamics induced by these spectra, we run GD on the auxiliary quadratic function
\begin{align}
f(\theta) \;=\; \tfrac{1}{2}\,(\theta - \theta^\star)^\top\, \Theta^{(\gamma)}\,(\theta - \theta^\star),
\label{eq:quad}
\end{align}
unique minimizer of which is~$\theta^\star$. 
The Hessian of~\eqref{eq:quad} is exactly the weighted NTK $\Theta^{(\gamma)}$. 
Letting $e_k:= \theta_k - \theta^\star$ denote the iteration error and using the optimal constant step size for quadratic minimization~\cite[Section~3.2]{nocedal2006numerical}
$\alpha \;=\; 2/(\lambda_{\max}\!\bigl(\Theta^{(\gamma)}\bigr) + \lambda_{\min}\!\bigl(\Theta^{(\gamma)}\bigr))$,
the GD iteration $\theta_{k+1} = \theta_k - \alpha\, \nabla f(\theta_k)$ yields the following estimates for the error reduction and convergence rate:
\begin{align}
\|e_k\|/\|e_0\| \;\le\; \rho^k, \qquad \rho \;=\; \frac{\kappa\bigl(\Theta^{(\gamma)}\bigr)-1}{\kappa\bigl(\Theta^{(\gamma)}\bigr)+1}.
\label{eq:GD_estimate}
\end{align}}

\newtext{Using the convergence metrics~\eqref{eq:GD_estimate}, we now compare the expected convergence behavior of GD under three weighting strategies: uniform weights~(U), GN balancing, and a combined approach pairing GN with spatially adaptive weighting~(GN+SA).}

\newtext{\paragraph{Uniform (U) weights ($\gamma_\Omega = \gamma_\Gamma = 1$):}
The composite spectrum of $\Theta^{(\gamma)}$ is the union of the two block spectra, i.e., 
\[
\Lambda\bigl(\Theta^{(\gamma)}\bigr) = \{10,\, 5,\, 0.1,\, 0.04,\, 0.02,\, 0.01\},
\qquad
\text{and}
\qquad
\kappa\bigl(\Theta^{(\gamma)}\bigr) = \frac{10}{0.01} = 1000.
\]
The boundary block dominates while the residual block is essentially frozen. 
From  \eqref{eq:GD_estimate}, the contraction factor is $\rho \approx 0.998$, GD thus requires approximately $1,000$ iterations to reduce the error by a factor of ten.
}

\newtext{\paragraph{GN balancing ($\gamma_\Omega = 1,\ \gamma_\Gamma = 0.01$):}
GN drives the weights toward $\gamma_\Omega \lambda_{\max}^{(\Omega)} \approx \gamma_\Gamma \lambda_{\max}^{(\Gamma)}$, i.e., $\gamma_\Gamma \approx 0.01$. 
The composite spectrum then becomes
\[
\Lambda\bigl(\Theta^{(\gamma_{\text{GN}})}\bigr) = \{0.1,\, 0.1,\, 0.05,\, 0.04,\, 0.02,\, 0.01\},
\quad
\text{and}
\quad
\kappa\bigl(\Theta^{(\gamma_{\text{GN}})}\bigr) = \frac{0.1}{0.01} = 10.
\]
Hence, the condition number is reduced by two orders of magnitude and matches exactly $\max(\kappa^{(\Omega)},\kappa^{(\Gamma)}) = 10$, confirming the prediction that GN collapses the overall conditioning to the worst within-block stiffness. The contraction factor improves to $\rho \approx 9/11 \approx 0.818$, so that only about $12$ GD iterations suffice to reduce the error by an order of magnitude.
}

\newtext{\paragraph{GN with spatially adaptive weighting (GN+SA):}
On top of the inter-block balance, considered in the previous paragraph, we now apply spatially adaptive weights within the residual block.
After such a rescaling, the residual block spectrum is compressed, e.g.,\ from $\{0.1,\, 0.04,\, 0.02,\, 0.01\}$ to a more uniform $\{0.1,\, 0.06,\, 0.05,\, 0.04\}$, so that the within-block condition number drops to $\kappa(\widetilde{\Theta}_{\Omega\Omega}) \approx 2.5$. The overall condition number then satisfies
\[
\kappa\bigl(\Theta^{(\gamma)}\bigr) \;\approx\; \max\bigl(\kappa(\widetilde{\Theta}_{\Omega\Omega}),\, \kappa^{(\Gamma)}\bigr) \;\approx\; 2.5.
\]
As a consequence, only around $3$ iterations are needed to obtain the same tenfold error reduction.
}

\newtext{Figure~\ref{fig:NTK_balancing} visualizes the same experiment. The left panel shows the GD error $\|e_k\|/\|e_0\|$ on the auxiliary quadratic~\eqref{eq:quad}, for the three strategies. 
The empirical curves match theoretical rates closely, and the order-of-magnitude reduction in iteration count predicted by the analysis is clearly visible. The right panel displays the composite NTK spectrum for each strategy.
As we can see, GN collapses the inter-block gap, while SA weighting additionally compresses the residual block toward a more uniform spectrum.
}

\begin{figure}[h!]
\centering
\includegraphics{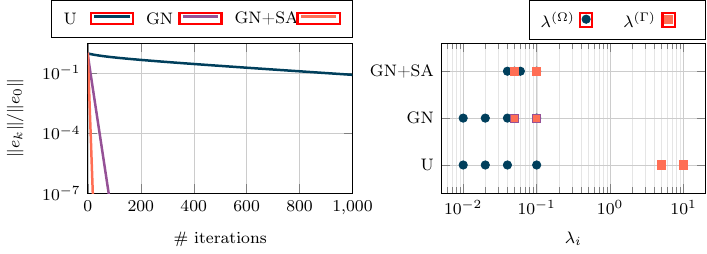}
\caption{\newtext{Effect of  U, GN, and GN+SA strategies on the convergence of GD.
\emph{Left:} Relative GD error $\|e_k\|/\|e_0\|$, with $e_k := \theta_k - \theta^\star$, as a function of the iteration number.
\emph{Right:} Eigenvalues of $\Theta^{(\gamma)}$ at initialization, with residual-block eigenvalues marked by circles and boundary-block eigenvalues by squares.}}
\label{fig:NTK_balancing}
\end{figure}
\end{num_example}

\begin{num_example}{\textbf{\emph{Loss balancing in practice.}}}
\label{sec:pinn_loss_balancing_example}
\newtext{
Building on the insights from the synthetic quadratic experiment of Numerical Example~\ref{sec:weighting_example}, we now examine how the different weighting strategies perform in practical PINN training.
More precisely, we consider the one-dimensional diffusion equation, given as
\begin{align}
u_t \;=\; u_{xx} + e^{-t}\bigl(\pi^2 - 1\bigr)\sin(\pi x), \quad (x,t)\in(-1,1)\times(0,1),
\label{eq:ex7_pde}
\end{align}
with homogeneous Dirichlet conditions $u(\pm1,t)=0$ and initial condition $u(x,0)=\sin(\pi x)$, whose manufactured exact solution is ${u^\star(x,t) = \sin(\pi x)\, e^{-t}}$.
The PINN is an MLP with $4$ hidden layers of width $32$, and tanh activation functions. 
We utilize $256$ residual collocation points and $80$ initial and $80$ boundary points.
Training is performed using the full-batch Adam with learning rate $10^{-3}$ for $5{,}000$ iterations. 
Boundary and initial conditions are enforced softly, so the loss has the three terms, i.e., 
$\widehat R_m = \gamma_\Omega \widehat R_\Omega + \gamma_{ IC}\widehat R_{IC} + \gamma_{BC}\widehat R_{BC}$.
}

\begin{lstlisting}[language=Python, caption={\newtext{Example of computing weights for loss balancing strategies.}},label=lst:loss_weights]
def comp_losses(model, xt_pde, xt_ic, xt_bc, xi=None):
    # Computes per-component losses (R_Omega, R_IC, R_BC)
    # Inputs xt_pde, xt_ic, xt_bc are collocation points for residual, IC and BC
    # Argument `xi` is a per-point weight vector on the residual collocation set;
    # If None is provided, then the residual term is plain mean-square, i.e., xi_j = 1

    # Evaluate residual loss
    r = residual(model, xt_pde).squeeze()
    if xi is None:
        L_pde = (r ** 2).mean()
    else:
        L_pde = (xi.detach() * r ** 2).mean()
        
    # Evaluate IC loss
    u_ic_pred = model(xt_ic).squeeze()
    u_ic_true = torch.sin(PI * xt_ic[:, 0])
    L_ic = ((u_ic_pred - u_ic_true) ** 2).mean()
    
    # Evaluate BC loss    
    L_bc = (model(xt_bc) ** 2).mean()
    return L_pde, L_ic, L_bc


def gradnorm_targets(model, xt_pde, xt_ic, xt_bc, xi=None):
    # Return the target weights based on GN rule
    params = [p for p in model.parameters() if p.requires_grad]
    Lp, Li, Lb = comp_losses(model, xt_pde, xt_ic, xt_bc, xi=xi)
    gp = grad_norm(Lp, params)
    gi = grad_norm(Li, params)
    gb = grad_norm(Lb, params)
    gmax = max(gp, gi, gb) + 1e-30
    return {'pde': gmax / (gp + 1e-30),
            'ic':  gmax / (gi + 1e-30),
            'bc':  gmax / (gb + 1e-30)}, (gp, gi, gb)


def rba_target(model, xt_pde, beta=1.0):
    # Return the residual-based per-point weights based on size of residual
    r = residual(model, xt_pde).detach().abs().squeeze()
    rb = r.pow(beta)
    return r.shape[0] * rb / (rb.sum() + 1e-30)
\end{lstlisting}

\begin{lstlisting}[language=Python, caption={\newtext{Snippet of training loop with adaptive loss scaling.}},label=lst:adaptive_loss_scaling]
    opt = optim.Adam(model.parameters(), lr=lr)
    gamma = {'pde': 1.0, 'ic': 1.0, 'bc': 1.0}
    xi = torch.ones(xt_pde.shape[0])
    xi_snaps = [(0, xi.detach().cpu().numpy().copy())]

    for k in range(iters):
        # Weight updates every gn_every iterations
        # Update xi using SA/RBA first so the subsequent GradNorm target sees it
        if strategy == 'gn_sa' and (k + 1) xi_new = rba_target(model, xt_pde, beta=rba_beta)
            # Smoothing to deal with oscillatory residuals during the initial phase of training
            xi = (1 - gn_alpha) * xi + gn_alpha * xi_new 
            xi_snaps.append((k + 1, xi.detach().cpu().numpy().copy()))
        if strategy in ('gn', 'gn_sa') and (k + 1) xi_used = xi if strategy == 'gn_sa' else None
            tgt, _ = gradnorm_targets(model, xt_pde, xt_ic, xt_bc, xi=xi_used)
            for kk in gamma:
                gamma[kk] = (1 - gn_alpha) * gamma[kk] + gn_alpha * tgt[kk]

        opt.zero_grad()
        xi_used = xi if strategy == 'gn_sa' else None
        Lp, Li, Lb = comp_losses(model, xt_pde, xt_ic, xt_bc, xi=xi_used)
        total = gamma['pde'] * Lp + gamma['ic'] * Li + gamma['bc'] * Lb
        total.backward()
        opt.step()

\end{lstlisting}

\newtext{
We now compare the performance of three weighting strategies discussed earlier in this section: U, GN, and GN+SA ($\xi_j \propto |r_{\theta}(x_j)|^\beta$, $\beta=1$).
An example of their implementation can be found in Code snippets~\ref{lst:loss_weights} and~\ref{lst:adaptive_loss_scaling}.
GN and GN+SA adaptive schemes update weights every $100$ Adam iterations.
For all three strategies, Figure~\ref{fig:PINN_loss_balancing_grads} reports the per-component gradient norms $\|\nabla_\theta\widehat R_j\|,$ where $j=\{\Omega, IC, BC \} $. 
As we can see, with uniform weights, $\|\nabla_\theta\widehat R_{\Omega}\|$ dominates the composite gradient, especially during the second half of the training. 
Utilizing GN weighting, the residual norms have a more comparable scale, which allows for more uniform residual reduction across all iterations. 
The GN+SA approach takes this a step further by addressing both inter-block and intra-block imbalance simultaneously.
In particular, by upweighting collocation points where the residual is currently large, it compresses the intra-block spectrum while sharpening the inter-block balance, allowing for a more substantial reduction of the IC/BC gradients.}

\begin{figure}[h!]
\centering
\includegraphics{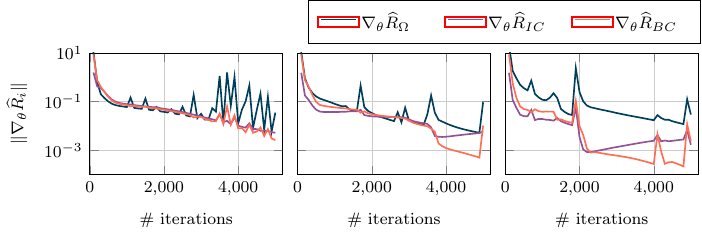}
\caption{\newtext{Per-component gradient norms during Adam training of the diffusion PINN model. 
The results for three different weighting strategies are reported: U, GN, and GN+SA (from left to right).}}
\label{fig:PINN_loss_balancing_grads}
\end{figure}

\newtext{The improvement in the optimization directly translates into a reduction in the test error. 
To observe this, Figure~\ref{fig:PINN_loss_balancing_l2} reports the relative $L^2$ test error
$\|h_{\theta_k} - u^\star\|_{L^2}\,/\,\|u^\star\|_{L^2}$ as a function of the Adam iteration. 
As we can see, using weighting strategies not only allows for faster error reduction but also for achieving lower overall error.
Compared to the quadratic analysis of Numerical  Example~\ref{sec:weighting_example}, which predicted a reduction in the condition number by two orders of magnitude, the practical gains are more modest. 
This is due to the fact that the NTK evolves throughout training,  the block-diagonal structure is only approximate, and Adam's adaptive step sizes partially compensate for gradient imbalance even under uniform weights. 
Nevertheless, the qualitative ranking of the three strategies is preserved, confirming that the NTK conditioning analysis provides a useful predictor of practical training behavior.
}

\begin{figure}[h!]
\centering
\includegraphics{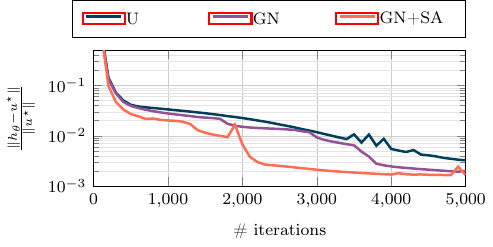}
\caption{\newtext{Relative $L^2$ test error of the diffusion PINN against the manufactured solution $u^\star$ during Adam training with different loss weighting strategies.}}
\label{fig:PINN_loss_balancing_l2}
\end{figure}

\end{num_example}

\begin{takeawaybox}
\newtext{
In SciML, the loss function is often composed of multiple heterogeneous terms that encode different regularizations, physical constraints, and data. 
This is especially prevalent in PINNs, where PDE residuals, initial conditions, boundary conditions, and optional data terms must all be satisfied simultaneously.
The loss terms typically differ in magnitude by several orders of magnitude.
Treating them with equal weights produces a severely ill-conditioned optimization problem in which the gradient of the dominant term monopolizes the rest, so part of the loss is effectively never minimized.
Loss-balancing strategies reweight the terms during training to keep gradient contributions comparably scaled.
The resulting well-balanced loss yields a much better-conditioned empirical NTK, leading to faster convergence.
}
\end{takeawaybox}

 \section{\newtext{Input embeddings and multiscale architectures}}
\label{sec:data_precond}
\sloppy

\newtext{
Recall from Section~\ref{sec:spectral_bias} that standard DNNs often exhibit strong \emph{spectral bias}. 
As a result, the eigenvalues of the empirical NTK $\Theta_k$ can decay rapidly with frequency, so high-frequency components of the residual contract slowly, at a rate $(1-\alpha\lambda_i)^k$ with $\lambda_i \ll \lambda_{\max}$, and may effectively stagnate during training.
}

\newtext{The strategies discussed in this section aim to mitigate spectral bias \emph{before} the optimizer is applied by reshaping the eigenstructure of $\Theta_k$. 
The goal is to flatten the NTK spectrum, or at least to shift the relevant frequency range of the target function into the well-conditioned part of $\Theta_k$.
We illustrate this idea with two concrete examples: Fourier feature (FF) embeddings~\cite{tancik2020fourier} and multiscale network architectures~\cite{wang2021eigenvector}.
From the preconditioning viewpoint adopted in Section~\ref{sec:precond_lens},  FF embeddings act in data space by transforming the inputs before they enter the network, whereas multiscale architectures act in function space by changing the function representation itself.
}

\subsection{\newtext{Fourier feature embeddings}}
\label{ssec:fourier_features}

\newtext{FF embeddings were developed based on the observation that spectral bias originates in the input layer~\cite{tancik2020fourier}.
In particular,  the first layer $\phi_1(x) = \sigma(W_1 x + b_1)$, with random initialization $W_1 \sim \pazocal{G}(0, \sigma_W^2)$, maps low-dimensional inputs to a narrow band of frequencies determined by $\sigma_W$. 
The subsequent layers, therefore, cannot create frequencies that are not already present at the input, and the resulting NTK has rapidly decaying eigenvalues at high frequencies, irrespective of network depth or width.
}

\newtext{The FF embeddings address this limitation by replacing the raw input $x \in \mathbb{R}^d$ with the explicit sinusoidal encoding}
\begin{equation}
  \Phi(x)
  =
  \begin{pmatrix}
    \cos(2\pi B x) \\[2pt]
    \sin(2\pi B x)
  \end{pmatrix}
  \in \mathbb{R}^{2 m_F}.
  \label{eq:fourier_encoding}
\end{equation}
\newtext{Here, $B \in \mathbb{R}^{m_F\times d}$ is a frequency matrix whose rows $b_1,\ldots,b_{m_F}$ are sampled from a fixed distribution, typically $
b_i \sim \pazocal{N}(0,\sigma_B^2 I)$, and $m_F$ is the number of Fourier features. The network then operates on the lifted representation, i.e.,}
\begin{equation}
  h_\theta(x) = \mathrm{DNN}_\theta\bigl(\Phi(x)\bigr).
  \label{eq:ff_network}
\end{equation}
Frequency components up to the scale $\sigma_B$ are therefore injected into the network.
\newtext{Code snippet~\ref{lst:FF} provides an example of a PyTorch implementation of an MLP with Fourier feature embeddings.} 

\begin{figure}[h]
\centering
\resizebox{\linewidth}{!}{\includegraphics{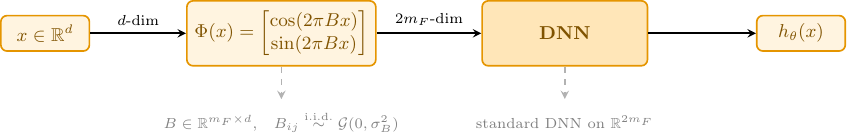}
}
\caption{\newtext{FF architecture. The input $x \in \mathbb{R}^d$ is first mapped to a $2 m_F$-dimensional vector of cosine and sine features by a \emph{fixed} random matrix $B$. A DNN then operates on this lifted representation. The embedding explicitly injects sinusoids up to frequency $\sim \sigma_B$.}}
\label{fig:ff_arch}
\end{figure}

\begin{lstlisting}[language=Python, caption={\newtext{Example of implementation of an MLP with Fourier features.}},label=lst:FF]

class FF_MLP(nn.Module):
    # MLP with Fourier-features
    # Matrix B ~ N(0, sigma_B^2 I) is fixed at initialization
    def __init__(self, sigma_B, m_F=M_F):
        super().__init__()
        # Sample values of B 
        B = torch.randn(m_F, 1) * float(sigma_B) 
        self.register_buffer('B', B)
        # MLP can be any network, but with increased input dimension
        self.core = MLP(in_dim=2 * m_F)
        
    def embed(self, x):
        Bx = 2.0 * float(np.pi) * x @ self.B.t()
        return torch.cat([torch.cos(Bx), torch.sin(Bx)], dim=1)
        
    def forward(self, x):
        # Embedding transforms input before passing it through MLP
        return self.core(self.embed(x)) 

\end{lstlisting}

\subsubsection{\newtext{NTK of a FF network}}

\newtext{FF embeddings transform the NTK into a kernel that can represent frequencies up to the scale $\sigma_B$ without the strong low-frequency bias exhibited by standard NTKs.
To understand this phenomenon, we first show that the NTK depends on $x$ only through $\Phi(x)$. 
In particular, for $h_\theta(x) = \mathrm{DNN}_\theta(\Phi(x))$, the embedding $\Phi$ does not depend on the network parameters $\theta$, so the parameter gradient is just the DNN's parameter gradient evaluated with $\Phi(x)$, and the NTK is the DNN's NTK evaluated at the embedded points:
\[
\Theta_{\mathrm{FF}}(x, x') \;=\; \Theta_{\mathrm{DNN}}\bigl(\Phi(x),\,\Phi(x')\bigr).
\]
A standard result for wide MLPs with random isotropic weights is that  $\Theta_{\mathrm{DNN}}(z,z')$ depends on the two inputs only through their dot product $z^\top z'$ and their norms~\cite{Jacot2018_NTK}. 
For FFs, the norm is fixed, since
\[
\|\Phi(x)\|^2 
= \sum_{f=1}^{m_F}\bigl(\cos^2(2\pi b_f^\top x) 
+ \sin^2(2\pi b_f^\top x)\bigr) = m_F,
\]
so the only quantity that varies is $\Phi(x)^\top \Phi(x')$.
In other words, the FF network's NTK is fully determined by the  inner products $\Phi(x)^\top\Phi(x')$, which we now evaluate explicitly.

Using the relation $\cos a \cos b + \sin a \sin b = \cos(a-b)$, a direct computation gives  $\Phi(x)^\top\Phi(x') = \sum_{f=1}^{m_F}\cos\bigl(2\pi b_f(x-x')\bigr)$. 
With large $m_F$ and $b_f \sim \pazocal{N}(0,\sigma_B^2 I)$, the law of large numbers  replaces the empirical average by its expectation, yielding a Gaussian kernel via the  standard random-features identity~\cite{rahimi2007random}:
\[
\frac{1}{m_F}\,\Phi(x)^\top\Phi(x')
\;\approx\;
\mathbb{E}_b\bigl[\cos(2\pi b(x-x'))\bigr]
\;=\;
\exp\bigl(-2\pi^2\,\sigma_B^2\,\|x-x'\|^2\bigr).
\]

In the infinite-width limit, it follows that the NTK of DNN with FF is a Gaussian kernel
\[
\;\Theta_{\mathrm{FF}}(x, x') \;\propto\; \exp\bigl(-2\pi^2\,\sigma_B^2\,\|x - x'\|^2\bigr),\;
\]
whose bandwidth is controlled directly by $\sigma_B$~\cite{tancik2020fourier}.
The eigenvalues of any translation-invariant kernel are given by the Fourier transform of its profile~\cite{rahimi2007random}, so for $\Theta_{\mathrm{FF}}$  they decay as a Gaussian centered at zero frequency.
More precisely, denoting by $\omega \in \mathbb{R}^d$ the spatial 
frequency variable, the eigenvalues satisfy
\begin{align}
\lambda_\omega^{\mathrm{FF}}
\;\propto\;
\exp\bigl(-\tfrac{\|\omega\|^2}{2\sigma_B^2}\bigr).
\label{eq:eig_NTK}
\end{align}
Thus, inside the band $\|\omega\|\le\sigma_B$, the Gaussian profile is nearly flat.
To see this, we substitute $\|\omega\| = \sigma_B$ in~\eqref{eq:eig_NTK}. 
The ratio of the eigenvalue at the band edge $\|\omega\| = \sigma_B$ to the eigenvalue at zero frequency $\omega = 0$ is therefore
\[
\frac{\lambda_{\sigma_B}^{\mathrm{FF}}}{\lambda_0^{\mathrm{FF}}} 
= \frac{\exp(-\sigma_B^2/2\sigma_B^2)}{\exp(0)} = \exp({-1/2}) \approx 0.6.
\]
As a consequence,  all frequencies within the band $\|\omega\|\le\sigma_B$ have eigenvalues of comparable magnitude and no frequency is favored over another.
Outside the band $\|\omega\|>\sigma_B$, eigenvalues decay rapidly as $\exp(-\|\omega\|^2/2\sigma_B^2)$, so modes far from the band are effectively invisible to the optimizer.

The appropriate choice of bandwidth $\sigma_B$ is crucial in practice. 
If $\sigma_B$ is too small, the DNN with FF encodings behaves almost like a standard DNN, and high-frequency components remain hard to learn. 
If $\sigma_B$ is too large, the model may over-emphasize high frequencies at the expense of large-scale structure. 
Setting $\sigma_B$ to match the highest frequency present in the target ensures that every relevant mode lies in the well-conditioned part of the spectrum. The effective condition number of the NTK on that band drops from $\pazocal{O}(e^{\omega_{\max}})$ to $\pazocal{O}(1)$. 
In this sense, the Fourier feature embedding acts as a  \emph{spectral preconditioner} of the NTK in data space.
In practice, $\sigma_B$ should be chosen on the order of the largest relevant frequency in the target, informed by prior knowledge of the target or a short hyper-parameter sweep.
}

\subsection{\newtext{Multiscale architectures}}
\label{ssec:multiscale_dnn}
\begin{figure}[h]
\centering
\resizebox{\linewidth}{!}{\includegraphics{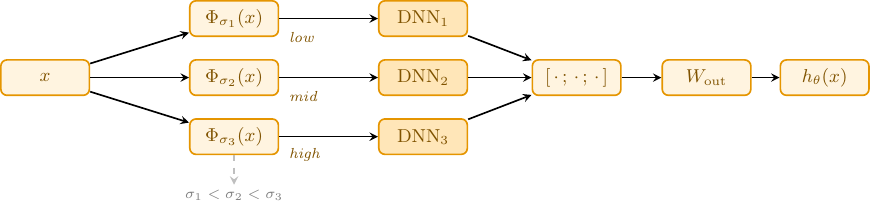}
}
\caption{\newtext{Multiscale architecture with $Q=3$ parallel subnetworks, each applying 
a Fourier feature map $\Phi_{\sigma_q}$ at a different scale  $\sigma_1 < \sigma_2 < \sigma_3$. 
The subnetwork outputs are concatenated and fused by a final linear layer.}}
\label{fig:multiscale_arch}
\end{figure}

\newtext{FF embeddings require knowing $\sigma_B$ in  advance.
Moreover, a single bandwidth often struggles to simultaneously resolve features at very different scales. 
Multiscale architectures address both limitations by employing \emph{parallel subnetworks} operating at  different frequency scales~\cite{wang2021eigenvector}. 
Each subnetwork $q = 1, \ldots, Q$ applies a FF map $\Phi_{\sigma_q}$ at its own scale and feeds the transformed input features through a given subnetwork. 
The subnetwork outputs are then linearly combined as
}
\begin{equation}
\begin{aligned}
  h_\theta(x)
  &= W_{\mathrm{out}}\,\bigl[y_1(x);\,y_2(x);\,\cdots;\,y_Q(x)\bigr] + b_{\mathrm{out}}, \\
  y_q(x) &= \mathrm{DNN}_q\!\bigl(\Phi_{\sigma_q}(x)\bigr), \quad q = 1, \ldots, Q,
\end{aligned}
\label{eq:multiscale_dnn}
\end{equation}
\newtext{see Figure \ref{fig:multiscale_arch} for an example of such architecture. The frequency scales satisfy  $\sigma_1 < \sigma_2 < \cdots < \sigma_Q$ and  are chosen to cover the frequency range of the target solution. 
Each subnetwork acts as a frequency specialist, i.e., $q$-th subnetwork captures features with frequencies of order $\sigma_q$, and the final linear layer combines outputs of all subnetworks. 
The advantage over a single Fourier-feature network is that each subnetwork operates in a well-conditioned regime. Its $\sigma_q$ matches its target frequency band, thereby avoiding tension between low and high frequencies within a single network. 
}

\newtext{Code snippet~\ref{lst:multiscale} provides an example of a  PyTorch implementation of the multiscale architecture~\eqref{eq:multiscale_dnn}  with $Q$ parallel subnetworks.}

\begin{lstlisting}[language=Python, caption={\newtext{Example of multiscale architecture with Q parallel subnetworks.}},label=lst:multiscale]
class MultiscaleFF(nn.Module):
    # Q parallel FF subnetworks with bandwidths `sigmas`, combined by a linear layer
    # Each subnetwork: FF(sigma_q) + MLP -> R^{width}
    def __init__(self, sigmas, m_F=M_F, width=WIDTH, depth=DEPTH):
        super().__init__()
        self.Q = len(sigmas) # Number of subnetworks
        self.subnetworks = nn.ModuleList()
        for s in sigmas: 
            # Generation of a fixed embedding matrix with particular sigma
            B = torch.randn(m_F, 1) * float(s) 
            # Each subnetwork is MLP with a specific FF embeddings 
            subnetwork = nn.Module()
            subnetwork.B = nn.Parameter(B, requires_grad=False)
            layers = [nn.Linear(2 * m_F, width), nn.Tanh()]
            for _ in range(depth - 1):
                layers += [nn.Linear(width, width), nn.Tanh()]
            subnetwork.mlp = nn.Sequential(*layers)
            self.subnetworks.append(subnetwork)

        # The last layer used for concatenation of subnetworks outputs
        self.head = nn.Linear(self.Q * width, 1)
        
    def forward(self, x):
        feats = []
        # Forward propagation is carried out through all subnetworks
        for br in self.subnetworks:
        	    # Each subnetwork is using specific FF encoding
            Bx = 2.0 * float(np.pi) * x @ br.B.t()
            phi = torch.cat([torch.cos(Bx), torch.sin(Bx)], dim=1)
            feats.append(br.mlp(phi))
        cat = torch.cat(feats, dim=1)
        # Head is combining outputs of all subnetworks
        return self.head(cat)
\end{lstlisting}

\subsubsection{\newtext{NTK of multiscale networks}}
\newtext{The preconditioning mechanism behind the multiscale networks architectures becomes transparent in the NTK regime. 
While FFs reshape $\Theta$ into a kernel concentrated at a chosen frequency band, the multiscale architecture factors it into independently conditioned bands.
Both are concrete instances of \emph{spectral preconditioning}, the first acting on the input (data-space) and the second on the function representation (function-space).
Their effectiveness in PINN training has been demonstrated on a wide range of stiff and oscillatory benchmarks~\cite{wang2021eigenvector,tancik2020fourier}.

Because the $Q$ subnetworks share no parameters apart from the output weights $W_{\mathrm{out}}$, the empirical NTK of \eqref{eq:multiscale_dnn} admits an approximate 
block decomposition, with cross-block coupling entering only through  $W_{\mathrm{out}}$, i.e.,
\[
  \Theta_{\mathrm{MS-FF}}(x, x')
  \;\approx\;
  \sum_{q=1}^Q \|w_q\|_2^2\,\Theta_{\mathrm{FF},\,\sigma_q}(x, x'),
\]
where $w_q \in \mathbb{R}^{\mathrm{width}}$ denotes the $q$-th  block of the output weight vector $W_{\mathrm{out}} \in  \mathbb{R}^{Q \cdot \mathrm{width}}$, corresponding to the  contribution of the $q$-th subnetwork to the scalar output.
Moreover,  $\Theta_{\mathrm{FF},\,\sigma_q}$ denotes a single-scale FF kernel at bandwidth $\sigma_q$.
Neglecting the cross-block coupling, this gives the following condition number for $\Theta_{\mathrm{MS-FF}}$:
\[
  \kappa(\Theta_{\mathrm{MS-FF}})
  \;\approx\;
  \max_{1\le q\le Q}\,\kappa\bigl(\Theta_{\mathrm{FF},\,\sigma_q}\bigr).
\]
Hence, $\kappa(\Theta_{\mathrm{MS-FF}})$ is bounded by the worst-conditioned single subnetwork. 
Since each $\sigma_q$ is matched to its target frequency band,  each within-band NTK $\Theta_{\mathrm{FF},\sigma_q}$ is  well-conditioned, i.e., $\kappa(\Theta_{\mathrm{FF},\sigma_q}) 
\approx \pazocal{O}(1)$, and consequently  $\kappa(\Theta_{\mathrm{MS-FF}}) \approx \pazocal{O}(1)$ as well.
Consequently, the overall conditioning is governed by these favorable within-band condition numbers, rather than by the large condition number induced by spanning all frequencies at once. 
Thus, stratifying the network across scales effectively replaces a  single ill-conditioned NTK with a collection of well-conditioned blocks,  one per frequency band.
}

\begin{num_example}{\textbf{\emph{Fourier features and multiscale architectures.}}}
\label{sec:ff_multiscale_example}

\newtext{To make the spectrum-reshaping mechanism of Sections~\ref{ssec:fourier_features}--\ref{ssec:multiscale_dnn} concrete, we revisit the spectral-bias setup of Numerical Example~\ref{num_ex:spectral_bias} with a multi-frequency target, given as
\[
u^\star(x) \;=\; 0.7\,\sin(2\pi x) \;+\; 0.5\,\sin(4\pi x) \;+\; 0.3\,\sin(8\pi x).
\]}

\newtext{We train a small MLP (three hidden layers of width $64$, $\tanh$ 
activations) to fit the target function under different input encodings.
We consider the following architectures: i) a \emph{plain MLP} operating  on raw inputs $x \in \mathbb{R}$;  ii) an MLP with FF  embeddings~\eqref{eq:fourier_encoding} with $m_F = 32$ and bandwidths 
$\sigma_B \in \{2, 4, 8\}$; iii) a multiscale FF network~(MS-FF)  with three parallel subnetworks with the bandwidths $\sigma_B \in \{2, 4, 8\}$.
All networks are trained on $N = 1{,}024$ points sampled uniformly  in $[0,1]$ using Adam for $20{,}000$ iterations.
}

\newtext{
As we can see from Figure~\ref{fig:ff_convergence} on the left, the plain-MLP NTK spectrum decays sharply. 
In MLPs with FF embeddings, FFs inject a fixed band of frequencies and roughly flatten the spectrum within that band.
The bandwidth $\sigma_B$ directly controls which frequencies are emphasized. 
MS-FF concatenates the three bands and, in turn, gives rise to the smallest condition number $\kappa(\Theta)$, see Table~\ref{tab:ff_results}.
}

\newtext{
Figure~\ref{fig:ff_convergence} on the right, and Table~\ref{tab:ff_results} also demonstrate that MS-FF matches the best single-FF in $L^2$ error.  
Its three subnetworks partition the active band, so the result is dominated by within-band conditioning rather than the global frequency range.
Here, we emphasize that the benefit of using MS-FF is not that it always beats the best matched single-FF (it does not, and on this target FF $\sigma_B=4$ is essentially tied with it). 
Rather, it is a safe choice when the target's active band (the range of spatial frequencies that significantly contribute to the solution $u^\star(x)$) is unknown a priori.
Note that every other single-FF in the table is worse than MS-FF on the considered target.
}

\begin{table}[h!]
\centering
\renewcommand{\arraystretch}{1.15}
\begin{tabular}{lcc}
\toprule
\emph{Architecture} & $L^2$ error & $\kappa(\Theta_k)$ at $k=0$ \\
\midrule
MLP                       & $9.58 \times 10^{-4}$ & $1.4\times 10^{11}$ \\
FF, $\sigma_B=2$                & $9.99 \times 10^{-5}$ & $2.5\times 10^{10}$ \\
FF, $\sigma_B=4$                & $\mathbf{5.54 \times 10^{-5}}$ & $1.4\times 10^{10}$ \\
FF, $\sigma_B=8$                & $1.62 \times 10^{-4}$ & $3.0\times 10^{8}$  \\
MS-FF $\sigma_B= \{2,4,8\}$       & $\mathbf{5.59 \times 10^{-5}}$ & $\mathbf{3.0\times 10^{7}}$ \\
\bottomrule
\end{tabular}
\caption{\newtext{$L^2$ error and condition number of the empirical NTK at initialization.}}
\label{tab:ff_results}
\end{table}

\begin{figure}[h!]
\centering
\includegraphics{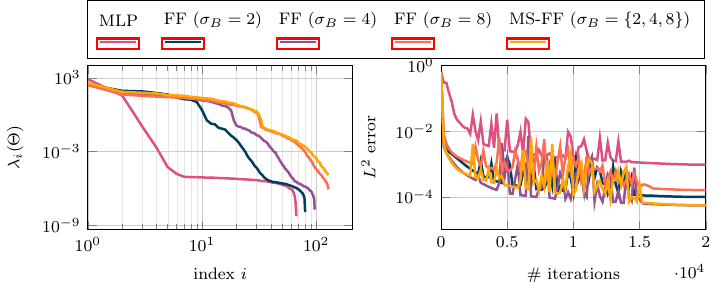}
\caption{\newtext{\emph{Left:} Empirical NTK eigenvalue spectrum at initialization.
\emph{Right:} $L^2$ error of $h_{\theta_k}$ against the target $u^\star$ during training. 
Plain MLP plateaus at $\sim 10^{-3}$, while the FF and MS-FF networks enable convergence to a more accurate solution.}}
\label{fig:ff_convergence}
\end{figure}

\end{num_example}

\begin{takeawaybox}
\newtext{Solving PDEs with sharp features, oscillatory solutions, and  multiscale physics is particularly challenging due to \emph{spectral  bias}: DNNs tend to learn low-frequency content first, while  high-frequency components are learned significantly more slowly. 
This occurs because the eigenvalues of the empirical NTK decay rapidly  with the frequency of the target function. 
Both FF embeddings and multiscale network architectures mitigate this  difficulty by acting as \emph{spectral preconditioners}. 
FF embeddings act as a data-space preconditioner, reshaping the NTK's spectrum by transforming the input data.
Multiscale architectures act as a function-space preconditioner, reshaping the spectrum further by partitioning  the frequency range across independent subnetworks. 
In both cases, the objective is to align the relevant frequencies  with the well-conditioned part of the NTK spectrum, thereby reducing  the effective condition number and accelerating convergence.
}
\end{takeawaybox}

 \section{Sampling and mini-batching strategies for PINNs}
\label{sec:collocation_minibatch}
The approximation quality of PINNs depends crucially on  \emph{how collocation points are sampled} (discretization error)  and \emph{how mini-batches are constructed} during stochastic optimization process (mini-batch sampling error). 
Although largely orthogonal, both reshape the spectrum of the empirical NTK and therefore influence the convergence of the underlying optimizer.
\newtext{From the preconditioning viewpoint of  Section~\ref{sec:precond_lens}, the collocation point sampling strategies discussed in this section act as data-space preconditioners that reshape the effective NTK spectrum.}

\subsection{Choice of collocation points}
As discussed in Section~\ref{sec:PINNs_definition}, PINNs enforce the PDE residual and boundary conditions only at a finite set of collocation points, so the resulting optimization problem is intrinsically determined by the sampling strategy. 
Recall from Section~\ref{sec:NTK_illconditioning} that each collocation point contributes a rank-one term to the empirical NTK, given as
$
\Theta_k 
= \frac{1}{m} \sum_{j=1}^m 
  \nabla_\theta h_{\theta_k}(x_j)\,
  \nabla_{\theta} h_{\theta_k}(x_j)^\top,$
so the distribution of collocation points directly shapes its spectrum.

Undersampling relevant regions (e.g., boundary layers or sharp interfaces) leaves the associated Jacobian directions poorly represented, leading to small or vanishing eigenvalues of the NTK. 
Conversely,  oversampling smoother regions concentrates the spectrum in a few dominant directions. 
Thus, both effects inflate $\kappa(\Theta_k)$, resulting in ill-conditioning.
The sampling distribution, therefore, directly controls the spectral properties of the NTK and, by extension, the difficulty of solving the optimization problem. 
Even with a fixed architecture, convergence can vary dramatically depending solely on the choice of collocation points. 
The design of the collocation set is thus not only a modeling choice, but also a key optimization parameter, particularly in problems with strong spatial heterogeneity.

\subsubsection{Adaptive choice of collocation points}
\label{sec:adaptive_sampling}
A practical way to mitigate NTK/Hessian ill-conditioning caused  by poor sampling is to adaptively relocate or enrich the collocation set in regions identified as problematic.
This approach is closely related to adaptive mesh refinement in classical numerical methods for PDEs~\cite{verfurth1994posteriori}.

In practice, adaptive resampling of collocation points is performed every $T$ iterations, after which training resumes from the current  parameters~$\theta$, creating a feedback loop between optimization and sampling. 
At each refinement step $t+1$, a new collocation set  $\pazocal{D}_\Omega^{(t+1)} = \{ x_j^{(\Omega,t+1)} \}_{j=1}^{m_\Omega}$  is constructed by sampling from the distribution with density
\begin{align}
p^{(t+1)}(x)
= \frac{\bigl(\eta^{(t+1)}(x)\bigr)^{\beta}}
       {\int_{\Omega} \bigl(\eta^{(t+1)}(\tilde x)\bigr)^{\beta}
       \,\mathrm{d}\tilde x},
\label{eq:sampling_density}
\end{align}
where $\eta^{(t+1)}(x)$ is a user-specified indicator function, for example, the pointwise PDE residual  $|r_{\theta}(x)|$, with $r_{\theta}(x)$ defined in \eqref{eq:pde_res}.
The parameter $\beta > 0$  controls the sharpness of the sampling distribution. 
Larger values of~$\beta$ concentrate new points in regions where  $\eta^{(t+1)}$ is large, while smaller values promote broader exploration of the domain.

This adaptive sampling strategy can be interpreted as importance sampling for the residual empirical risk. 
By weighting each collocation point with $\frac{1}{p^{(t+1)}}$, we obtain  an unbiased estimator~\cite{nabian2021efficient} of the continuous residual integral ${\int_{\Omega}   \big|\pazocal{N}[h_\theta](x) - q(x)\big|^2         \, \mathrm{d}\pazocal{P}_\Omega(x)}$, given as 
\begin{align}
\widehat{R}_\Omega(\theta)
= \frac{1}{m_\Omega}
  \sum_{j=1}^{m_\Omega}
  \frac{1}{p^{(t+1)}(x_j^{(\Omega)})}
  \big|\pazocal{N}[h_\theta](x_j^{(\Omega)}) - q(x_j^{(\Omega)})\big|^2.
\label{eq:importance_sampling_estimator}
\end{align}

\paragraph{Conditioning of the NTK}
From a spectral viewpoint, adaptive sampling improves the conditioning of the empirical NTK/Hessian. 
Recall from Section~\ref{sec:NTK_illconditioning} that the NTK is defined as a sum of positive semidefinite rank-one terms, i.e.,~$\Theta^{(m)} = \tfrac{1}{m}\sum_{j=1}^m J_jJ_j^\top$, where $m$ denotes the number of collocation points and $J_j:=J(x_j)$.

Adding a new collocation point yields
\[
\Theta^{(m+1)}
= \frac{m}{m+1}\Theta^{(m)}
 + \frac{1}{m+1}J_{m+1}J_{m+1}^\top,
\qquad J_{m+1}J_{m+1}^\top \succeq 0.
\]
Now, let $A := \frac{m}{m+1} \Theta^{(m)}$, and  $B := \frac{1}{m+1} J_{m+1} J_{m+1}^\top$, so that $\Theta^{(m+1)} = A + B$ with $B \succeq 0$.
By Weyl’s inequality for Hermitian matrices~\cite[Theorem 4.3.1]{horn2012matrix}, it follows that
\[
\lambda_i(A+B) \ge \lambda_i(A),
\qquad i = 1,\dots,n.
\]
Therefore,
\[
\lambda_i(\Theta^{(m+1)})
\ge
\lambda_i\!\left(\frac{m}{m+1} \Theta^{(m)}\right)
=
\frac{m}{m+1}\,\lambda_i(\Theta^{(m)})\geq 0,
\qquad i = 1,\dots,n.
\]
Hence, in practice, when the new Jacobian $J_{m+1}$ introduces a direction not yet represented in  $\{J(x_j)\}_{j=1}^m$, the smallest eigenvalue $\lambda_{\min}(\Theta^{(m)})$ increases sharply, while $\lambda_{\max}(\Theta^{(m)})$ is nearly unaffected.
\newtext{
To understand why, we therefore assume that the new direction is approximately aligned with the eigenvector $q_{\min}$ associated with $\lambda_{\min}(\Theta^{(m)})$,  i.e., with the direction that is most deficient in the current Jacobian span. 
By applying the standard first-order perturbation formula for symmetric  matrices~\cite[Section~6.3]{horn2012matrix}, the shift in each eigenvalue $\lambda_i$ satisfies
\begin{equation}
    \lambda_i(\Theta^{(m+1)}) \approx \frac{m}{m+1}\lambda_i(\Theta^{(m)}) 
    + \frac{1}{m+1}\bigl(q_i^\top J_{m+1}\bigr)^2.
\end{equation}
For $\lambda_{\min}$, the product $q_{\min}^\top J_{m+1}$ is large by assumption, so the correction term dominates and $\lambda_{\min}$  increases sharply. 
For $\lambda_{\max}$, however, the leading eigenvector  $q_{\max}$ corresponds to directions already well-represented in the existing Jacobian. 
Since $J_{m+1}$ was selected for its novelty,  it is nearly orthogonal to $q_{\max}$, i.e., $q_{\max}^\top J_{m+1} \approx 0$. 
The correction term $q_{\max}^\top J_{m+1}$ is therefore negligible, and $\lambda_{\max}$ is increased only by the mild rescaling factor $\frac{m}{m+1}$, leaving $\lambda_{\max}$ nearly  unchanged for large $m$.}

Since the smallest eigenvalue $\lambda_{\min}(\Theta^{(m)})$ increases sharply, while $\lambda_{\max}(\Theta^{(m)})$ is nearly unaffected, the condition number $\kappa(\Theta^{(m)})$ is improved. 
Geometrically, the quadratic form $\theta^\top \Theta^{(m)} \theta$ determines the level sets of the loss landscape. 
Improving~$\kappa(\Theta^{(m)})$ corresponds to making these level sets more isotropic, which in turn permits larger stable step sizes and faster convergence of the GD method.

\subsubsection{Example: Adaptive Sampling for 1D Poisson problem}
\label{sec:toy_example}
To illustrate how adaptive sampling can improve conditioning, we consider the one-dimensional Poisson problem on $\Omega=(0,1)$, given as
${
-\,u''(x) = g(x), \  u(0)=u(1)=0,}
$
with a localized forcing term
\[{
g(x) = 10\, e^{-100 (x-0.9)^2},}
\]
which induces sharp gradients in the solution near $x\approx 0.9$.

We approximate the solution by a simple two-parameter surrogate,
\[
h(x;\theta) = \theta_1\,\phi_1(x) + \theta_2\,\phi_2(x),
\qquad 
\phi_1(x) = x(1-x),\qquad \phi_2(x) = x^2(1-x)^2.
\]
The residual and its Jacobian with respect to the parameters~$\theta$ are given by
\[
r(x;\theta) = -h''(x;\theta) - g(x), 
\qquad 
J(x) = \nabla_\theta r(x;\theta)
      = \begin{bmatrix} -\phi_1''(x) \\[3pt] -\phi_2''(x) \end{bmatrix}.
\]
The NTK  is defined as ${\Theta^{(m)} = \frac{1}{m}\sum_{j=1}^m J(x_j)J(x_j)^\top}$, where $m$ denotes the number of collocation points.

\paragraph{Uniform sampling}
Let the collocation points be uniformly spaced as
\[
x_j \in \{0.1,\,0.3,\,0.5,\,0.7,\,0.9\}.
\]
Most of these points lie in regions where $g(x)$ and $r(x)$ are small, resulting in Jacobians that are nearly collinear.
Consequently, the resulting NTK, ${\Theta^{(5)} = \tfrac{1}{5}\sum_{j=1}^5 J(x_j)J(x_j)^\top}$, has eigenvalues $\lambda_{\max}\approx 4.1$ and $\lambda_{\min}\approx 0.02$, yielding a condition number $\kappa(\Theta^{(5)})\!\approx\!200$.

\begin{figure}[ht]
\centering
\includegraphics{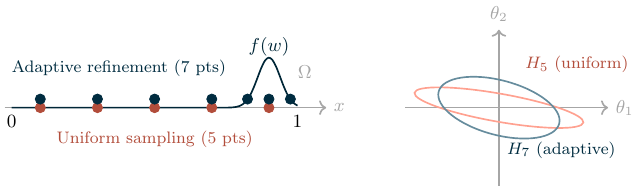}
\caption{
\emph{Left:} Domain $\Omega=(0,1)$ with localized forcing near $x\approx0.9$.
Uniform sampling is depicted in orange, while adaptive refinement is depicted in blue.
\emph{Right:} Level sets of $\theta^\top H_m \theta=1$ with $H_m:=H(x_m)$ become more isotropic after refinement,
reflecting improved conditioning.
}
\label{fig:adaptive_sampling_NTK_poisson_final}
\end{figure}

\paragraph{Adaptive refinement}
Now, let us suppose that a residual-based sampling identifies the high-error region near $x\approx 0.9$
and enriches the collocation set to
\[
x_j \in \{0.1,\,0.3,\,0.5,\,0.7,\,0.85,\,0.9,\,0.95\},
\]
as illustrated in Figure~\ref{fig:adaptive_sampling_NTK_poisson_final}.
The additional points contribute Jacobian directions aligned with the localized feature in~$g(x)$.
The refined NTK,  $\Theta^{(7)} = \tfrac{1}{7}\sum_{j=1}^7 J(x_j)J(x_j)^\top$, exhibits $\lambda_{\max}\approx 4.0$ and $\lambda_{\min}\approx 0.30$,
in turn reducing the condition number to $\kappa(\Theta^{(7)})\!\approx\!13$.
Note that the smallest eigenvalue increases significantly, while the largest remains nearly unchanged.

\subsubsection{Examples of practical refinement strategies}
\newtext{We now briefly review two refinement strategies commonly used in practice~\cite{wu2023comprehensive}. 
Both strategies utilize the indicator function $\eta^{(t)}(x)=|r_{\theta}(x)|$, and are specified as follows:
\begin{itemize}
\item \emph{RAD (Residual-based Adaptive Distribution):} 
Every $T$ iterations, the entire collocation set $\pazocal{D}_\Omega^{(t)}$ is \emph{redrawn} from
\[
p^{(t+1)}(x)\;\propto\; \frac{|r_{\theta}(x)|^{\beta}}{\mathbb{E}_{x\sim\pazocal{U}(\Omega)}\!\left[|r_{\theta}(x)|^{\beta}\right]} \;+\; c,
\]
where $\beta>0$ controls how aggressively the density concentrates on high-residual regions and $c\ge 0$ adds a uniform floor that preserves exploration.
As the collocation points are completely redrawn, the budget of $m_\Omega$ collocation points is held fixed across all refinement steps. 
Note that this sampling strategy is stochastic, and when combined with the importance weights $1/p^{(t+1)}(x_j^{(\Omega)})$ in $\widehat R_\Omega$, it yields an unbiased estimator~\eqref{eq:importance_sampling_estimator}.

\item \emph{RAR-G (Residual-based Adaptive Refinement - Greedy):} 
This approach corresponds to taking the limit $\beta\to\infty$ in~\eqref{eq:sampling_density}.
More precisely, every $T$ iterations, the residual $|r_{\theta}(x)|$ is evaluated on a dense pool of candidates drawn from $\Omega$. 
The $q$ points with the largest $|r_{\theta}(x)|$ are then appended to $\pazocal{D}_\Omega^{(t)}$, yielding $\pazocal{D}_\Omega^{(t+1)} = \pazocal{D}_\Omega^{(t)} \cup \{x_{m_\Omega^{(t)}+1},\ldots,x_{m_\Omega^{(t)}+q}\}$ and $m_\Omega^{(t+1)} = m_\Omega^{(t)} + q$.
In this case, the collocation set grows monotonically, hence the name greedy. 
By the Weyl inequality, each appended $J(x_{m+i}) J(x_{m+i})^\top$ to $\Theta^{(m)}$ can only increase the eigenvalues of $\Theta^{(m)}$, associated with the largest gain along the worst-residual direction.
Since the points are selected greedily rather than by importance sampling, the resulting estimator of $\widehat R_\Omega$ is therefore biased toward high-residual regions.
\end{itemize}
Code snippet~\ref{lst:adaptive_sampling} provides an example of implementing the RAD and RAR-G sampling strategies in Python.
For more advanced indicators and adaptive schemes, we refer to~\cite{wu2023comprehensive, visser2026pacmann}.}

\begin{lstlisting}[language=Python, caption={\newtext{Example of residual-based adaptive sampling strategies (RAR-G and RAD). Every \texttt{resample\_freq} iterations, RAR-G appends one new point to the dataset $\pazocal D_\Omega$ while RAD redraws a whole new dataset.}},label=lst:adaptive_sampling]
def rar_g_select(model, k, pool=10000):
    # Greedy top-k: pick the candidates with the largest |residual|
    cand = sample_uniform(pool) # Sample candidate pool of samples
    r = residual(model, cand).detach().abs().squeeze() # Evaluate residual
    idx = torch.topk(r, k=k).indices 
    return cand[idx] # Return points with top-k largest residual

def rad_select(model, m, pool=10000, k_exp=1.0, c=1.0):
    # Importance sampling: draw m points without replacement from
    # p(x) proportional to |r(x)|^k_exp / mean(|r|^k_exp) + c
    cand = sample_uniform(pool) # Sample candidate pool of samples
    r = residual(model, cand).detach().abs().squeeze() # Evaluate residual
    rk = r.pow(k_exp)
    p = rk / (rk.mean() + 1e-30) + c
    p = p / p.sum() # Construct distribution density
    # Sample m points using distribution $p$
    idx = torch.multinomial(p, m, replacement=False) 
    return cand[idx] # Return points chosen by importance sampling

# Inside the training loop, we use adaptive sampling strategy every resample_every steps
if (k + 1) if   strategy == 'rar': 
		# Expand the existent set of collocation points
		xt = torch.cat([xt, rar_g_select(model, n_add=1)], 0).detach()
	elif strategy == 'rad': 
		# Replace the existent set of collocation points
		xt = rad_select(model, m=N_RAD).detach()
\end{lstlisting}

\begin{num_example}{\newtext{\textbf{\emph{Adaptive sampling on a 1D diffusion PINN.}}}}
\label{sec:pinn_adaptive_sampling_example}
\newtext{The toy example discussed in Section~\ref{sec:toy_example} is intentionally low-dimensional: with only two parameters and a hand-picked feature near $x\!\approx\!0.9$, the conditioning argument is essentially algebraic. 
We now repeat the same comparison using PINN, trained to approximate the solution of the same one-dimensional time-dependent diffusion problem as considered in Numerical Example~\ref{sec:pinn_loss_balancing_example}. 
Thus, we consider the following problem:
\begin{align*}
u_t \;=\; u_{xx} + e^{-t}\bigl(\pi^2 - 1\bigr)\sin(\pi x), \quad (x,t)\in(-1,1)\times(0,1),
\end{align*}
with $u(\pm 1, t) = 0$, $u(x, 0) = \sin(\pi x)$ and 
manufactured exact solution ${u^\star(x,t) = \sin(\pi x)\,e^{-t}}$. 
The setup follows~\cite[Section 3.2]{wu2023comprehensive}: a tanh-MLP with $4$ hidden layers of width $32$, trained by full-batch Adam (learning rate $10^{-3}$), with $80$ initial-condition points and $80$ boundary-condition points held fixed during training.}

\newtext{We compare three strategies for constructing the residual collocation set. 
Uniform sampling uses $m_\Omega = 30$ collocation points drawn once at $t=0$. 
We also consider the RAR-G and RAD strategies described above with an update interval of $1{,}000$ iterations. 
All three runs use the same network initialization and the same initial $30$ interior collocation points.

Figure~\ref{fig:pinn_adaptive_l2} reports the relative $L^2$ test error with respect to the manufactured solution $u^\star$. 
Three approaches behave identically until iteration $1{,}000$, since no refinement has yet been triggered. 
After approximately $2{,}000$ iterations, Adam with uniform sampling plateaus, since with only $30$ fixed collocation points, the residual loss has already been minimized along the directions spanned by the corresponding Jacobians, and Adam can no longer extract additional signal.

In contrast, RAR-G continues to make progress because each refinement step adds a collocation point whose Jacobian is  aligned with the largest-residual direction, expanding the  lower end of the empirical NTK spectrum via a rank-one update.
For RAD, the entire set of collocation points is refreshed every $1{,}000$ iterations, so the empirical NTK is continually enriched with new directions throughout training.}

\begin{figure}[h!]
\centering
\includegraphics{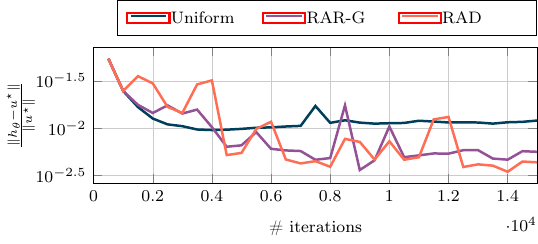}
\caption{\newtext{Relative $L^2$ test error of the PINN solution with respect to the exact solution $u^\star$ during Adam training. Uniform sampling is compared with the adaptive RAR-G and RAD sampling approaches.}}
\label{fig:pinn_adaptive_l2}
\end{figure}
\end{num_example}

\subsection{Mini-batch sampling}
\newtext{Having discussed how the global set of collocation  points should be chosen, we now turn to a finer question: how should  subsets of these points be selected at each iteration of the stochastic optimization process?}
As discussed in Section~\ref{sec:first_order}, SGD relies on stochastic gradient estimates computed using mini-batches.
For classical ML tasks with i.i.d.\ data, these mini-batch gradients provide an unbiased estimator of the full gradient, with variance controlled by the batch size. 
In PINNs, however, selecting informative and stable mini-batches is substantially more delicate, for several reasons:

\begin{enumerate}
\item \emph{Non-i.i.d.\ sampling:} 
When collocation points are generated using structured sampling strategies such as quasi-Monte Carlo or adaptive schemes, the samples are no longer independent. 
As a consequence, classical variance-based analyses of SGD do not apply directly. 
In the non-i.i.d.\ setting, the accuracy of mini-batch gradient estimates is no longer controlled by the variance of the estimator  (as in the classical i.i.d.\ case), but rather by the discrepancy  (a measure of uniformity of the point set over the domain) and the dependence properties of the sampling scheme~\cite{Niederreiter1992_RandomNumber}.

\item \emph{Spatially correlated and heterogeneous sampling:}
The PDE residual is spatially correlated and typically highly nonuniform across the domain.
Mini-batches that contain points in steep-gradient regions generate disproportionately large updates, while batches drawn from smoother regions produce much smaller ones.
This heterogeneity leads to inconsistent optimization steps and can induce oscillations or even divergence.

\item \emph{Highly unbalanced variance:}
For a squared-residual objective, the mini-batch gradient estimator can be written as
\begin{equation}\label{eq:mini_batch}
g_{\pazocal I_k}(\theta)
= \frac{1}{|\pazocal I_k|}
\sum_{x_j \in \pazocal I_k} \nabla_\theta r_\theta(x_j)\, r_\theta(x_j),
\qquad
\phi_\theta(x) := \nabla_\theta r_\theta(x)\, r_\theta(x).
\end{equation}
Under i.i.d.\ Monte Carlo sampling, the variance scales as
\[
\mathrm{Var}\!\left[g_{\pazocal I_k}(\theta)\right]
= \frac{1}{|\pazocal I_k|}\,\mathrm{Var}\!\left[\phi_\theta(x)\right],
\qquad
\mathbb{E}[\|\phi_\theta(x)\|^2]
=
\mathbb{E}\!\left[r_\theta(x)^2 \, \|\nabla_\theta r_\theta(x)\|^2\right].
\]
In PDE-constrained losses, $\nabla_\theta r_\theta$ involves derivatives of order $p$ of the network's output $h_\theta$ and may therefore grow rapidly with the frequency content of the solution, recall Section~\ref{sec:NTK_diff_op}.
As a result, higher-order PDEs or solutions with sharp layers can produce large batch-to-batch fluctuations, destabilizing early training and slowing convergence.

\item \emph{Interaction with NTK conditioning:}
Even if the stochastic gradient variance is moderate, the mini-batch update is still filtered through the NTK.
From~\eqref{eq:lin_dyn}, we know that each update is projected onto the eigenbasis of the NTK.
Directions corresponding to small NTK eigenvalues are thus effectively frozen:
\[
(I - \alpha\Theta_k)\, q_i = (1-\alpha\lambda_i)\, q_i
      \approx q_i
      \qquad\text{for }\lambda_i \ll 1,
\]
so that SGD makes little progress in these components.
However, stochastic gradient noise still projects onto these poorly conditioned subspaces, where curvature is nearly zero, leading to noise accumulation and spurious oscillations.
Since the gradient noise  $\epsilon_k = g_{\pazocal{I}_k} - \nabla\widehat{R}_m$ is isotropic  on average, it projects equally onto all eigendirections of the NTK. 
In well-conditioned directions, the signal dominates the noise.  
In poorly conditioned directions, where the signal is nearly zero,  the noise dominates and causes spurious oscillations.
As a consequence, optimization may appear to stagnate despite rapid progress in better-conditioned directions.

\end{enumerate}

To avoid these difficulties, practitioners often use the \emph{entire} set of collocation points during training, so PINNs are typically trained in a \emph{deterministic} or \emph{very large-batch} regime.
Nevertheless, several works report that mini-batch training can improve generalization in PINNs~\cite{wang2021gradientpathologies}.

\subsubsection{\newtext{Practical mini-batching strategies}}
\label{sec:practical_minibatching}

\newtext{
We now discuss several practical mini-batching strategies for PINNs.
To this aim, let us recall Theorem~\ref{th:nonconvex}. 
Adapted to the mini-batch setting, i.e., taking into account that $\mathrm{Var}\!\left[g_{\pazocal I_k}(\theta)\right]=\sigma^2(\theta)/|\pazocal I_k|$, it states that for a smooth nonconvex objective $\widehat R$ and SGD with \emph{bounded-variance} (Assumption~\ref{hp_noise}) and \emph{unbiased} gradient estimates, we have
\begin{align}
\frac{1}{K}\sum_{k=0}^{K-1}\mathbb{E}[\|\nabla \widehat R(\theta_k)\|^2]
\;\le\;
\frac{2\,\bigl(\widehat R(\theta_0)-\widehat R^\star\bigr)}{\sqrt{K}}
\;+\;
\frac{L_f\,\sigma^2}{\sqrt{K}\,|\pazocal I_k|},
\label{eq:bcn_bound}
\end{align}
with $L_f$ the Lipschitz constant of the gradient of $\widehat{R}$ and  $\widehat R^\star=\inf \widehat R$.
The first term in~\eqref{eq:bcn_bound} is the initialization gap, while the second is the mini-batch variance contribution that scales as $1/|\pazocal I_k|$.

Utilizing standard i.i.d.~mini-batching while training PINNs can break both assumptions of this theorem.
First, even with uniformly randomly sampled collocation points, the stochastic gradient is generally biased (i.e., ${\mathbb E[g_{\pazocal I_k}(\theta)] \ne \nabla \widehat R(\theta_k)}$) if the proportions of sampled points from $\{\Omega, \Gamma_D, \Gamma_N\}$ 
 do not match the corresponding loss weights in $\widehat R = \gamma_{\Omega} \widehat R_\Omega + \gamma_{\Gamma_D} \widehat R_{\Gamma_D} + \gamma_{\Gamma_N} \widehat R_{\Gamma_N}$. 
Second, the PDE residual $r_{\theta}$ involves up-to-$p$-th order derivatives of $h_\theta$, so the term $\nabla_\theta r_{\theta}(x)\,r_{\theta}(x)$ grows rapidly with the frequency content of the solution (see also~\eqref{eq:mini_batch}) and inflates $\sigma^2$ in ways the bounded-variance assumption cannot capture.
}

\newtext{We now discuss some practical mini-batch construction strategies, which can be viewed as tools for 
i) ensuring that the unbiasness assumption of Theorem~\ref{th:nonconvex} is satisfied; 
ii) tightening the initialization gap and 
iii) tightening the variance bound. 
Throughout this discussion, the mini-batch $\pazocal I_k$ is constructed from a given collocation pool $\pazocal D = \pazocal D_\Omega \cup \pazocal D_{\Gamma_D} \cup \pazocal D_{\Gamma_N}$.}

\begin{itemize}
\item \newtext{\emph{Unbiased gradients estimators:} The simplest way to obtain unbiased gradient estimates is to sample the mini-batches according to the loss weights defining 
$\widehat{R}$. 
Thus, given $\widehat R = \sum_S \gamma_S\,\widehat R_S$, $S\in\{\Omega,\Gamma_D,\Gamma_N\}$, we shall sample mini-batches using fixed proportions from $\{\Omega, \Gamma_D, \Gamma_N\}$ matching the weights $\{ \gamma_S \}$.  
Taking expectations over uniform sampling within each subset yields
${\mathbb E[g_{\pazocal I_k}(\theta)] = \sum_S \gamma_S\,\nabla_\theta\widehat R_S(\theta)=\nabla_\theta\widehat R(\theta)}$, so the stochastic gradient estimator is unbiased.

For example, for an equally weighted PINN loss over $\{\Omega, \Gamma_D, \Gamma_N\}$, the mini-batch should contain one third of its samples from each subset.
}

\item \newtext{\emph{Initialization-gap reduction:} A complementary strategy to improve convergence is to target the  initialization gap $\widehat{R}(\theta_0) - \widehat{R}^\star$ 
directly, by decomposing the problem into a sequence of progressively more difficult subproblems. 
This is analogous to continuation and homotopy methods in classical numerical  analysis~\cite{allgower2003introduction}, as well as curriculum learning in classical ML~\cite{bengio2009curriculum}. 
In practice, this approach can be realized, for example, by gradually increasing the complexity of the forcing terms or boundary conditions~\cite{krishnapriyan2021characterizing}.

For time-dependent PDEs, the most natural realization is to respect the causal structure of the PDE, as described below.
The main idea is to fit earlier times before progressively introducing later ones~\cite{wang2024respecting, penwarden2023unified, mattey2022novel}. 
To this end, we restrict the temporal support of the interior/residual component of the mini-batch to
$\{(x,t)\in\pazocal I_k^\Omega \mid t\le T_{\text{cur}}(k)\}$ where $t$ is current time and $T_{\text{cur}}(k)$ increases during training. 
Thus, instead of solving the full problem from a random initialization, we construct a sequence of progressively more difficult subproblems, each initialized from the solution of the previous one. 
As a result, each stage starts from a significantly better parameter estimate, leading to a much smaller effective initialization gap than the naive quantity  $\widehat R(\theta_0^{\text{random}})-\widehat R^\star$, where $\theta_0^{\text{random}}$ is a random initial guess.
}

\item  \newtext{\emph{Variance reduction:} Reducing the variance of the stochastic gradient estimates $g_{\pazocal{I}_k}$ can substantially improve the performance of the underlying optimization algorithms. 
The most straightforward approach is to increase the mini-batch size $|\pazocal I_k|$ as much as computationally feasible, cf.~Numerical Example~\ref{sec:example_code}. 
More sophisticated strategies would include constructing mini-batches via random reshuffling methods~\cite{mishchenko2020random}, sampling with low-discrepancy sequences~\cite{Niederreiter1992_RandomNumber}, or employing importance sampling techniques~\cite{nabian2021efficient}.
 We also note that classical variance-reduction approaches, such as control-variate estimators widely used in Monte Carlo methods and reinforcement learning, have so far received little attention in the PINN literature. 
 Similarly, stochastic variance-reduction methods including SVRG~\cite{johnson2013accelerating} and SAGA~\cite{defazio2014saga} remain largely unexplored in PINN context. 
 Their effective integration into PINN training may therefore provide promising directions for future research.

Among various variance reduction strategies, \emph{stratified sampling} offers a particularly natural fit for PINNs, since the spatial structure  of the domain can be exploited directly to reduce gradient variance.
The idea is to partition the domain into a fixed grid of $C$ cells $\{\Omega_c\}_{c=1}^C$, with weights $w_c = |\Omega_c|/V$, where $V$ is the volume of the full domain, and to enforce that each mini-batch draws a fixed quota of samples from every cell. This removes the between-cell variance that arises from random fluctuations in spatial coverage under i.i.d.\ sampling.

More precisely, let $\phi(x) = \nabla_\theta r_\theta(x)\, r_\theta(x)  \in \mathbb{R}^p$ denote the per-sample gradient contribution, and let  $\sigma^2 = \mathrm{tr}(\mathrm{Cov}(\phi))$ denote its total variance. 
Denoting by $\mu_c = \mathbb{E}[\phi \mid x \in \Omega_c]$ and  $\sigma_c^2 = \mathrm{tr}(\mathrm{Cov}(\phi \mid x \in \Omega_c))$  the mean and total variance of $\phi$ over cell $\Omega_c$,  and by $\mu = \sum_c w_c \mu_c$ the global mean,  the total variance decomposes as
\begin{align}
\sigma^2
=
\underbrace{\sum_{c=1}^C w_c\sigma_c^2}_{\text{within-cell}}
+
\underbrace{\sum_{c=1}^C w_c\|\mu_c-\mu\|^2}_{\text{between-cell}}.
\label{eq:variance_decomp}
\end{align}
An i.i.d.\ mini-batch of size $|\pazocal{I}_k|$ yields variance  $\sigma^2/|\pazocal{I}_k|$, retaining both terms. 
A stratified estimator, by enforcing a fixed quota of samples from  each cell, eliminates the between-cell term  $\sum_c w_c\|\mu_c - \mu\|^2$ by construction, yielding
\[
\mathrm{Var}(g_{\pazocal{I}_k}) 
= \frac{1}{|\pazocal{I}_k|}\sum_{c=1}^C w_c\sigma_c^2,
\]
which is always smaller than or equal to $\sigma^2/|\pazocal{I}_k|$,  with equality only when all cell means $\mu_c$ are identical.

In PINNs, this can be particularly beneficial since $\phi(x)$ is typically highly heterogeneous in space (e.g., near boundary layers or underfitted regions). 
As a result, the between-cell variance is often large, and stratification can significantly reduce gradient variance at no additional computational cost per iteration.
}

\end{itemize}
\newtext{Code snippet~\ref{lst:minibatch} provides an example of implementing  the aforementioned mini-batching strategies in Python.}

\begin{lstlisting}[language=Python,
caption={\newtext{Implementation of mini-batch construction strategies: 
loss-stratified, spatial-stratified, and causal time-marching.}},
label=lst:minibatch]
# Using equal loss-component fractions yield an unbiased gradient estimator
# for the equally-weighted loss R = R_Omega + R_IC + R_BC
STRAT_FRAC = (1/3, 1/3, 1/3)

# Causal schedule: linearly grow the admissible time horizon from T_MIN
# at k=0 to t=1 once k >= T_WARM * K (fraction of total training)
T_MIN, T_WARM = 0.05, 0.3
def t_cur(k, K):
    # Current time horizon for causal time-marching at iteration k
    return T_MIN + (1.0 - T_MIN) * min(1.0, k / max(1, T_WARM * K))

def _split(B, frac):
    # Split batch budget B into (pde, ic, bc) counts
    # BC absorbs rounding so that bs_pde + bs_ic + bs_bc = B exactly
    bs_pde = round(B * frac[0])
    bs_ic  = round(B * frac[1])
    bs_bc  = B - bs_pde - bs_ic
    return bs_pde, bs_ic, bs_bc

def stratified_batch(xt_pde, xt_ic, xt_bc, B, **kwargs):
    # Loss-stratified sampling: fixed (1/3, 1/3, 1/3) split across interior part of the domain, ICs and BCs
    bs_pde, bs_ic, bs_bc = _split(B, STRAT_FRAC)
    return (xt_pde[torch.randperm(len(xt_pde))[:bs_pde]],
            xt_ic [torch.randperm(len(xt_ic)) [:bs_ic ]],
            xt_bc [torch.randperm(len(xt_bc)) [:bs_bc ]])


def sample_spatial(xt_pool, K):
    # Spatial-stratified sampling: partition the domain into a K_x x K_t
    # grid with K = K_x x K_t  cells and draw one point per cell
    # Choose grid dimensions to match the 2:1 aspect ratio of Omega=[-1,1]x[0,1]
    K_t = max(1, int(round((K / 2) ** 0.5)))
    K_x = max(1, int(round(K / K_t)))
    K   = K_x * K_t
    x, t = xt_pool[:, 0], xt_pool[:, 1]
    # Assign each pool point to a cell index in {0, ..., K-1}
    ix    = torch.clamp(((x + 1) / 2 * K_x).long(), 0, K_x - 1)
    it    = torch.clamp((t * K_t).long(), 0, K_t - 1)
    cells = ix * K_t + it
    # For each occupied cell, pick one point uniformly at random
    chosen = []
    for c in range(K):
        ci = (cells == c).nonzero(as_tuple=True)[0]
        if len(ci) > 0:
            chosen.append(ci[torch.randint(len(ci), (1,))].item())
    return torch.tensor(chosen, dtype=torch.long)


def loss_and_spatial_batch(xt_pde, xt_ic, xt_bc, B, **kwargs):
    # Loss and spatial-stratified sampling: loss split + spatial-cell stratification within the interior batch
    bs_pde, bs_ic, bs_bc = _split(B, STRAT_FRAC)
    idx_pde = sample_spatial(xt_pde, bs_pde)
    return (xt_pde[idx_pde],
            xt_ic[torch.randperm(len(xt_ic))[:bs_ic]],
            xt_bc[torch.randperm(len(xt_bc))[:bs_bc]])

def loss_and_causal_batch(xt_pde, xt_ic, xt_bc, B, k, K):
   # Loss-stratified sampling and causal time-marching: restrict interior and BC samples to t <= T_cur(k)
   # Progressively expanding the time horizon
    bs_pde, bs_ic, bs_bc = _split(B, STRAT_FRAC)
    Tk       = t_cur(k, K)
    # Find interior and BC points satisfying the causal constraint t <= Tk
    mask_pde = (xt_pde[:, 1] <= Tk).nonzero(as_tuple=True)[0]
    mask_bc  = (xt_bc [:, 1] <= Tk).nonzero(as_tuple=True)[0]
    # If fewer eligible points are available than needed, sample with replacement
    if len(mask_pde) >= bs_pde:
        idx_pde = mask_pde[torch.randperm(len(mask_pde))[:bs_pde]]
    else:
        idx_pde = mask_pde[torch.randint(len(mask_pde), (bs_pde,))]
    if len(mask_bc) >= bs_bc:
        idx_bc = mask_bc[torch.randperm(len(mask_bc))[:bs_bc]]
    else:
        idx_bc = mask_bc[torch.randint(len(mask_bc), (bs_bc,))]
    idx_ic = torch.randperm(len(xt_ic))[:bs_ic]
    return xt_pde[idx_pde], xt_ic[idx_ic], xt_bc[idx_bc]
\end{lstlisting}

\begin{num_example}{\textbf{\emph{Impact of mini-batch construction strategies on PINNs training.}}}
\label{sec:pinn_minibatch_example}
\newtext{
In this example, we consider an advection-diffusion PDE, given as
$$u_t + c\,u_x - \mu\,u_{xx} = g(x,t), \qquad (x,t) \in (-1,1)\times(0,1),$$
with $c=1$, $\mu=10^{-2}$, so the problem is strongly advection-dominated. 
A manufactured exact solution is $u^\star(x,t) = \exp\!\bigl(-50(x+0.5-ct)^2\bigr)$ that advects from $x=-0.5$ to $x=+0.5$ across the time horizon.
The source $g$ and the non-homogeneous IC and BC are chosen such that $u^\star$ exactly satisfies the PDE.
The PINN model is an  MLP with four hidden layers of width $32$, trained by Adam with a fixed learning rate $\alpha=10^{-3}$ for $15{,}000$ iterations.
The collocation pool $\pazocal D = \pazocal D_\Omega \cup \pazocal D_{\Gamma_D} \cup \pazocal D_{\Gamma_N}$ contains ${m_\Omega = 4{,}096}$ i.i.d.\ uniform interior points and $m_{\text{IC}} = m_{\text{BC}} = 80$ equispaced initial condition and boundary points.
On each iteration, Adam utilizes a mini-batch $\pazocal I_k \subset \pazocal D$ constructed either uniformly or by loss-component stratification.
In addition, we consider variants that incorporate causal time-marching (loss stratification plus interior $t \le T_{\text{cur}}(k)$ growing linearly from $t=0.05$ to $t=1$ over the first $30\%$ of training) and spatial stratification. 
}

\begin{figure}[h!]
\centering
\includegraphics{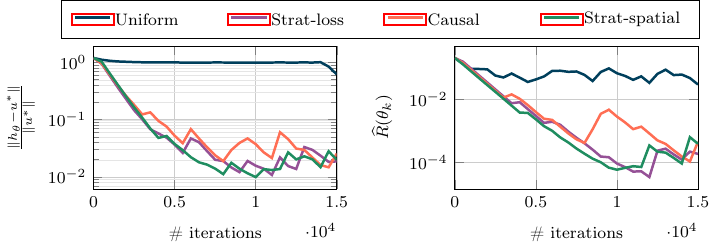}
\caption{\newtext{Training of PINN using Adam with $|\pazocal I_k|= 16$ for an advection-dominated PDE. 
The experiment is performed across four mini-batching strategies: uniform, loss-component stratified (Strat-loss), causal, and spatial stratification (Strat-spatial). \emph{Left:} Relative $L^2$ test error. \emph{Right:} Training loss.}}
\label{fig:ex11_expA}
\end{figure}
\newtext{
Figure~\ref{fig:ex11_expA} reports the relative $L^2$ test error and the training loss trajectories for all considered mini-batching strategies. 
As we can see, employing more advanced mini-batching strategies can improve the $L^2$ test error, compared to using the standard uniform approach.

Figure~\ref{fig:ex11_expB} demonstrates the performance of all considered mini-batching strategies for varying batch sizes.
As shown, increasing the batch size significantly improves the accuracy of the trained PINN. 
We can also see that for very small batches ($|\pazocal I_k|\in \{4, 8, 16\}$), the uniform approach diverges, whereas more advanced strategies enable Adam to converge. 
For small-to-medium batches ($|\pazocal I_k|\in \{16, 32, 64, 128\}$), the performance of the uniform mini-batching strategy becomes more competitive.
In the large-batch regime ($|\pazocal I_k|\in \{256, 512\}$), the uniform approach recovers and achieves comparable accuracy as the advanced mini-batch construction approaches.  
Overall, these results indicate that advanced mini-batching techniques are particularly important under tight computational and memory budgets, i.e., $|\pazocal I_k|\lesssim 64$.
Here, we note that the larger gains in accuracy are expected for stiffer evolution PDEs, such as those considered in~\cite{wang2024respecting, daw2022rethinking}.
}

\begin{figure}[h!]
\centering
\includegraphics{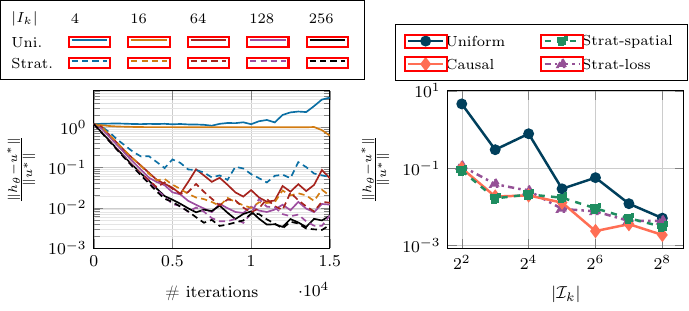}
\caption{\newtext{
Relative test $L^2$-error of PINN for advection-diffusion example with respect to varying batch sizes $|\pazocal I_k|$. 
\emph{Left:} Error observed during Adam's training with uniform and stratified-spatial mini-batching. 
\emph{Right:} Final test error obtained using Adam with varying mini-batching strategies.
}}
\label{fig:ex11_expB}
\end{figure}

\end{num_example}

\begin{takeawaybox}
\newtext{
The choice of sampling and mini-batching strategies significantly influences the effectiveness of PINN training. 
For collocation points, adaptive resampling enriches the collocation set in poorly resolved regions, in turn reducing the condition number of the empirical NTK and improving convergence of the underlying optimizer. 
For mini-batches, loss-stratified sampling can be used to construct unbiased gradient estimates by enforcing fixed proportions of samples from each loss component $\{\Omega, \Gamma_D, \Gamma_N\}$,  while spatial stratification can reduce gradient variance by ensuring uniform spatial coverage. 
In addition,  causality-aware training can reduce the effective initialization gap by decomposing the problem into a sequence of progressively more difficult subproblems.
Altogether, this leads to faster and more reliable convergence of stochastic optimizers applied to PINNs.
}
\end{takeawaybox}

  \section{Towards scalable training for SciML: Recent advances and open challenges}
\label{sec:outlook}
In the preceding sections, we reviewed foundational optimization methods (SGD and its accelerated, adaptive, and second-order variants) together with SciML-specific design choices in loss formulation, model architecture, and sampling strategy. 
In this final section, we survey recent trends and open challenges  in optimization for large-scale SciML training, organized according to whether the preconditioning strategy acts in parameter, data,  or function space (cf.~Section~\ref{sec:precond_lens}).
Our review focuses on how classical scientific computing and ML  techniques transfer to SciML, identifies their current limitations,  and highlights open research directions that remain to be explored.
Given the breadth and rapid evolution of the field, we provide a structured overview of selected research directions rather than a comprehensive survey.

\subsection{Data-space preconditioning}
Data-space preconditioning encompasses all strategies that improve optimization by transforming the inputs, reweighting the loss, or reshaping the sampling distribution, without modifying the network architecture or the optimizer itself.

Classical data-space preconditioning techniques in classical ML that are also widely used by SciML practitioners include feature scaling and standardisation~\cite{zhongkai2024pinnacle}, data augmentation~\cite{maharana2022review}, and equation non-dimensionalisation~\cite{Wang2023_ExpertGuide}.
However, in SciML settings, these ubiquitous techniques often lack a rigorous theoretical foundation.
For instance, in contrast to standard ML, normalization may interact in subtle ways with the PDE structure, or the constraint enforcement, potentially leading to undesirable optimization dynamics~\cite{karniadakis2021physics}.

Other input transformations specifically designed for SciML have been motivated by taking into account physical considerations.  
In particular, input representations that enrich the coordinate system by introducing high-frequency components have been successfully proposed to mitigate the spectral bias of DNNs.  
Prominent examples include Fourier feature mappings~\cite{tancik2020fourier} \newtext{(see Section \ref{ssec:fourier_features})}, positional encodings~\cite{quarteroni2025combining}, or graph embeddings~\cite{makarov2021survey}.  
While these techniques substantially accelerate the learning of oscillatory or multiscale solutions, the choice of frequencies, bandwidths, or embedding structures remains largely problem-dependent and empirical, leaving systematic selection criteria as an open research direction~\cite{rahaman2019spectral,wang2021gradientpathologies}.

In the context of PINNs, \textit{sampling strategies}  also play a central role in determining how residual and BC's collocation points influence the conditioning of the underlying optimization problem (see Section~\ref{sec:collocation_minibatch}).  
Consequently, many advanced sampling strategies can be interpreted as data-space preconditioners.  
Popular examples include variance-reduced and importance sampling strategies~\cite{nabian2021efficient,daw2022rethinking}, adaptive and residual-based sampling~\cite{wu2023comprehensive}, or failure-informed or active sampling approaches~\cite{gao2023failure}.  
Similarly, \textit{adaptive weighting} strategies for balancing loss terms play a key role in stable and efficient PINN training (see Section~\ref{sec:problem_def}).
Existing methods dynamically adjust weights based on gradient magnitudes, NTK spectra, or residual statistics to alleviate stiffness and imbalance during optimization (e.g.,~\cite{mcclenny2023self,subramanian2022adaptive,anagnostopoulos2024residual,Wang2021_whenPINNsFail}).
Despite their practical effectiveness, both sampling and weighting strategies remain largely heuristic, and developing principled,  theoretically grounded alternatives constitutes an open research direction.

\subsection{Parameter-space preconditioning}
Parameter-space preconditioners improve optimization by  modifying how parameter updates are computed, for example, by 
incorporating curvature information to transform the gradient direction, or by restricting updates to lower-dimensional, better-conditioned subspaces of the full parameter space.

Established optimization techniques such as AdaGrad~\cite{duchi2011adaptive}, RMSProp~\cite{tieleman2012lecture}, and Adam~\cite{kingma2017adammethodstochasticoptimization} are computationally efficient and robust in many large-scale ML applications, yet their behavior in SciML settings remains poorly understood.  
In particular, recent studies have reported gradient pathologies and unexpected optimizer dynamics in stiff PDE-constrained problems, suggesting nontrivial interactions between adaptive scaling and PDE-induced spectral properties~\cite{scott2025designing}.

A complementary class of methods attempts to approximate the inverse Hessian or Fisher information matrix.
Algorithms such as natural gradient descent~\cite{amari1998natural}, K-FAC~\cite{martens2015optimizing}, Shampoo~\cite{gupta2018shampoo}, SOAP~\cite{vyas2024soap}, or MUON~\cite{liu2025muon} maintain structured curvature approximations, often exploiting Kronecker or layerwise factorizations to remain tractable.
However, due to their awareness of curvature, these approaches remain prohibitive in large-scale settings owing to substantial memory and computational overhead~\cite{zampini2024petscml}.

Specific to SciML are preconditioners that exploit its connections to numerical linear algebra and scientific computing. 
Several recent approaches adapt large-scale preconditioning techniques from PDE solvers to DNN training, with a particular emphasis on \textit{domain-decomposition and multilevel methods}.  
These methods provide implicit preconditioning by restricting updates to lower-dimensional, better-conditioned parameter subspaces, which is especially appealing in SciML, where parameter groups often align with physical modes or operator components.

In domain-decomposition approaches~\cite{klawonn2024machine}, the parameter space is often partitioned into blocks (e.g., layers, or channels), and optimization proceeds via additive or multiplicative Schwarz-type iterations; see, e.g.,~\cite{kopanivcakova2024enhancing, lee2026two,gu2022decomposition,amid2022locoprop,siegel2023greedy}.  
Closely related multilevel approaches equip the parameter space with a hierarchy of coarse-to-fine DNN representations, analogous to mesh hierarchies in multigrid methods~\cite{trottenberg2000multigrid}.
\newtext{Networks with fewer parameters}, obtained for example by reducing network width, or depth, provide inexpensive approximations of curvature and enable the construction of coarse search directions that can accelerate convergence~\cite{Kopanicakova_2020c,kopanicakova_22_1,chang2017multi,ahamed2025multiscale,ponce2023multilevel,Calandra02012022}.

Despite the promising empirical performance of several domain-decomposition and multilevel methods, the design of effective coarse spaces,  inter-level transfer operators, and principled parameter  decompositions with convergence guarantees in stochastic regimes 
remains an important open problem, with recent progress including multilevel AdaGrad~\cite{Kopanicakova_2023a, gratton2025recursive} and multilevel adaptive regularization techniques \cite{marini2024multilevel}.

\subsection{Function-space preconditioning}
As introduced in Section~\ref{sec:precond_lens},  function-space preconditioners modify the geometry of the 
optimization problem directly in the space of functions represented by the DNN, for example, by changing the inner product or norm in which the loss is measured.

One example of function-space preconditioning is Sobolev training~\cite{czarnecki2017sobolev}, which augments the loss with derivative information, thereby modifying the inner product in function space. 
In the SciML context, variants of this idea have been developed to guide PINNs toward solutions that minimize errors in Sobolev norms, often improving convergence rates and robustness relative to standard $L_2$ losses~\cite{son2021sobolev}.
Nevertheless, the choice of derivative orders, relative weightings, and their interaction with PDE stiffness and BC conditions remains largely experimental, and a systematic theory linking Sobolev loss design to operator spectra and optimization conditioning remains an active area of research.

Beyond loss modification, function-space preconditioning can also be achieved by restructuring the functional basis itself,  drawing again inspiration from domain-decomposition and multilevel methods in classical numerical analysis. 
\newtext{Here, we note that some of the techniques discussed below admit a dual interpretation as they act simultaneously in function space and parameter space, and are discussed here primarily from the function-space viewpoint.}
In particular, in the context of PINNs, several domain-decomposition variants explicitly partition the physical domain and assign separate DNNs to each subregion, as in XPINNs~\cite{jagtap2020extended}, cPINNs~\cite{jagtap2020conservative}, FBPINNs~\cite{moseley2023finite}, or APINNs~\cite{hu2023augmented}.
From a function-space perspective, this implicitly constructs a basis of locally supported functions, yielding a block-structured representation that improves conditioning and stabilizes the optimization process. 
\newtext{From a parameter-space perspective, the domain partitioning induces a block-diagonal structure on the parameter space,  analogous to a Schwarz-type preconditioners. 
It is worth noting that these methods also enable parallelization and facilitate the representation of localized or multiscale features, often making PDE problems that are otherwise inaccessible to monolithic PINNs tractable.}

Similarly, multiscale architectures such as those discussed in Section~\ref{ssec:multiscale_dnn} can be interpreted as function-space preconditioners since by combining subnetworks operating at different frequency scales, they partition the functional basis across frequency bands, reducing spectral bias and improving NTK conditioning. 
More broadly, several recent multilevel strategies construct hierarchies of function spaces to train SciML models more efficiently. 
By first learning a coarse approximation and then refining the function space, these approaches also reduce spectral bias and accelerate convergence.
Representative examples include multilevel training for data-driven models~\cite{lye2021multi} and multifidelity PINNs~\cite{riccietti2022multilevel,meng2020composite}, block-structured coarse-to-fine training~\cite{gratton2024block}, and multilevel domain-decomposition methods~\cite{dolean2024multilevel}, while related multiresolution and wavelet- or Fourier-based architectures~\cite{he2024multi} achieve similar effects by separating low- and high-frequency components.

Despite the empirical success of these approaches, systematic design principles for subdomain partitioning, interface conditions, coarse space selection, and convergence guarantees in stochastic training regimes remain important open problems.

\section{\newtext{Conclusion and} outlook}
\newtext{
In this chapter, we reviewed optimization techniques for training  SciML models, emphasizing that efficient training requires not only 
a suitable optimizer, but also problem-aware design choices in the loss formulation, model architecture, and sampling strategy. 
These two aspects were addressed in complementary parts of the chapter. 
Foundational optimization tools were covered in  Sections~\ref{sec_ERM}--\ref{sec:second_order}, while  problem-aware design choices were discussed in  Sections~\ref{sec:problem_def}--\ref{sec:collocation_minibatch}.}

\newtext{The preconditioning viewpoint introduced in  Section~\ref{sec:precond_lens} provides a unifying framework 
for organizing the existing algorithms and connecting them to ideas from classical scientific computing. 
The large-scale approaches reviewed in Section~\ref{sec:outlook} illustrate how modern SciML training increasingly draws on principles from numerical PDE solvers while retaining the flexibility of DNNs. 
Together, data-space, parameter-space, and function-space preconditioning strategies tackle complementary challenges of the underlying optimization problem, including spectral bias, ill-conditioning, and the multiscale nature of PDE  solutions, offering a unified perspective for improving the efficiency and reliability of large-scale SciML training.}

Nevertheless, many open questions remain, particularly in nonconvex and stochastic settings, and the relationships between these three preconditioning viewpoints are still only partially understood.
As a result, the design of efficient optimization techniques remains an active area of research, with continued progress needed to develop scalable and reliable methods for increasingly complex SciML applications.

\section*{Acknowledgements}
The work of A.K.~benefited from Artificial and Natural Intelligence Toulouse Institute (ANITI), funded by the  France 2030 program under Grant Agreement No.~ANR-23-IACL-0002. 
The work of E.R.~was supported by the ANR project MEPHISTO (ANR-24-CE23-7039), the PEPR IA SHARP and the Fondation Simone et Cino Del Duca.

\bibliographystyle{siam}
\bibliography{biblio}

\end{document}